\documentclass[a4paper]{amsart}
\usepackage{amssymb}
\usepackage{amsmath,amssymb,mathrsfs,latexsym,tikz,cite,hyperref}
\usepackage{amsfonts}
\usepackage{genyoungtabtikz}
\usetikzlibrary{tikzmark}
\usetikzlibrary{decorations.pathreplacing}
\usetikzlibrary{arrows,shapes,calc}
\usepackage{rotating, graphicx}
\usepackage{lscape,marginnote}
\usepackage{mathtools}
\usepackage{comment}

\usepackage{amsaddr}
\usepackage{subcaption}
\usepackage{amssymb}
\usepackage[cmtip,all]{xy}

\usepackage{fancyhdr} 
\usepackage{indentfirst}
\usepackage{newlfont}
\usepackage{hyperref}

\usepackage{amsthm}
\usepackage{bm}
\usepackage{dsfont}
\usepackage{tikz-cd}
\tikzcdset{every label/.append style = {font = \small}}
\usepackage{nicematrix}

\usepackage{ifthen}
\usepackage{xkeyval}
\usepackage{xcolor}
\usepackage{calc}

\usepackage[colorinlistoftodos]{todonotes} 

\title{Full runner removal theorem for Ariki-Koike algebras}
\author[A.~Dell'Arciprete]{Alice Dell'Arciprete}
\address{Department of Mathematics, \\University of York, \\York YO10 5DD, UK.}
\email{alice.dellarciprete@york.ac.uk}
\thanks{}

\numberwithin{equation}{section}
\numberwithin{figure}{section}
\newtheorem{thm}{Theorem}[section]
\newtheorem{prop}[thm]{Proposition}
\newtheorem{lem}[thm]{Lemma}
\newtheorem{cor}[thm]{Corollary}
\newtheorem{conj}[thm]{Conjecture}

\theoremstyle{definition}
\newtheorem{defn}[thm]{Definition}
\newtheorem{exe}[thm]{Example}

\theoremstyle{rem}
\newtheorem{rmk}[thm]{Remark}

\newcommand{\Z}{\mathbb{Z}}
\newcommand{\la}{\lambda}

\newcommand{\C}{\mathbb{C}}

\newcommand\zez{\mathbb{Z}/e\mathbb{Z}}
\newcommand\UU{\mathcal{U}}

\definecolor{bead}{gray}{0.2}

\newcommand\tl{\begin{picture}(8,4)
\put(4,-1){\line(0,1){3}}
\put(4,2){\line(1,0){7}}
\end{picture}}

\newcommand{\bd}{\begin{picture}(8,6)
\put(4,-1){\line(0,1){8}}
\put(4,3){\circle*{6}}
\end{picture}}
\newcommand{\nb}{\begin{picture}(8,6)
\put(4,-1){\line(0,1){8}}
\put(3,3){\line(1,0){2}}
\end{picture}}

\newcommand\bigtl{\begin{picture}(16,8)
\put(8,-2){\line(0,1){6}}
\put(8,4){\line(1,0){14}}
\end{picture}}
\newcommand\bigtr{\begin{picture}(16,8)
\put(8,-2){\line(0,1){6}}
\put(-6,4){\line(1,0){14}}
\end{picture}}
\newcommand\bigtm{\begin{picture}(16,8)
\put(8,-2){\line(0,1){6}}
\put(-6,4){\line(1,0){28}}
\end{picture}}

\newcommand{\bigbd}{\begin{picture}(16,12)
\put(8,-2){\line(0,1){22}}
\put(8,6){\circle*{12}}
\end{picture}}
\newcommand{\bignb}{\begin{picture}(16,12)
\put(8,-2){\line(0,1){22}}
\put(6,6){\line(1,0){4}}
\end{picture}}

\newcommand{\bigvd}{\begin{picture}(16,20)
\put(8,10){\circle*{2}}
\put(8,4){\circle*{2}}
\put(8,16){\circle*{2}}
\end{picture}}

\begin{document}
\begin{abstract}
We consider the representation theory of the Ariki-Koike algebra, a $q$-deformation of the group algebra of the complex reflection group $C_r \wr S_n$. 
We define the addition of a runner full of beads for the abacus display of a multipartition and investigate some combinatorial properties of this operation. We focus our attention on the $v$-decomposition numbers, i.e. the polynomials arising from the Fock space representation of the quantum group $U_v(\widehat{\mathfrak{sl}}_e)$. Using Fayers' LLT-type algorithm for Ariki-Koike algebras, we relate $v$-decomposition numbers for different values of $e$ for the class of $e$-multiregular multipartitions, by adding a full runner of beads to each component of the abacus displays for the labelling multipartitions.
\end{abstract}
\maketitle
\section{Introduction}
Let $n$ be a positive integer. Let $\mathfrak{S}_n$ denote the symmetric group of degree $n$. This has the famous Coxeter presentation with generators $T_1,\ldots, T_{n-1}$ and relations
\begin{align*}
T_i^2&= 1, & \text{ for } 1 &  \leq i \leq n - 1,\\
T_i T_j &= T_j T_i, & \text{ for }1 &  \leq  i < j - 1 \leq n - 2,\\
T_i T_{i+1} T_i &= T_{i+1} T_i T_{i+1}, &  \text{ for } 1 &  \leq  i \leq  n-2.
\end{align*}
If we view this as a presentation for a (unital associative) algebra over a field $\mathbb{F}$, then the algebra we get is the group algebra $\mathbb{F}\mathfrak{S}_n$. Let $q$ be a non-zero element of $\mathbb{F}$. Now we can introduce a ‘deformation’, by replacing the relation $T_i^2=1$ with
$$(T_i + q)(T_i - 1) = 0$$ for each $i$. The resulting algebra is the Iwahori-Hecke algebra $H_n = H_{\mathbb{F},q}(\mathfrak{S}_n)$ of the symmetric group $\mathfrak{S}_n$. This algebra (of which the group algebra
$\mathbb{F}\mathfrak{S}_n$ is a special case) arises naturally, and its representation theory has been extensively studied. An excellent introduction to this theory is provided by Mathas's book \cite{Mat99}. As long as $q$ is non-zero, the representation theory of $H_n$ bears a remarkable resemblance to the representation theory of $\mathfrak{S}_n$. Indeed, there are many theorems concerning $H_n$ which reduce representation-theoretic notions to statements about the combinatorics of partitions.

In this paper, we consider the representation theory of the Ariki-Koike algebra, that can be thought as well as a generalisation of the symmetric group algebra since it includes $\mathbb{F}\mathfrak{S}_n$ as a special case. This algebra is a deformation of the group algebra of the complex reflection group $C_r \wr \mathfrak{S}_n$, defined using parameters $q, Q_1, \ldots, Q_r \in \mathbb{F}$ and it is usually denoted by $\mathcal{H}_{\mathbb{F},q,\bm Q}(C_r \wr \mathfrak{S}_n)$ where $\bm Q = \{Q_1, \ldots, Q_r\}$. For brevity, we will often write $\mathcal{H}_{r,n}$ instead of $\mathcal{H}_{\mathbb{F},q,\bm Q}(C_r \wr \mathfrak{S}_n)$.

The Ariki-Koike algebras were introduced by Ariki and Koike in \cite{AK94}. The representation theory of these algebras is beginning to be well understood. For example, the simple modules of the Ariki-Koike algebras have been classified; the blocks are known; and, in principle, the decomposition matrices of the Ariki-Koike algebras can be computed in characteristic zero. A comprehensive review of the representation theory of the Ariki-Koike algebras can be found in Mathas's paper \cite{Mat04}. In many respects it seems that the Ariki-Koike algebra behaves in the same way as the Iwahori-Hecke algebra $H_n$; many of the combinatorial theorems concerning $H_n$ have been generalised to the Ariki-Koike algebra, with the role of partitions being played by multipartitions. In fact, much of the difficulty of understanding the Ariki-Koike algebra seems to lie in finding the right generalisations of the combinatorics of partitions to multipartitions - very simple combinatorial notions (such as the definition of an $e$-regular partition) can have rather nebulous generalisations (such as ‘dual Kleshchev’ multipartitions). 

Similarly to Iwahori-Hecke algebras, the main problem of interest in the representation theory of Ariki-Koike algebras is the \textit{decomposition number problem}, which asks for the composition multiplicities of simple modules in the so called Specht modules. The \textit{decomposition matrix} records these multiplicities.

It is known that computing the decomposition numbers in the case $\mathbb{F} = \mathbb{C}$ is an important first step in working out the decomposition numbers over any field (see \cite{Gec92, Gec98}). 
In fact, the following result holds. Let $\bm D_p$ be the decomposition matrix of $\mathcal{H}_{\mathbb{F}_p,q,\bm Q}(C_r \wr \mathfrak{S}_n)$ with $\mathbb{F}_p$ a field of characteristic $p>0$ and $\bm D$ be the decomposition matrix of $\mathcal{H}_{\mathbb{C},\zeta,\bm Q}(C_r \wr \mathfrak{S}_n)$. Then there exists a square unitriangular matrix $\bm A$, called \textit{adjustment matrix}, such that $$\bm D_p=\bm D \bm A.$$

Fortunately, the decomposition numbers $d_{\bm\la\bm\mu}$ can be computed when $\mathbb{F} = \mathbb{C}$; they are the values at $v = 1$ of certain polynomials $d_{\bm\lambda \bm\mu}(v)$, which have accordingly become known as `$v$-decomposition numbers'. This result has first been conjectured for Iwahori–Hecke algebras by Lascoux, Leclerc and Thibon \cite{LLT96} and proved for the wider class of Ariki-Koike algebras by Ariki \cite{A96}. It is by far the most significant theorem in this regard. The $v$-decomposition numbers arise from the Fock space representation of the quantum group $U_v(\widehat{\mathfrak{sl}}_e)$. This has a natural basis indexed by the set of partitions for $H_n$ (respectively, of multipartitions for $\mathcal{H}_{r, n}$), and a `canonical basis' which is invariant under the bar involution. The $v$-decomposition numbers are the entries of the transition matrix between these two bases.

For Iwahori-Hecke algebras of $\mathfrak{S}_n$, there is a fast algorithm due to Lascoux, Leclerc and Thibon \cite{LLT96} for computing the canonical basis and so the $v$-decomposition numbers.

For Ariki-Koike algebras, there are different generalisations of this algorithm due to Jacon \cite{Jac05}, Yvonne \cite{Yvo07} and Fayers \cite{Fay10}. We will use the one presented in \cite{Fay10} because it adapts better to our purposes.

Since the decomposition numbers in characteristic 0 can be computed by the LLT algorithm and its generalisations, in effect the problem of determining the decomposition matrices in arbitrary characteristic is equivalent to computing adjustment matrices. However, not a great deal is known about adjustment matrices; the most general statement we have about adjustment matrices is James’s Conjecture. 

\begin{conj}[James’s Conjecture]
Let $\mathbb{F}$ be a field of characteristic $p>0$ and suppose that $\bm D_p=\bm D \bm A$. If $n < pe$, then the adjustment matrix $\bm A$ is the identity matrix.
\end{conj}

For Iwahori–Hecke algebras, this conjecture has been verified for blocks of weight at most four, thanks to the work of Richards \cite{Ric96} and Fayers \cite{Fay07conj, Fay08conj}. However, after being a central focus of research in representation theory for thirty years, this conjecture was finally shown to be false by Williamson in \cite{will17}.

In the attempt of tackling the decomposition number problem, in \cite{JM02} James and Mathas proved the so called `empty' runner removal theorem in which they relate $v$-decomposition numbers of Iwahori-Hecke algebras for different values of $e$, by adding empty runners to the abacus displays for the labelling partitions. After that, in \cite{Frunrem} Fayers proves a similar theorem, which involves adding full runners to these abacus displays. In this paper, for a class of multipartitions, called $e$-multiregular, we generalise Fayers' theorem to the Ariki-Koike algebras showing that the $v$-decomposition numbers $d_{\bm\lambda\bm\mu}(v)$ and $d_{\bm\lambda^+\bm\mu^+}(v)$ coincide, where $\bm\lambda^+$ and $\bm\mu^+$ are the multipartitions obtained from the $e$-abacus display of $\bm\lambda$ and $\bm\mu$ by adding a runner full of beads in each of their components.

\begin{thm}[Theorem \ref{canbasis_runrem}]\label{thm2}
Let $\bm\lambda, \bm\mu$ be $r$-multipartitions in a block $B$ of $\mathcal{H}_{\mathbb{F},q,\bm Q}(C_r \wr \mathfrak{S}_n)$ with $\bm\mu$ $e$-multiregular. 
If the new inserted runners defining $\bm\la^+$ and $\bm\mu^+$ are `long enough', then
$$d_{\bm\lambda\bm\mu}(v)=d_{\bm\lambda^{+}\bm\mu^{+}}(v).$$
\end{thm}

In Section \ref{sec:basicdef}, we define the algebras that we will work with, along with giving an overview of any background material that we will need in order to study their representation theory. This will include both the algebraic setup we require and some combinatorial definitions such as partitions, tableaux and abacuses together with their generalisations for  Ariki-Koike algebras. We also define Kleshchev and dual Kleshchev multipartitions and give an explicit description of the relation between these two definitions.

In Section \ref{sec:LLTalg}, we introduce the Fock space representation of the quantum group $U_v(\widehat{\mathfrak{sl}}_e)$ and present the LLT-type algorithm for Ariki-Koike algebras given in \cite{Fay10}. This algorithm allows us to generalise the `full' runner removal theorem of Iwahori-Hecke algebras in \cite{Frunrem} to Ariki-Koike algebras. 

Thus, in Section \ref{sec:addFullRun} we define the addition of a runner full of beads in each component of an abacus display of a multipartition. We show that adding a runner full of beads to an abacus display of the empty partition corresponds to a precise sequence of induction operators. We then prove some results that describe how the addition of a full runner interacts with the induction operators. Finally, we use all these properties together with the Fayers' LLT-type algorithm for Ariki-Koike algebras \cite{Fay10} to show that the coefficients of the canonical basis element corresponding to the $e$-multiregular multipartition $\bm\mu$ coincide with the coefficients of the canonical basis element corresponding to the multipartition $\bm\mu^+$.

{Recent work has given us a new line of attack. The cyclotomic quiver Hecke algebras of type $A$, known as KLR algebras (defined independently by Khovanov and Lauda and by Rouquier \cite{KL09,Rou08}), have been shown to be isomorphic to Ariki-Koike algebras by Brundan and Kleshchev in \cite{BK09}. Via this isomorphism, the $\Z$-grading of the KLR algebras can be used in the setting of Ariki-Koike algebras, and thus graded Specht modules \cite{BKW11} and graded decomposition numbers can be studied. However, this is beyond the scope of what we are going to consider in this paper.}

\section{Basic definitions}\label{sec:basicdef}

\subsection{The Ariki-Koike algebras}\label{AKdef}
Let $r\geq 1$ and $n \geq0$. Let $W_{r,n}$ be the complex reflection group $C_r \wr~\mathfrak{S}_n$. 
We can define the Ariki-Koike algebra as a deformation of the group algebra $\mathbb{F}W_{r,n}$.
\begin{defn}
Let $\mathbb{F}$ be a field and $q,Q_1,\ldots,Q_r$ be elements of $\mathbb{F}$, with $q$ non-zero. Let $\bm{Q}=(Q_1, \ldots, Q_r)$. The \textbf{Ariki-Koike algebra} $\mathcal{H}_{\mathbb{F},q,\bm{Q}}(W_{r,n})$ of $W_{r,n}$ is defined to be the unital associative $\mathbb{F}$-algebra with generators $T_0, \ldots, T_{n-1}$ and relations
\begin{align*}
(T_0-Q_1)\cdots(T_0-Q_r)&=0,\\
T_0T_1T_0T_1&= T_1T_0T_1T_0,\\
(T_i+1)(T_i-q)&=0, & \text{for }1 &\leq i \leq n - 1,\\
T_iT_j &= T_jT_i, &\text{for }0 &\leq i < j-1 \leq n - 2,\\
T_iT_{i+1}T_i &= T_{i+1}T_iT_{i+1}, &\text{for } 1&\leq  i \leq n - 2.\\
\end{align*}
\end{defn}
For brevity, we may write $\mathcal{H}_{r,n}$ for $\mathcal{H}_{\mathbb{F},q, \bm{Q}}(W_{r,n})$. Define $e$ to be minimal such that $1 + q + \ldots + q^{e-1}= 0$, or set $e = \infty$ if no such value exists. Throughout this paper we shall assume that $e$ is finite and we shall refer to $e$ as the \textbf{quantum characteristic}. Set $I = \{0, 1, \ldots, e-1\}$: we will identify $I$ with $\mathbb{Z}/e\mathbb{Z}$. We refer to any $r$-tuple of integers $(a_1, \ldots, a_r) \in \mathbb{Z}^r$ as a \textbf{multicharge}. We say $\bm{Q}$ is \textbf{$q$-connected} if, for each $j\in\{1, \ldots, r\}$, $Q_j = q^{a_j}$ for some $a_j \in \mathbb{Z}$. In \cite{dm02}, Dipper and Mathas prove that any Ariki-Koike algebra is Morita equivalent to a direct sum of tensor products of smaller Ariki-Koike algebras, each of which has $q$-connected parameters. Thus, we may assume that we are always working with a Ariki-Koike algebra with each $Q_j$ being an integral power of $q$. So we assume that we can find an $r$-tuple of integers $\bm{a} = (a_1, \ldots, a_r)$ such that $Q_j = q^{a_j}$ for each $j$ where $q \neq 1$, so that $q$ is a primitive $e^{\text{th}}$ root of unity in $\mathbb{F}$. We call such $\bm{a}$ a \textbf{multicharge} of $\mathcal{H}_{r,n}$.
Since $e$ is finite then we may change any of the $a_j$ by adding a multiple of $e$, and we shall still have $Q_j = q^{a_j}$. 

For the rest of this section we will fix a multicharge $\bm{a}=(a_1, \ldots, a_r)\in \mathbb{Z}^r$. Notice that at this stage we are not requiring that $\bm{a}$ is a multicharge of the Ariki-Koike algebra.

\subsection{Multipartitions}\label{multpar}
A \textbf{partition} of $n$ is defined to be a non-increasing sequence $\lambda = (\lambda_1, \lambda_2, \dots)$ of non-negative integers whose sum is $n$. {The integers $\lambda_b$, for $b\geq1$, are called the \textbf{parts} of $\lambda$. We write $|\lambda| = n$.\\
Since $n < \infty$, there is a $k$ such that $\lambda_b = 0$ for $b > k$ and we write $\lambda = (\lambda_1, \dots , \lambda_k)$. We write $\varnothing$ for the unique empty partition $(0, 0, \dots, 0)$. If a partition has repeated parts, for convenience we group them together with an index. For example, $$(4,4,2,1,0,0, \dots) = (4,4,2,1) = (4^2,2,1).$$
If $\lambda$ is a partition, we define the \textbf{conjugate partition} $\lambda'$ of $\lambda$ to be the partition with $b^{\text{th}}$ part $\lambda'_b = \#\{c \geq 1\text{ }|\text{ } \lambda_c \geq b\}.$
The \textbf{Young diagram} of a partition $\lambda$ is the subset
$$[\lambda] := \{(b,c) \in \mathbb{N}_{>0} \times \mathbb{N}_{>0
}\text{ } | \text{ } c\leq \lambda_b\}.$$
Let $2\leq e<\infty$ be a natural number. For each $l\geq 1$, we define the $l^{\text{th}}$ \textbf{ladder} in $\mathbb{N}^2$ to be the set
$$\mathcal{L}_l = \{(b,c) \in \mathbb{N}^2_{>0} \text{ }| \text{ }b + (e - 1)(c - 1) = l\}.$$
All the nodes in $\mathcal{L}_l$ have the same residue (namely, $ 1 - l \mod e$), and we define the residue of $\mathcal{L}_l$ to be
this residue. If $\la$ is a partition, we define the $l^{\text{th}}$ ladder $\mathcal{L}_l(\la)$ of $\la$ to be the intersection of $\mathcal{L}_l$ with the Young diagram of $\la$.

\begin{exe} Suppose $e = 3$, and $\la = (4, 3, 1)$. Consider the Young diagram of $\la$. Then in the first diagram we label each node of $[\la]$ with the number of the ladder in which it lies, while in the second one we filled the nodes with their residues:
\Yvcentermath1
$$\young(1357,246,3) \quad \text{ and } \quad \young(0120,201,1).$$
\end{exe}
Other useful definitions about partitions we recall are the following.
\begin{itemize}
    \item A partition $\lambda$ is \textbf{$e$-singular} if $\lambda_{i+1} = \lambda_{i+2} = \ldots = \lambda_{i+e} > 0$ for some $i$. Otherwise, $\lambda$ is \textbf{$e$-regular}. For example, for $e=3$ the partition $(4, 3, 1)$ is $3$-regular, while the partition $(4,1^4)$ is $3$-singular.
    \item A partition $\lambda$ is \textbf{$e$-restricted} if $\lambda_i -\lambda_{i+1} < e$ for every $i \geq1$.
\end{itemize}
Notice that a partition $\lambda$ is $e$-regular if and only if its conjugate $\lambda'$ is $e$-restricted.

\begin{defn}\label{multipar}
An $r$-\textbf{multipartition} of $n$ is an ordered $r$-tuple $\bm{\lambda} = (\lambda^{(1)}, \dots, \lambda^{(r)})$  of partitions such that $$|\bm{\lambda}| := |\lambda^{(1)}|+ \ldots +|\lambda^{(r)}| = n.$$
If $r$ is understood, we shall just call this a multipartition of $n$.
We write $\mathcal{P}^r$ for the set of $r$-multipartitions.
\end{defn}

As with partitions, we write the unique multipartition of $0$ as $\bm\varnothing$. The \textbf{Young diagram} of a multipartition $\bm\lambda$ is the subset
$$
[\bm{\lambda}]:=\{(b,c,j)\in \mathbb{N}_{>0}\times \mathbb{N}_{>0}\times\{1, \ldots, r\} \text{ }|\text{ } c \leq \lambda_b^{(j)}\}.
$$
We may abuse notation by not distinguishing a multipartition from its Young diagram.}
The elements of $[\bm{\lambda}]$ are called \textbf{nodes} of $\bm\lambda$. We say that a node $\mathfrak{n} \in [\bm{\lambda}]$ is \textbf{removable} if $[\bm{\lambda}]\setminus \{\mathfrak{n}\}$ is also the Young diagram of a multipartition. We say that an element $\mathfrak{n} \in \mathbb{N}^2_{>0}\times \{1, \ldots, r\}$ is an \textbf{addable node} if $\mathfrak{n} \notin [\bm{\lambda}]$ and $[\bm{\lambda}] \cup \{\mathfrak{n}\}$ is the Young diagram of a multipartition. 
Define a bijection $'$ from $\mathbb{N}_{>0}^{2}\times \{1,\ldots, r\}$ to itself by
$$(b,c,j)'=(c,b,r+1-j).$$
\begin{defn}
Given a multipartition $\bm\la =(\la^{(1)}, \ldots, \la^{(r)})$, define the \textbf{conjugate multipartition} to have Young diagram
$$[\bm\la']=\{\mathfrak{n}' \text{ }|\text{ } \mathfrak{n}\in \bm\la\};$$
that is, $\bm\la'=({\la^{(r)}}' ,\ldots, {\la^{(1)}}')$, where ${\la^{(j)}}'$ is the usual conjugate partition to $\la^{(j)}$.
\end{defn}

\begin{exe}
Consider $\bm\la = ((2^2, 1), (1^2), (3, 1))$. Then $\bm\la'=((2,1^2),(2), (3,2))$ and its Young diagram is
\Yvcentermath1
$$\begin{matrix*}[l]
\yng(2,1^2) \\
\\
\yng(2)\\
\\
\yng(3,2) \end{matrix*}.$$
\end{exe}

{Given a multicharge $\bm{a}=(a_1, \ldots, a_r)$, to each node $(b,c,j) \in [\bm\lambda]$ we associate its \textbf{residue} $\mathrm{res}_{\bm{a}}(b,c,j) = a_j + c-b$ $(\text{mod }e)$. We draw the residue diagram of $\bm\lambda$ by replacing each node in the Young diagram by its residue.}
\begin{exe}\label{exsameblock}
Suppose $r=3$ and $\bm{a}=(1,0,2)$. Let $\bm{\lambda}=((1^2),(2),(2,1))$ and $\bm\mu=((1),(2,1),(1^3))$ be two multipartitions of $7$. If $e=4$, then the $4$-residue diagram of $[\bm{\lambda}]$ and of $[\bm\mu]$ are
$$
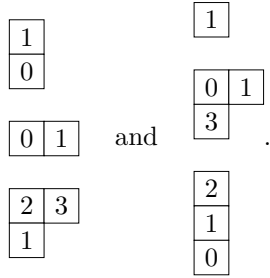
\begin{matrix*}[l]\young(1,0) \\ \\ \young(01) \\ \\ \young(<2><3>,1)\end{matrix*} \quad \text{ and } \quad\begin{matrix*}[l] \young(1) \\ \\ \young(01,3) \\  \\\young(2,1,0)\end{matrix*}.
$$
\end{exe}

{
\subsection{Dual Kleshchev multipartitions}\label{ssec:Kl_mpt}
Residues of nodes are also useful in classifying the simple $\mathcal{H}_{r,n}$-modules. Indeed, the notion of residue helps us to describe a certain subset $\mathcal{K}$ of the set of all multipartitions, which index the simple modules for $\mathcal{H}_{r,n}$.
We impose partial order $\succ$ on the set of nodes of residue $i \in I$ of a multipartition by saying that $(b,c,j)$ is \textbf{below} $(b',c',j')$ (or $(b',c',j')$ is \textbf{above} $(b,c,j)$) if either $j>j'$ or ($j=j'$ and $b>b'$). In this case we write $(b,c,j)\succ(b',c',j')$. Note this order restricts to a total order on the set of all addable and removable nodes of residue $i \in I$ of a multipartition.

Suppose $\bm\lambda$ is a multipartition, and given $i \in I$ define the \textbf{$i$-signature of $\bm\lambda$ with respect to $\succ$} by examining all the addable and removable $i$-nodes of $\bm\lambda$ in turn from lower to higher, and writing a $+$ for each addable node of residue $i$ and a $-$ for each removable node of residue $i$. Now construct the \textbf{reduced $i$-signature with respect to $\succ$} by successively deleting all adjacent pairs $-+$. If there are any $-$ signs in the reduced $i$-signature of $\bm\lambda$, the lowest of these nodes is called a \textbf{good node of $[\bm\lambda]$ with respect to $\succ$}.

\begin{defn}
We say that $\bm\la$ is a \textbf{dual Kleshchev multipartition} if and only if there is a sequence
$$\bm\lambda = \bm\lambda(n), \bm\lambda(n-1), \ldots, \bm\lambda(0)= \bm\varnothing$$
of multipartitions such that for each $k$, $[\bm\lambda(k-1)]$ is obtained from $[\bm\lambda(k)]$ by removing a good node with respect to $\succ$.
\end{defn}

This definition depends on the multicharge $\bm a=(a_1, \ldots, a_r)$. We write $\mathcal{K}'(a_1, \ldots, a_r)$ for the set of dual Kleshchev multipartitions with multicharge  $(a_1, \ldots, a_r)$.
The reason why we call these multipartitions \textit{dual} Kleshchev is because they are the dual version of the multipartitions usually called Kleshchev or restricted in \cite{BK09}. There is a canonical bijection between the Kleshchev and the dual Kleshchev multipartitions; for details see \cite{mythesis, BK09}. 

\begin{exe}\label{Kl_exe}
Suppose $r=2$ and $e=4$. Consider the multicharge $\bm a =(1,0)$. Then the multipartition $\bm\la=((2,1), (1^3))$ is dual Kleshchev. Indeed, we have the following sequence of multipartitions obtained from $\bm\la$ by removing each time a good node - we write $[\bm\la(k)]\xleftarrow{i}_{\succ} [\bm\la(k-1)]$ to denote that $[\bm\la(k-1)]$ is obtained from $[\bm\la(k)]$ by removing a good $i$-node for $i\in I$ with respect to the total order $\succ$:
$$\Yinternals1\Yaddables1 \begin{matrix*}[l]\Ycorner1\yngres(4,2,1) \\ \\ \Ycorner0\yngres(4,1^3) \end{matrix*} \xleftarrow{2}_{\succ} \begin{matrix*}[l]\Ycorner1\yngres(4,2,1) \\ \\ \Ycorner0\yngres(4,1^2) \end{matrix*} \xleftarrow{2}_{\succ} \begin{matrix*}[l]\Ycorner1\yngres(4,1^2) \\ \\ \Ycorner0\yngres(4,1^2) \end{matrix*} \xleftarrow{0}_{\succ} \begin{matrix*}[l]\Ycorner1\yngres(4,1) \\ \\ \Ycorner0\yngres(4,1^2) \end{matrix*}
\xleftarrow{3}_{\succ}\begin{matrix*}[l]\Ycorner1\yngres(4,1) \\ \\ \Ycorner0\yngres(4,1) \end{matrix*}
\xleftarrow{1}_{\succ}\begin{matrix*}[l] 1 \\ \\ \Ycorner0\yngres(4,1) \end{matrix*}
\xleftarrow{0}_{\succ}
\begin{matrix*}[l]\varnothing \\ \\ \varnothing \end{matrix*}.$$
\end{exe}

\begin{rmk}\cite{BK09}\label{kl=restr}
If $r = 1$, then the multipartition $(\la)$ is dual Kleshchev if and only if $\la$ is $e$-regular.
\end{rmk}

%
\subsection{Specht module and simple modules}\label{specht}
The algebra $\mathcal{H}_{r,n}$ is a cellular algebra \cite{DJM98,GL96} with the cell modules indexed by $r$-multipartitions of $n$. For each $r$-multipartition $\bm\lambda$ of $n$ we define a $\mathcal{H}_{r,n}$-module $S'(\bm\lambda)$ called  \textbf{dual Specht module}; these modules are the cell modules defined in \cite{Mat03} arising from the cellular basis $\{ n_{\mathfrak{st}}\}$ of $\mathcal{H}_{r,n}$. The details of this definition can be found in \cite[\S 4]{Mat03}. When $\mathcal{H}_{r,n}$ is semisimple, the dual Specht modules form a complete set of non-isomorphic irreducible $\mathcal{H}_{r,n}$-modules. However, we are mainly interested in the case when $\mathcal{H}_{r,n}$ is not semisimple. In this case, we set $D'(\bm\la) = S'(\bm\la)/ \mathrm{rad}\text{ } S'(\bm\la)$ with $\mathrm{rad}\text{ } S'(\bm\la)$ the radical of the bilinear form of $ S'(\bm\la)$ (the form is defined in terms of the structure constants of the basis $\{n_{\mathfrak{s}\mathfrak{t}}\}$). Then, as shown in \cite[Theorem 4.2]{A01}, a complete set of non-isomorphic simple $\mathcal{H}_{r,n}$-modules is given by
$$\{D'(\bm\la)\text{ }|\text{ } \bm\la\text{ is a dual Kleshchev multipartition}\}.$$

}

If $\bm\lambda$ and $\bm\mu$ are $r$-multipartition of $n$ with $\bm\mu$ dual Kleshchev, let $d_{\bm\lambda \bm\mu} = [S'(\bm\lambda):D'(\bm\mu)]$ denote the multiplicity of the simple module $D'(\bm\mu)$ as a composition factor of the Specht module $S'(\bm\lambda)$. The matrix $D = (d_{\bm\lambda \bm\mu})$ is called the \textbf{decomposition matrix} of $\mathcal{H}_{r,n}$ and determining its entries is one of the most important open problems in the representation theory of the Ariki-Koike algebras. It follows from the cellularity of $\mathcal{H}_{r,n}$ that the decomposition matrix is `triangular'; to state this we need to define the dominance order on multipartitions. Given two multipartitions $\bm{\lambda}$ and $\bm{\mu}$ of $n$, we say that $\bm{\lambda}$ \textbf{dominates} $\bm{\mu}$, and write $\bm{\lambda} \trianglerighteq \bm{\mu}$, if {
$$\sum_{a=1}^{j-1}|\lambda^{(a)}| + \sum_{b=1}^{i} \lambda_b^{(j)} \geq \sum_{a=1}^{j-1}|\mu^{(a)}| + \sum_{b=1}^{i}\mu_b^{(j)}$$
for $j=1,2, \ldots, r$ and for all $i \geq 1$}. Then we have the following.
\begin{thm}\cite{DJM98,GL96}\label{decmtxHQ}
Let $\mathbb{F}$ be a field. Suppose $\bm\lambda$ and $\bm\mu$ are $r$-multipartitions of $n$ with $\bm\mu$ dual Kleshchev.
\begin{enumerate}
\item If $\bm\mu=\bm\lambda$, then $[S'(\bm{\lambda}):D'(\bm{\mu})]=1$.
\item If $[S'(\bm{\lambda}):D'(\bm{\mu})]>0$, then $\bm{\lambda} \trianglelefteq \bm{\mu}$.
\end{enumerate}
\end{thm}

Note that the dominance order is a partial order on the set of multipartitions. The dominance order is certainly the `correct' order to use for multipartitions, but it is sometimes useful to have a total order, $>$, on the set of multipartitions. The one we use is given as follows. 
\begin{defn}
Given two multipartitions $\bm\lambda$ and $\bm\mu$ of $n$ with $\bm\lambda \neq \bm\mu$, we write $\bm\lambda > \bm\mu$ if and only if the minimal $j \in \{1, \ldots, r\}$ for which $\lambda^{(j)} \neq \mu^{(j)}$ and the minimal $i\geq 1$ such that $\lambda^{(j)}_i \neq \mu^{(j)}_i$ satisfy $\lambda^{(j)}_i > \mu^{(j)}_i$. This is called the \textbf{lexicographic order} on multipartitions.
\end{defn}

\subsection{$\beta$-numbers and the abacus}\label{betaAbacus}
Given the assumption that the cyclotomic parameters of $\mathcal{H}_{r,n}$ are all powers of $q$, we may conveniently represent multipartitions on an abacus display. Fix $\bm{a}=(a_1, \ldots, a_r)\in \mathbb{Z}^{r}$
to be a multicharge of $\mathcal{H}_{r,n}$.
\begin{defn}
Let $\bm{\lambda}=(\lambda^{(1)}, \ldots, \lambda^{(r)})$ be a multipartition of $n$. For every $i\geq1$ and for every $j \in \{1, \ldots, r\}$, we define the \textbf{$\beta$-number} $\beta_i^j$ to be
$$\beta_i^j := \lambda_i^{(j)} + a_j - i.$$
\end{defn}
The set $B_{a_j}^j =\{ \beta_1^j, \beta_2^j, \ldots \}$ is the set of $\beta$-numbers (defined using the integer $a_j$) of partition $\lambda^{(j)}$. It is easy to see that any set $B_{a_j}^j =\{ \beta_1^j, \beta_2^j, \ldots \}$ is a set containing exactly $a_j + N$ integers greater than or equal to $-N$, for sufficiently large $N$. 

For each set $B_{a_j}^j$ we can define an abacus display in the following way. We take an abacus with $e$ infinite vertical runners, which we label $0, 1, \ldots, e-1$ from left to right, 
and we mark positions on runner $l$ and label them with the integers congruent to $l$ modulo $e$, so that 
position $(x+1)e+l$ lies immediately below position $xe + l$, for each $x$. Now we place a bead at position $\beta_i^j$, for each $i$. The resulting configuration is called $e$-abacus display, or $e$-abacus configuration, for $\lambda^{(j)}$ with respect to $a_j$. Moreover, we say that the bead corresponding to the $\beta$-number $xe+l$ is at \textbf{level} $x$ for $x \in \mathbb{Z}$.

Hence, we can now define the $e$-\textbf{abacus display}, or the $e$-\textbf{abacus configuration}, for a multipartition $\bm{\lambda}$ with respect to $\bm{a}$ to be the $r$-tuple of $e$-abacus displays associated to each component $\lambda^{(j)}$. If it is clear which $e$ we are referring to, we simply say abacus configuration. When we draw abacus configurations we will draw only a finite part of the runners and we will assume that above this point the runners are full of beads and below this point there are no beads. 

\begin{exe}
Suppose that $r = 3$, $\bm{a} = (-1, 0, 1)$ and $\bm{\lambda}= ((1), \varnothing, (1^2))$. Then we have
\begin{align*}
B^1_{-1} &= \{-1,-3,-4,-5, \ldots\};\\
B^2_0 &= \{-1,-2,-3,-4\ldots\};\\
B^3_1 &= \{1,0,-2,-3,-4,\ldots\}.
\end{align*}
So, the abacus display with respect to the multicharge $\bm{a}$ for $\bm{\lambda}$ when $e = 4$ is
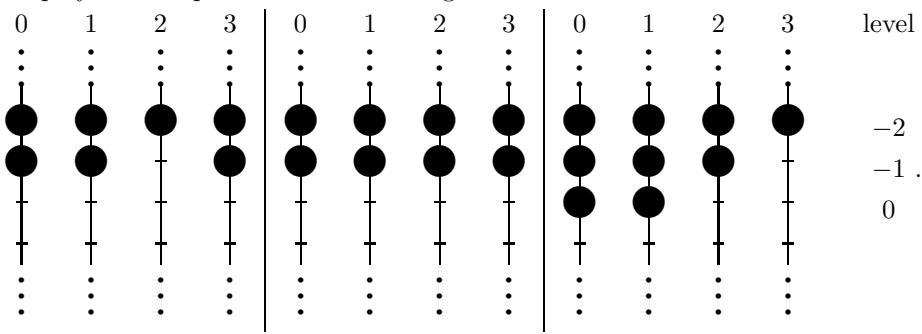
\begin{center}
\begin{tabular}{c|c|c}
$      
        \begin{matrix}

        0 & 1 & 2 & 3 \\

        \bigvd & \bigvd & \bigvd & \bigvd \\

        \bigbd & \bigbd & \bigbd & \bigbd \\
        
        \bigbd & \bigbd & \bignb & \bigbd \\
     
        \bignb & \bignb & \bignb & \bignb \\
        
        \bignb & \bignb & \bignb & \bignb \\

        \bigvd & \bigvd & \bigvd & \bigvd \\

        \end{matrix}
$
&
$     \begin{matrix}

        0 & 1 & 2 & 3 \\

        \bigvd & \bigvd & \bigvd & \bigvd \\

        \bigbd & \bigbd & \bigbd & \bigbd \\

        \bigbd & \bigbd & \bigbd & \bigbd \\

        \bignb & \bignb & \bignb & \bignb \\
        
        \bignb & \bignb & \bignb & \bignb \\

        \bigvd & \bigvd & \bigvd & \bigvd \\

        \end{matrix}
$
&
$     \begin{matrix}

        0 & 1 & 2 & 3 & & \text{level}\\

        \bigvd & \bigvd & \bigvd & \bigvd & &\\

        \bigbd & \bigbd & \bigbd & \bigbd & & -2\\

        \bigbd & \bigbd & \bigbd & \bignb & & -1\\

        \bigbd & \bigbd & \bignb & \bignb & & 0\\
        
        \bignb & \bignb & \bignb & \bignb & & \\

        \bigvd & \bigvd & \bigvd & \bigvd & & \\

        \end{matrix}
.$\\
\end{tabular}
\end{center}
\end{exe}

\subsection{Rim $e$-hooks and $e$-core}
We recall some other definitions about the diagram of a partition. 

\begin{defn}
Suppose $\lambda$ is a partition of $n$ and $(i,j)$ is a node of $[\lambda]$.
\begin{enumerate}
\item The \textbf{rim} of $[\lambda]$ is defined to be the set of nodes
$$\{(i,j) \in [\lambda] \text{ }|\text{ }(i+1,j+1) \notin [\lambda]\}.$$
\item Define an \textbf{$e$-rim hook} to be a connected subset $R$ of the rim containing exactly $e$ nodes such that $[\lambda]  \setminus R$ is the diagram of a partition.
\item  If $\lambda$ has no $e$-rim hooks then we say that $\lambda$ is an \textbf{$e$-core}.
\end{enumerate}
\end{defn}


An abacus display for a partition is useful for visualising the removal of $e$-rim hooks. If we are given an abacus display for $\lambda$ with $\beta$-numbers in a set $B$, then $[\lambda]$ has a $e$-rim hook if and only if there is a $\beta$-number $\beta_i \in B$ such that $\beta_i-e \notin B$. Furthermore, removing a $e$-rim hook corresponds to reducing such a $\beta$-number by $e$. On the abacus, this corresponds to sliding a bead up one position on its runner. So, $\lambda$ is an $e$-core if and only if every bead in the abacus display has a bead immediately above it. Using this, we can see that the definition of 
$e$-core of $\lambda$ is well defined.
%
{Finally, we can also notice that each bead corresponds to the end of a row of the diagram of $\lambda$ (or to a row of length $0$).

For multipartitions, by the definition of the $\beta$-numbers the node at the end of the row (if it exists) has residue $i$ if and only if the corresponding bead is on runner $i$ of the abacus for $i \in I$.
Thus, if we increase any $\beta$-number by one, this is equivalent to moving a bead from runner $i$ to runner $i + 1$ which is equivalent to adding a node of residue $i + 1$ to the diagram of a component $\lambda^{(j)}$. Similarly, decreasing a $\beta$-number by one is equivalent to moving a bead from runner $i$ to runner $i - 1$ which is equivalent to removing a node of residue $i$ from the diagram of $\lambda^{(j)}$.}

\section{An LLT-type algorithm for Ariki-Koike algebras}\label{sec:LLTalg}

In this section we consider the integrable representation theory of the quantised enveloping algebra $\UU=U_v(\widehat{\mathfrak{sl}}_e)$. For any dominant integral weight $\Lambda$ for $\UU$, the irreducible highest-weight module $V(\Lambda)$ for $\UU$ can be constructed as a submodule $M^{\bm s}$ of a \emph{Fock space} $\mathcal{F}^{\bm s}$ (which depends not just on $\Lambda$ but on an ordering of the fundamental weights involved in $\Lambda$). Using the standard basis of the Fock space, one can define a \emph{canonical basis} 
 for $M^{\bm s}$. There is considerable interest in computing this canonical basis (that is, computing the transition coefficients from the canonical basis to the standard basis) because of Ariki's theorem, which says that these coefficients, evaluated at $v=1$, give the decomposition numbers for Ariki-Koike algebras.  In the case where $\Lambda$ is of level $1$, the LLT algorithm due to Lascoux, Leclerc and Thibon \cite{LLT96} computes the canonical basis. The purpose of this section is to present the generalisation of this algorithm to higher levels given by Fayers in \cite{Fay10}.
Fayers' algorithm computes the canonical basis for an intermediate module $M^{\otimes\bm s}$, which is defined to be the tensor product of level $1$ highest-weight irreducible modules.  It is then straightforward to discard unwanted vectors to get the canonical basis for $M^{\bm s}$. In order to describe this algorithm we need to introduce some notation as we do in the next subsection.

%
%

\subsection{The quantum algebra $U_v(\widehat{\mathfrak{sl}}_e)$ and the Fock space}


%
%
%
%
%
We introduce the following notation for multipartitions. If $\bm\la=(\la^{(1)},\dots,\la^{(r)})$ is an $r$-multipartition for $r>1$, then we write $\bm\la_-$ for the ($r-1$)-multipartition $(\la^{(2)},\dots,\la^{(r)})$.  If $\bm\nu$ is an ($r-1$)-multipartition, we write $\bm\nu_+$ for the $r$-multipartition $(\varnothing,\nu^{(1)},\dots,\nu^{(r-1)})$.  Finally, if $\bm\mu$ is an $r$-multipartition, we write $\bm\mu_0$ for the $r$-multipartition $(\bm\mu_-)_+=(\varnothing,\mu^{(2)},\dots,\mu^{(r)})$.

%
\begin{defn}
We say that 
a multipartition $\bm\la$ is \textbf{$e$-multiregular} if $\la^{(j)}$ is $e$-regular for each $j$.  We write $\mathcal{R}$ for the set of $e$-regular partitions and $\mathcal{R}^r$ for the set of all $e$-multiregular $r$-multipartitions, if $e$ is understood.
\end{defn}

%
%
%

Denote with $\UU$ the quantised enveloping algebra $U_v(\widehat{\mathfrak{sl}}_e)$. This is a $\mathbb{Q}(v)$-algebra with generators $e_i,f_i$ for $i\in I$ and $v^h$ for $h\in P^\vee$, where $P^\vee$ is a free $\mathbb{Z}$-module with basis $\{h_i\mid i\in I\}\cup\{d\}$. Denote the dual basis of $(P^{\vee})^*$ by $\{\Lambda_0, \ldots, \Lambda_{e-1}, \delta\}$. These generators are subjected to the following relations
$$v^hv^{h'}=v^{h+h'}, \quad v^0 = 1,$$
$$v^h e_j v^{-h}=v^{\langle\alpha_j,h\rangle}e_j,$$
$$v^h f_j v^{-h}=v^{-\langle\alpha_j,h\rangle}f_j,$$
$$[e_i,f_j]= \delta_{ij}\dfrac{v^{h_i}-v^{-h_i}}{v-v^{-1}},$$
$$\sum_{k=0}^{1-\langle\alpha_i,h_j\rangle}(-1)^k \left[\begin{matrix}
1-\langle\alpha_i,h_j\rangle\\
k
\end{matrix}\right]e_i^{1-\langle \alpha_i,h_j \rangle -k}e_je_i^k =0\quad (i \neq j),
$$
$$\sum_{k=0}^{1-\langle\alpha_i,h_j\rangle}(-1)^k \left[\begin{matrix}
1-\langle\alpha_i,h_j\rangle\\
k
\end{matrix}\right]f_i^{1-\langle \alpha_i,h_j \rangle -k}f_jf_i^k=0 \quad (i \neq j),
$$
where $\alpha_j=2\Lambda_i - \Lambda_{i-1} - \Lambda_{i+1} +\delta_{i0}\delta$ for $i=0,1, \ldots, e-1$.
Here we follow the usual notation for $v$-integers, $v$-factorials and $v$-binomial coefficients:
$$[k]=\dfrac{v^{k}-v^{-k}}{v-v^{-1}}, \quad [k]!=[k][k-1] \cdots [1], \quad \left[\begin{matrix}
m\\k
\end{matrix}\right]=\dfrac{[m]!}{[m-k]![k]!}.$$
For any integer $m>0$, we write $f_i^{(m)}$ to denote the quantum divided power $f_i^m/[m]!$.

Moreover, $\UU$ is a Hopf algebra with comultiplication denoted $\Delta$ in \cite{Kas02} and hence the tensor product of two $\UU$-modules can be regarded as a $\UU$-module.

The $\mathbb{Q}$-linear ring automorphism $\overline{\phantom{o}}:\UU\to\UU$ defined by
\[\overline{e_i}=e_i,\qquad \overline{f_i}=f_i,\qquad \overline v=v^{-1},\qquad \overline{v^h}=v^{-h}\]
for $i\in I$ and $h\in P^\vee$ is called the \textbf{bar involution}.

Now we fix $\bm s\in I^r$ for some $r\geq1$, and define the \textbf{Fock space} $\mathcal{F}^{\bm s}$ to be the $\mathbb{Q}(v)$-vector space with a basis $\{\bm\la \mid \bm\la\in\mathcal{P}^r\}$, which we call the \textbf{standard basis}.  This has the structure of a $\UU$-module: for a full description of the module action, we refer to \cite{Fay10}. Here, we describe the action of the generators $f_0, \ldots, f_{e-1}$.

Given $\bm\la, \bm\xi \in\mathcal{P}^r$, then we write $\bm\la\xrightarrow{m:i}\bm\xi$ to indicate that $\bm\xi$ is obtained from $\bm\la$ by adding $m$ addable $i$-nodes. If this is the case, then we consider the total order $>$ on addable and removable nodes and so we define the integer
\begin{align}\label{N_i(la,xi)}
N_i(\bm\la,\bm\xi)= \sum_{\mathfrak{n}\in\bm\xi\setminus\bm\la}&((\text{number of addable $i$-nodes of $\bm\xi$ above $\mathfrak{n}$}) \nonumber \\
&- (\text{number of removable $i$-nodes of $\bm\la$ above $\mathfrak{n}$})).
\end{align}
Now the action of $f_i^{(m)}$ is given by
$$f_i^{(m)}\bm\la = \sum_{\bm\la\xrightarrow{m:i}\bm\xi}v^{N_i(\bm\la,\bm\xi)}\bm\xi.$$


\begin{prop}\label{Ncomps}
Let $i\in I$. Suppose $\bm\la$ and $\bm\xi$ are $r$-multipartitions such that $\bm\la\xrightarrow{m:i}\bm\xi$. Then
$$N_i(\bm\la, \bm\xi)= \sum_{\mathfrak{n}\in \bm\xi\setminus\bm\la}\left(N_i(\la^{(J_{\mathfrak{n}})},\xi^{(J_{\mathfrak{n}})})+\sum_{j=1}^{J_{\mathfrak{n}}-1} N_i(\la^{(j)}, \xi^{(j)})\right),$$
where $J_{\mathfrak{n}}$ is the component of $\mathfrak{n}$ in $\bm\xi$. 
\end{prop}
\begin{proof}
This follows from the definition of $N_i(\bm\la, \bm\xi)$ and from the total order $>$ on the set of all addable and removable nodes in a multipartition (Subsection \ref{ssec:Kl_mpt}). Indeed, for each $\mathfrak{n}\in \bm\xi\setminus\bm\la$, a term of $N_i(\bm\la, \bm\xi)$ consists of \begin{equation*}
\#\{\text{addable $i$-nodes of $\bm\xi$ above }\mathfrak{n}\}-\#\{\text{removable $i$-nodes of $\bm\la$ above }\mathfrak{n}\},
\end{equation*}
and the nodes above $\mathfrak{n}$ are exactly those above $\mathfrak{n}$ in the component $J_{\mathfrak{n}}$ and all the nodes in components $j$ with $j<J_{\mathfrak{n}}$. So, for each $\mathfrak{n}\in \bm\xi\setminus\bm\la$, a term of $N_i(\bm\la, \bm\xi)$ consists of
$$N_i(\la^{(J_{\mathfrak{n}})},\xi^{(J_{\mathfrak{n}})})+\sum_{j=1}^{J_{\mathfrak{n}}-1} N_i(\la^{(j)}, \xi^{(j)})$$
where $N_i(\la^{(J_{\mathfrak{n}})},\xi^{(J_{\mathfrak{n}})})$ is given by \eqref{N_i(la,xi)}, and for $j<J_{\mathfrak{n}}$ 
$$N_i(\la^{(j)}, \xi^{(j)}) =\#\{\text{addable $i$-nodes of $\xi^{(j)}$}\}-\#\{\text{removable $i$-nodes of $\la^{(j)}$}\}.$$
\end{proof}

The reason why we are interested in the Fock space is because the submodule $M^{\bm s}$ generated by $\bm \varnothing$ is isomorphic to the irreducible highest-weight module $V(\Lambda_{s_1}+\dots+\Lambda_{s_r})$. This submodule inherits a bar involution from $\UU$: this is defined by $\overline{\bm\varnothing}=\bm\varnothing$ and $\overline{um} = \overline u\,\overline m$ for all $u\in\UU$ and $m\in M^{\bm s}$.  This bar involution allows one to define a \textbf{canonical basis} for $M^{\bm s}$; this consists of vectors $G^{\bm s}(\bm\mu)$, for $\bm\mu$ lying in some subset of $\mathcal{P}^r$ (with our conventions, this is the set of dual Kleshchev multipartitions). These canonical basis vectors are characterised by the following properties:
\begin{itemize}
\item
$\overline{G^{\bm s}(\bm\mu)}=G^{\bm s}(\bm\mu)$;
\item
if we write $G^{\bm s}(\bm \mu) = \sum_{\bm\la\in\mathcal{P}^r}d^{\bm s}_{\bm\la\bm\mu}(v)\bm\la$ with $d^{\bm s}_{\bm\la\bm\mu}(v)\in\mathbb{Q}(v)$, then $d^{\bm s}_{\bm\mu\bm\mu}(v)=1$, while $d^{\bm s}_{\bm\la\bm\mu}(v)\in v\mathbb{Z}[v]$ if $\bm\la\neq\bm\mu$; in particular, $d^{\bm s}_{\bm\la\bm\mu}(v)= 0$ unless $\bm \mu \trianglerighteq \bm\la$.
\end{itemize}

{A lot of effort has been put in computing the canonical basis elements (i.e. computing the transition coefficients $d^{\bm s}_{\bm\la\bm\mu}(v)$), because of the following theorem. 

\begin{thm}\cite[Theorem 4.4]{A96}\label{arithm}
Let $\mathbb{F}$ be a field of characteristic $0$ and $\bm s \in I^r$ be a multicharge. Suppose $\bm\la$, $\bm\mu$ are $r$-multipartitions of $n$ with $\bm\mu$ dual Kleshchev. 
Then
$$[S'(\bm\la)\colon D'(\bm\mu)]=d_{\bm\la\bm\mu}^{\bm s}(1).$$
\end{thm}

This theorem, indeed, says that the coefficients $d^{\bm s}_{\bm\la\bm\mu}(v)$ specialised at $v = 1$ give the decomposition numbers of the corresponding cyclotomic Hecke algebras. In fact, thanks to the work of Brundan and Kleshchev \cite{BK09}, the coefficients $d^{\bm s}_{\bm\la\bm\mu}(v)$ (with $v$ still indeterminate) can be regarded as \textit{graded decomposition numbers} where `graded' is referred to the $\Z$-grading that Ariki-Koike algebras inherits by being isomorphic to certain quotients of KLR algebras. 
}

In \cite{BK09}, they also extend the bar involution on $M^{\bm s}$ to the whole of $\mathcal{F}^{\bm s}$. 
The extension of the bar involution to $\mathcal{F}^{\bm s}$ yields a canonical basis for the whole of $\mathcal{F}^{\bm s}$, indexed by the set of all $r$-multipartitions. In particular we get the following result.
\begin{thm}\label{decdom}\cite{Fay10}
For each multipartition $\bm\mu$, there is a unique vector
$$G^{\bm s}(\bm \mu) = \sum_{\bm\la\in\mathcal{P}^r}d^{\bm s}_{\bm\la\bm\mu}(v)\bm\la \in \mathcal{F}^{\bm s} \text{ with } d^{\bm s}_{\bm\la\bm\mu}(v)\in\mathbb{Q}(v)$$
such that
\begin{itemize}
\item
$\overline{G^{\bm s}(\bm\mu)}=G^{\bm s}(\bm\mu)$;
\item
$d^{\bm s}_{\bm\mu\bm\mu}(v)=1$, while $d^{\bm s}_{\bm\la\bm\mu}(v)\in v\mathbb{Z}[v]$ if $\bm\la\neq\bm\mu$;
\item
$d^{\bm s}_{\bm\la\bm\mu}(v)= 0$ unless $\bm \mu \trianglerighteq \bm\la$.
\end{itemize}
\end{thm}

In principle, the extension of the bar involution gives an algorithm for computing the canonical basis of $M^{\bm s}$. However, Fayers' algorithm works better for the purposes of this paper. So, we introduce and use Fayers' algorithm.

Fayers' approach is to compute the canonical basis for a module lying in between $M^{\bm s}$ and $\mathcal{F}^{\bm s}$. The way $\mathcal{F}^{\bm s}$ is defined and the choice of coproduct on $\UU$ mean that there is an isomorphism
\begin{alignat*}2
\mathcal{F}^{\bm s}&\overset{\sim}\longrightarrow\,\,&\mathcal{F}^{(s_1)}&\otimes\dots\otimes\mathcal{F}^{(s_r)}\\
\intertext{defined by linear extension of}
\bm\la&\longmapsto&(\la^{(1)})&\otimes\dots\otimes(\la^{(r)}).
\end{alignat*}
Thus, we will identify $\mathcal{F}^{\bm s}$ and $\mathcal{F}^{(s_1)}\otimes\dots\otimes\mathcal{F}^{(s_r)}$ via this isomorphism. Since each $\mathcal{F}^{(s_k)}$ contains a submodule $M^{(s_k)}$ isomorphic to $V(\Lambda_{s_k})$, $\mathcal{F}^{\bm s}$ contains a submodule $M^{\otimes\bm s}=M^{(s_1)}\otimes\dots\otimes M^{(s_r)}$ isomorphic to $V(\Lambda_{s_1})~\otimes~\dots~\otimes~V(\Lambda_{s_r})$.  This algorithm will compute the canonical basis of $M^{\otimes\bm s}$.

For presenting Fayers' LLT-type algorithm, we need the following result on canonical basis coefficients, since it will allow us to apply one of the steps of the algorithm. Recall that for any $r$-multipartition $\bm\la$ we define $\bm\la_{-}=(\la^{(2)}, \ldots, \la^{(r)})$; we also define $\bm s_- = (s_2, \ldots, s_r)$ for $\bm s \in I^r$.

\begin{cor}\label{fcomp}\cite[Corollary 3.2]{Fay10}
Suppose $\bm s\in I^r$ for $r>1$ and $\bm\mu\in\mathcal{P}^r$ with $\mu^{(1)}=\varnothing$.  If we write
\begin{align*}
G^{\bm s_-}(\bm\mu_-) &= \sum_{\bm\nu\in\mathcal{P}^{r-1}}d^{\bm s_-}_{\bm\nu\bm\mu_-}\bm\nu,\\
\intertext{then}
G^{\bm s}(\bm\mu) &= \sum_{\bm\nu\in\mathcal{P}^{r-1}}d^{\bm s_-}_{\bm\nu\bm\mu_-}\bm\nu_+.
\end{align*}
\end{cor}

\subsection{The LLT algorithm for $r=1$}
Before presenting the Fayers' LLT-type algorithm for Ariki-Koike algebras, we restrict attention to the case $r=1$, and explain the LLT algorithm for computing canonical basis elements $G^{(s_1)}(\mu)$. 
The LLT algorithm was first described in the paper \cite{LLT96}, to which we refer for more details and examples.


The canonical basis elements for $M^{(s_1)}$ are indexed by the $e$-regular partitions.  To construct $G^{(s_1)}(\mu)$ when $\mu$ is $e$-regular, we begin by constructing an auxiliary vector $A(\mu)$.  Let $l_1<\dots<l_t$ be the values of $l$ for which $\mathcal{L}_l(\mu)$ is non-empty.  For each $k$, let $m_k$ denote the number of nodes in $\mathcal{L}_{l_k}(\mu)$, and let $i_k$ denote the residue of $\mathcal{L}_{l_k}$.  Then the vector $A(\mu)$ is defined by
\[A(\mu) = f_{i_t}^{(m_t)}\dots f_{i_1}^{(m_1)}\cdot\varnothing.\]
$A(\mu)$ is obviously bar-invariant, and a lemma due to James \cite[6.3.54 \& 6.3.55]{JK81} implies that when we expand $A(\mu)$ as
\[A(\mu) = \sum_{\nu\in\mathcal{P}}a_\nu \nu,\]
we have $a_\mu=1$, while $a_\nu=0$ unless $\mu\trianglerighteq\nu$.  This means that $A(\mu)$ must equal $G^{(s_1)}(\mu)$ plus a $\mathbb{Q}(v+v^{-1})$-linear combination of canonical basis vectors $G^{(s_1)}(\nu)$ with $\mu\triangleright\nu$.  By induction on the dominance order we assume that these $G^{(s_1)}(\nu)$ have been computed, and so it is straightforward to subtract the appropriate multiples of these vectors from $A(\mu)$ to recover $G^{(s_1)}(\mu)$.  Moreover, the fact that the coefficients of the standard basis elements in $A(\mu)$ all lie in $\Z[v,v^{-1}]$ means that the coefficients of the canonical basis elements in $A(\mu)$ lie in $\Z[v+v^{-1}]$.  A more precise description of the procedure to strip off these canonical basis elements is given in Subsection \ref{algsec}.
\subsection{An LLT-type algorithm for $r\geq 1$}\label{algsec}

{Now, following \cite{Fay10} we give an algorithm for Ariki-Koike algebras which generalises the LLT algorithm for $r=1$.}  As mentioned above, this algorithm actually computes the canonical basis of $M^{\otimes \bm s}\cong M^{(s_1)}\otimes\dots\otimes M^{(s_r)}$.

Since the canonical basis elements $G^{(s_k)}(\mu)$ indexed by $e$-regular partitions $\mu$ form a basis for $M^{(s_k)}$, the tensor product $M^{(s_1)}\otimes\dots\otimes M^{(s_r)}$ has a basis consisting of all vectors $G^{(s_1)}(\mu^{(1)})\otimes\dots\otimes G^{(s_r)}(\mu^{(r)})$, where $\mu^{(1)},\dots,\mu^{(r)}$ are $e$-regular partitions. Translating this to the Fock space $\mathcal{F}^{\bm s}$, we find that $M^{\otimes \bm s}$ has a basis consisting of vectors
\[H^{\bm s}(\bm\mu) = \sum_{\bm\la\in\mathcal{P}^r}d^{(s_1)}_{\la^{(1)}\mu^{(1)}}\dots d^{(s_r)}_{\la^{(r)}\mu^{(r)}}\bm\la\]
for all $e$-multiregular multipartitions $\bm\mu$.  In fact, Fayers shows that 
\begin{prop}\cite[Proposition 4.2]{Fay07}\label{algo}
The canonical basis vectors $G^{\bm s}(\bm\mu)$ indexed by $e$-multiregular $r$-multipartitions $\bm\mu$ form a basis for the module $M^{\otimes\bm s}$.
\end{prop}
This implies in particular that the span of these vectors is a $\mathcal{U}$-submodule of $\mathcal{F}^{\bm s}$, which will enable the algorithm to work recursively.
Proposition \ref{algo} enables us to construct canonical basis vectors labelled by $e$-multiregular multipartitions recursively. As in the LLT algorithm, the idea is that to construct the canonical basis vector $G^{\bm s}(\bm\mu)$, we construct an auxiliary vector $A(\bm\mu)$ which is bar-invariant, and which we know equals $G^{\bm s}(\bm\mu)$ plus a linear combination of `lower' canonical basis vectors; the bar-invariance of $A(\bm\mu)$, together with dominance properties, allows these lower terms to be stripped off. Moreover, we have that $A(\bm\mu)$ lies in $M^{\otimes\bm s}$; then we know by Proposition \ref{algo} that all the canonical basis vectors occurring in $A(\bm\mu)$ are labelled by $e$-multiregular multipartitions, and therefore we can assume that these have already been constructed.

In fact, the proof of Proposition \ref{algo}, together with the LLT algorithm for partitions, gives us an LLT-type algorithm for Ariki-Koike algebras. We formalise this as follows.

First we define a partial order on multipartitions which is finer than the dominance order as follows: write $\bm\mu\succcurlyeq\bm\nu$ if either $|\mu^{(1)}|>|\nu^{(1)}|$ or $\mu^{(1)}\trianglerighteq\nu^{(1)}$.  Next we describe the steps of our recursive algorithm on the partial order $\succcurlyeq$. If $r>1$, we assume when computing $G^{\bm s}(\bm\mu)$ for $\bm\mu\in\mathcal{R}^r$ that we have already computed the vector $G^{\bm s_-}(\bm\mu_-)$, and that we have computed $G^{\bm s}(\bm\nu)$ for all $\bm\nu\in\mathcal{R}^r$ with $\bm\mu\succ\bm\nu$.

\begin{enumerate}
\item
If $\bm\mu=\bm\varnothing$, then $G^{\bm s}(\bm\mu)=\bm\varnothing$.
\item
If $\bm\mu\neq\bm\varnothing$ but $\mu^{(1)}=\varnothing$, then compute the canonical basis vector $G^{\bm s_-}(\bm\mu_-)$. Then $G^{\bm s}(\bm\mu)$ is given by
$$G^{\bm s}(\bm\mu) = \sum_{\bm\nu\in\mathcal{P}^{r-1}}d^{\bm s_-}_{\bm\nu\bm\mu_-}(v)\bm\nu_+.$$
\item
If $\mu^{(1)}\neq\varnothing$, then apply the following procedure.
\begin{enumerate}
\item
Let $\bm\mu_0=(\varnothing,\mu^{(2)},\dots,\mu^{(r)})$, and compute $G^{\bm s}(\bm\mu_0)$.
\item
Let $m_1,\ldots,m_t$ be the sizes of the non-empty ladders of $\mu^{(1)}$ in increasing order, and $i_1,\ldots,i_t$ be their residues.  Define $A(\bm \mu) = f_{i_t}^{(m_t)}\ldots f_{i_1}^{(m_1)}G^{\bm s}(\bm\mu_0)$. Write $A(\bm\mu)=\sum_{\bm\nu\in\mathcal{P}^r}a_{\bm\nu}\bm\nu$.  
\item\label{stepc}
If there is no $\bm\nu\neq\bm\mu$ for which $a_{\bm\nu}\notin v\mathbb{Z}[v]$, then stop.  Otherwise, take such a $\bm\nu$ which is maximal with respect to the dominance order, let $\alpha$ be the unique element of $\mathbb{Z}[v+v^{-1}]$ for which $a_{\bm\nu}-\alpha\in v\mathbb{Z}[v]$, replace $A(\bm\mu)$ by $A(\bm\mu)-\alpha G^{\bm s}(\bm\nu)$, and repeat. The remaining vector will be $G^{\bm s}(\bm\mu)$.
\end{enumerate}
\end{enumerate}

The vector $A(\bm\mu)$ computed in step 3 is a bar-invariant element of $M^{\bm s}$, because $G^{\bm s}(\bm\mu_0)$ is. Hence by Proposition \ref{algo} $A(\bm\mu)$ is a $\mathbb{Q}(v+v^{-1})$-linear combination of canonical basis vectors $G^{\bm s}(\bm\nu)$ with $\bm\nu\in\mathcal{R}^r$.  Furthermore, the rule for applying $f_i$ to a multipartition and the combinatorial results used in the LLT algorithm imply that $a_{\bm\mu}=1$, and that if $a_{\bm\la}\neq0$, then $\bm\mu\succcurlyeq\bm\la$.  
In particular, the partition $\bm\nu$ appearing in step 3(c) satisfies $\bm\mu\succ\bm\nu$; moreover, when $\alpha G^{\bm s}(\bm\nu)$ is subtracted from $A(\bm\mu)$, the condition that $a_{\bm\mu}=1$ and $a_{\bm\la}$ is non-zero only for $\bm\mu\succcurlyeq\bm\la$ remains true (because of Proposition \ref{decdom} and the fact that the order $\succcurlyeq$ refines the dominance order). So we can repeat, and complete step 3(c).

\begin{exe}
Let us take $e=r=2$, and write the set $I=\mathbb{Z}/2\mathbb{Z}$ as $\{0,1\}$.  Take $\bm s=(0,0)$. We want to compute the canonical basis element $G^{\bm s}(((1),(3,1)))$.

In the level $1$ Fock space $\mathcal{F}^{(0)}$, we need to compute $G^{(0)}((3,1))$. The non-empty ladders of the partition $(3,1)$ are $\mathcal{L}_1$, $\mathcal{L}_2$ and $\mathcal{L}_3$, of lengths $1,2$ and $1$ and residues $0,1,$ and $0$ respectively. So we have $$A((3,1))= f_0f_1^{(2)}f_0 \varnothing = (3,1)+v (2^2)+v^2 (2,1^2).$$ 
Since the coefficients in $A((3,1))$ (apart from the leading one) are divisible by $v$, we have $A((3,1))=G^{(0)}((3,1))$. Then
$$G^{\bm s}(\varnothing,(3,1))= (\varnothing,(3,1))+v (\varnothing,(2^2))+v^2 (\varnothing,(2,1^2)).$$
The only non-empty ladder of the partition $(1)$ is $\mathcal{L}_1$ of length $1$ and residue $0$. 
This time our auxiliary vector is 
\begin{align*}
    A((1),(3,1))=& f_0G^{\bm s}(\varnothing,(3,1))\\
    =& ((1),(3,1)) +v ((1),(2^2)) +v^2 ((1),(2,1^2)) +(v^2+1) (\varnothing, (3,2)) \\
    &+(v^3+v) (\varnothing, (3, 1^2)) + (v^4+v^2)(\varnothing,(2^2,1)).
\end{align*}
Notice that the coefficient of $(\varnothing, (3,2))$ is not a multiple of $v$.
So we compute $G^{\bm s}((\varnothing,(3,2)))$, that is:
$$G^{\bm s}((\varnothing,(3,2))) = (\varnothing,(3,2)) +v ((\varnothing,(3,1^2)) + v^2 (\varnothing,(2^2, 1)).$$
Hence, 
\begin{align*}
    G^{\bm s}((1),(3,1)) =& A((1),(3,1)) - G^{\bm s}((\varnothing,(3,2)))\\
    =& ((1),(3,1)) +v ((1),(2^2)) +v^2 ((1),(2,1^2)) +v^2 (\varnothing, (3,2)) \\
    &+v^3 (\varnothing, (3, 1^2)) + v^4(\varnothing,(2^2,1)).
\end{align*}
\end{exe}

\section{Addition of a full runner}\label{sec:addFullRun}
{Here, we recall the definition presented in \cite{Frunrem} of adding a runner `full' of beads to the abacus display of a partition. This will provide us the setting to generalise this definition to the case of a multipartition.}

\subsection{Truncated abacus configuration}\label{truncab}
To begin, consider the abacus configuration for a partition $\la=(\la_1, \ldots, \la_t)$.  The set $B_a(\la)$ of $\beta$-numbers used to define the abacus is an infinite set, and as such we have an infinite amount of beads in the abacus configuration, in particular there is a point where every row to the north of this point is completely full of beads. Instead we can consider a \textit{truncated} abacus configuration which has only finitely many beads on each runner, which we associate to a partition by filling in all the rows north of the highest beads with other beads. Conversely, if we are given a partition $\la$ we can fix a truncated abacus configuration associated to it. Let $N$ be an integer 
so that $x\in B_a(\la)$  whenever $x < Ne$. Then we define the truncated abacus configuration for $\la$ to be the one corresponding to the set $B_a(\la)\cap\{Ne,Ne+1,\ldots\}$. 
Notice that in terms of truncated abacus configuration, it makes sense to talk about the number of beads in the abacus. In particular, the truncated abacus configuration of $\la$ constructed in this way consists of $a-Ne$ beads. Indeed, 
$$B_a(\la)\cap\{Ne,Ne+1,\ldots\} = \{\la_1+a-1, \ldots, \la_t+a-t, 0+a-t-1, \ldots, Ne\}.$$
Hence, the number of beads of the truncated abacus configuration is equal to
$$|\{\la_1+a-1, \ldots, \la_t+a-t, 0+a-t-1, \ldots, Ne\}|,$$
that is $t+a-t-1-Ne+1=\la_1'+a-\la_1'-1-Ne+1=a-Ne$. {Thus, we will write ${\mathrm{Ab}_e(\la)}_{a-Ne}$ for the truncated $e$-abacus configuration of $\la$ with $a-Ne$ beads.}
\begin{exe}
Suppose $\lambda = (7, 4, 2^2)$, $e=5$ and $a = 0$. Then we take $N=-1$, so that 
$$B_0(\lambda) \cap \{-5,-4,-3, \ldots\}= \{6, 2,-1, -2,-5\}.$$
So, the truncated abacus display is
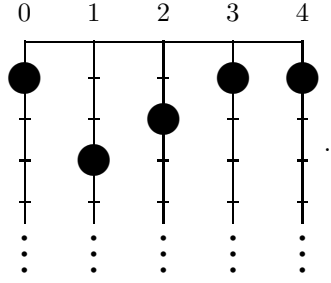
\[        \begin{matrix}

        0 & 1 & 2 & 3 & 4 \\

        \bigtl & \bigtm & \bigtm & \bigtm & \bigtr \\

        \bigbd & \bignb & \bignb & \bigbd & \bigbd \\

        \bignb & \bignb & \bigbd & \bignb & \bignb \\

        \bignb & \bigbd & \bignb & \bignb & \bignb \\
        
        \bignb & \bignb & \bignb & \bignb & \bignb \\

        \bigvd & \bigvd & \bigvd & \bigvd & \bigvd \\

        \end{matrix}.
\]
If we take $N=-3$, then we have
\begin{align*}
B_0(\la) &\cap \{-15, -14, -13,\ldots\}\\
&= \{6, 2,-1, -2,-5, -6,-7,-8,-9,-10,-11,-12,-13,-14,-15\}.
\end{align*}
This choice of $N$ gives the following truncated abacus configuration
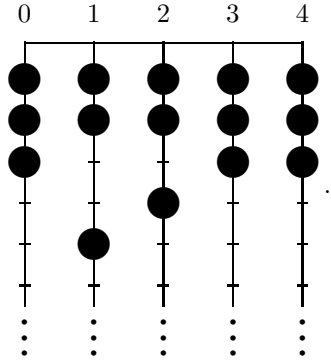
\[        \begin{matrix}

        0 & 1 & 2 & 3 & 4 \\

        \bigtl & \bigtm & \bigtm & \bigtm & \bigtr \\

        \bigbd & \bigbd & \bigbd & \bigbd & \bigbd \\

        \bigbd & \bigbd & \bigbd & \bigbd & \bigbd \\

        \bigbd & \bignb & \bignb & \bigbd & \bigbd \\

        \bignb & \bignb & \bigbd & \bignb & \bignb \\

        \bignb & \bigbd & \bignb & \bignb & \bignb \\
        
        \bignb & \bignb & \bignb & \bignb & \bignb \\

        \bigvd & \bigvd & \bigvd & \bigvd & \bigvd \\

        \end{matrix}.
\]
\end{exe}
Notice that for partitions the label of a runner does not correspond in general to the residue of the nodes represented by beads in that runner. Indeed, take $i, j \in I$, the beads on runner $j$ correspond to $i$-nodes of a partition $\la$ with $j\equiv i+b \mod e$ where $b$ is the number of beads of the truncated abacus configuration of $\la$.

Now, we want to have a closer look at the truncated $e$-abacus configuration of the empty partition $\varnothing$.
\begin{rmk}\label{emptyab}
Let $e\geq2$. Any truncated $e$-abacus configuration of $\varnothing$ has all the beads as high as possible with the runners from $0$ to $i$ consisting of $h+1$ beads and the runners from $i+1$ to $e-1$ consisting of $h$ beads for some $i\in I$ and some $h\geq 0$.
\end{rmk}

\subsection{Addition of a full runner for $r=1$ and empty partition}
Following \cite{Frunrem}, we define the addition of a `full' runner for the abacus display of a partition.

Given a partition $\la$ and a non-negative integer $k$, we construct a new partition $\la^{+k}$ as follows.
Let $a, N\in \mathbb{Z}$ such that $a \geq Ne$. Construct the truncated abacus configuration for $\la$ with $b:=a-Ne$ beads as in Section \ref{truncab}. Write $b+k=ce+d$, with $0 \leq d\leq e-1$, and add a runner to the abacus display immediately to the left of runner $d$; now put $c$ beads on this new runner, in the top $c$ positions, i.e. the position labelled $d, d+e+1, \ldots, d+(c-1)(e+1)$ in the usual labelling for an abacus with $e+1$ runners. The partition whose abacus display is obtained is $\la^{+k}$.
\begin{rmk}\label{residue_newrun}
The beads in the new inserted runner of $\la^{+k}$ correspond to nodes of residue $k$ mod $(e+1)$.
\end{rmk}
\begin{proof}
Construct the truncated $e$-abacus configuration for $\la$ with $b$ beads as above. Write $b+k=ce+d$ with $0 \leq d \leq e-1$. So the new inserted runner is labelled by $d$. We want to show that $d\equiv k+(b+c)$ mod $(e+1)$, since $b+c$ is the number of beads in the truncated $(e+1)$-abacus configuration of $\la^{+k}$. This is true because $k+b+c=ce+d+c=c(e+1)+d$.
\end{proof}

{We extend the operator $^{+k}$ linearly to the whole of the Fock space $\mathcal{F}$. We can now state the full runner removal theorem for the Iwahori-Hecke algebras of $\mathfrak{S}_n$ in terms of canonical bases.

\begin{thm}\cite[Theorem 3.1]{Frunrem}\label{Fullrunrem_level1} Suppose $\mu$ is an $e$-regular partition and $k \geq \mu_1$. Then $G_{e+1}(\mu^{+k}) = G_e(\mu)^{+k}$.
\end{thm}}

We want to focus our attention on the addition of a `full' runner for the empty partition $\varnothing$. 
\begin{rmk}
Let $k$ be a non-negative integer. Notice that $\varnothing^{+k}$ is an $(e+1)$-core. This follows by construction because in its abacus configuration all the beads are as high as possible.
\end{rmk}

\begin{prop}\label{empty+k}
Let $k$ be a non-negative integer.
\begin{enumerate}
\item If $k\in\{0, \ldots, e\}$, then $\varnothing^{+k}=\varnothing$.
\item Let $k>e$ and $k=k_1e+k_2$ with $k_1\geq1$ and $0 \leq k_2 \leq e-1$.
\begin{enumerate}
    \item If $k_2=0$, then $\varnothing^{+k}=((k_1-1)e, (k_1-2)e, \ldots, e)$.
    \item If $k_2\neq0$, then $\varnothing^{+k}=((k_1-1)e+k_2, (k_1-2)e+k_2,\ldots, k_2)$.
\end{enumerate}
\end{enumerate}
\end{prop}


\begin{proof}
Take $a$, $N \in \mathbb{Z}$ such that $a \geq Ne$ and construct the truncated abacus configuration of $\varnothing$ consisting of $b:=a-Ne$ beads. Let $i$ be the label of the runner of the last bead in $\mathrm{Ab}_e(\varnothing)$, we have $a-1=Me + i$ for some $M \in \mathbb{Z}$. Then $b=he+i+1$ with $h:=M-N$. By Remark \ref{emptyab}, the truncated abacus display of $\varnothing$ looks like the following:
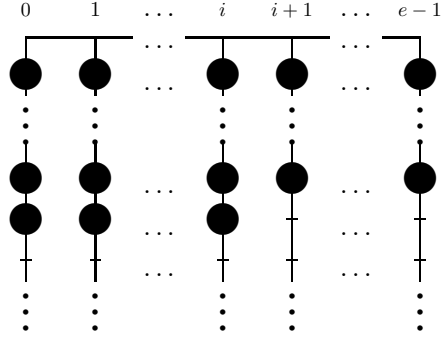
\[      \begin{matrix}

        \scalebox{0.75}{0} & \scalebox{0.8}{1} & \ldots & \scalebox{0.75}{$i$} & \scalebox{0.75}{$i+1$} & \ldots & \scalebox{0.75}{$e-1$}\\

        \bigtl & \bigtm & \ldots & \bigtm & \bigtm & \ldots & \bigtr \\
        \bigbd & \bigbd & \ldots & \bigbd & \bigbd & \ldots & \bigbd \\
        \bigvd & \bigvd &  & \bigvd & \bigvd &  & \bigvd \\
        \bigbd & \bigbd & \ldots & \bigbd & \bigbd & \ldots & \bigbd \\
        \bigbd & \bigbd & \ldots & \bigbd & \bignb & \ldots & \bignb \\
        \bignb & \bignb & \ldots & \bignb & \bignb & \ldots & \bignb  \\
        \bigvd & \bigvd &  & \bigvd & \bigvd &  & \bigvd \\

       \end{matrix}
\]
with $h+1$ beads in runners from $0$ to $i$ and $h$ beads in runners from $i+1$ to $e-1$. 
\begin{enumerate}
    \item If $k=\{0, \ldots, e\}$, then it is enough to show that the new inserted runner has either $h$ or $h+1$ beads because by Remark \ref{emptyab} the resulting truncated $(e+1)$-abacus configuration of $\varnothing^{+k}$ represents $\varnothing$.
    \begin{itemize}
        \item If $k\in \{0,\ldots,e-i-2\}$, then $b+k=he+i+1+k$ with $i+1+k \in \{i+1, \ldots, e-1\}$. Hence, the new inserted runner consists of $h$ beads and it is on the left of runner $i+1+k$. By Remark \ref{emptyab}, the runner $i$ of $\mathrm{Ab}_e(\varnothing)$ has $h+1$ beads, and all the runners of $\mathrm{Ab}_e(\varnothing)$ from $i+1$ to $e-1$ has exactly $h$ beads. So, $\mathrm{Ab}_{e+1}(\varnothing^{+k})$ represents $\varnothing$.
        \item If $k \in \{e-i-1, \ldots, e\}$, then $b+k=(h+1)e+i+1+k-e$ with $i+1+k-e\in \{0, \ldots, i+1\}$. Hence, the new inserted runner consists of $h+1$ beads and it is on the left of runner $i+1+k-e$. By Remark \ref{emptyab}, all the runners of $\mathrm{Ab}_e(\varnothing)$ from $0$ to $i$ has exactly $h+1$ beads. So, $\mathrm{Ab}_{e+1}(\varnothing^{+k})$ represents $\varnothing$.
        \end{itemize}
    \item If $k>e$, write $k=k_1e+k_2$ with $k_1\geq1$ and $0 \leq k_2 \leq e-1$.
    \begin{enumerate}
        \item If $k_2=0$ and $i\neq e-1$, then $b+k=(h+k_1)e+i+1$. Hence, we add the new runner to the left of runner $i+1$ with $h+k_1$ beads in the top positions. Reading off the partition from $\mathrm{Ab}_{e+1}(\varnothing^{+k})$ we get the partition $((k_1-1)e, (k_1-2)e, \ldots, e)$. Indeed, there are $k_1-1$ beads after the first empty position in $\mathrm{Ab}_{e+1}(\varnothing^{+k})$, and all of them are beads of the new inserted runner, so between two consecutive beads there are $e$ empty positions. \\
        If $k_2=0$ and $i=e-1$, then $b+k=(h+k_1+1)e$. Hence, we add the new runner to the left of runner $0$ with $h+k_1+1$ beads in the top positions. As above, reading off the partition from $\mathrm{Ab}_{e+1}(\varnothing^{+k})$ we get the partition $((k_1-1)e, (k_1-2)e, \ldots, e)$.
        \item If $k_2\neq0$ and $k_2\in \{0,\ldots,e-i-2\}$, then $b+k=(h+k_1)e+i+1+k_2$. Hence, we add the new runner to the left of runner $i+1+k_2$ with $h+k_1$ beads in the top positions. Reading off the partition from $\mathrm{Ab}_{e+1}(\varnothing^{+k})$ we get the partition $((k_1-1)e+k_2, (k_1-2)e+k_2, \ldots, k_2)$. Indeed, the first bead from the first empty position occurs after $(b+k)-b$ mod $e~= k_2$ empty spaces. Also, there are $k_1$ beads after the first empty position in $\mathrm{Ab}_{e+1}(\varnothing^{+k})$, and all of them are beads of the new inserted runner, so between two consecutive beads there are $e$ empty positions. \\
        If $k_2\neq0$ and $k_2\in \{e-i-1,\ldots,e-1\}$, then $b+k=(h+k_1+1)e+i+1+k_2-e$. As above, reading off the partition from $\mathrm{Ab}_{e+1}(\varnothing^{+k})$ we get the partition $((k_1-1)e+k_2, (k_1-2)e+k_2, \ldots, k_2)$.
    \end{enumerate}
\end{enumerate}
\end{proof}

\begin{exe}\label{exempty+k}
Consider the empty partition and $e=3$. Consider the abacus display of $\varnothing$ consisting of $b=7$ beads. 
\begin{enumerate}
\item Take $k=8=2\cdot3+2$. Then $\varnothing^{+8}$ has the following abacus configuration
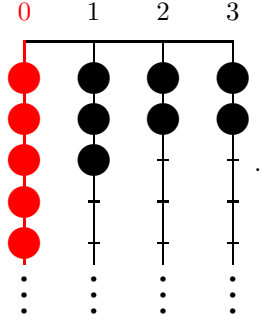
\[      \begin{matrix}

        \color{red}{0} & 1 & 2 & 3\\

        \bigtl & \bigtm & \bigtm & \bigtr \\

        \color{red}{\bigbd} & \bigbd & \bigbd & \bigbd \\

        \color{red}{\bigbd} & \bigbd & \bigbd & \bigbd \\

        \color{red}{\bigbd} & \bigbd & \bignb & \bignb \\
        
        \color{red}{\bigbd} & \bignb & \bignb & \bignb  \\
        
        \color{red}{\bigbd} & \bignb & \bignb & \bignb  \\

        \bigvd & \bigvd & \bigvd & \bigvd \\

        \end{matrix}.
\]
So, we have that $\varnothing^{+8} = (5,2)$.

\item Take $k=9=3\cdot3+0$. Then $\varnothing^{+9}$ has the following abacus configuration
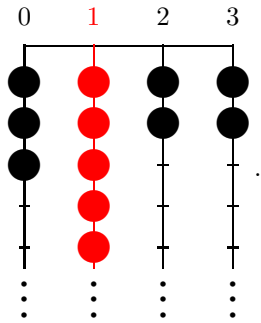
\[      \begin{matrix}

        0 & \color{red}{1} & 2 & 3\\

        \bigtl & \bigtm & \bigtm & \bigtr \\

        \bigbd & \color{red}{\bigbd} & \bigbd & \bigbd \\

        \bigbd & \color{red}{\bigbd} & \bigbd & \bigbd \\

        \bigbd & \color{red}{\bigbd} & \bignb & \bignb \\
        
        \bignb & \color{red}{\bigbd} & \bignb & \bignb  \\
        
        \bignb & \color{red}{\bigbd} & \bignb & \bignb  \\

        \bigvd & \bigvd & \bigvd & \bigvd \\

        \end{matrix}.
\]
So, we have that $\varnothing^{+9} = (6,3)$.
\item Take $k=10=3\cdot3+1$. Then $\varnothing^{+10}$ has the following abacus configuration
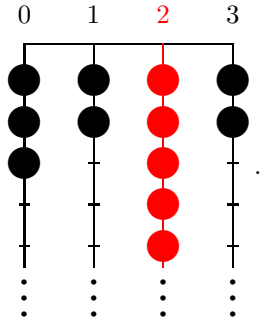
\[      \begin{matrix}

        0 & 1 & \color{red}{2} & 3\\

        \bigtl & \bigtm & \bigtm & \bigtr \\

        \bigbd & \bigbd & \color{red}{\bigbd} & \bigbd \\

        \bigbd & \bigbd & \color{red}{\bigbd} & \bigbd \\

        \bigbd & \bignb & \color{red}{\bigbd} & \bignb \\
        
        \bignb & \bignb & \color{red}{\bigbd} & \bignb  \\
        
        \bignb & \bignb & \color{red}{\bigbd} & \bignb  \\

        \bigvd & \bigvd & \bigvd & \bigvd \\

        \end{matrix}.
\]
So, we have that $\varnothing^{+10} = (7,4,1)$.
\end{enumerate}
\end{exe}

\begin{rmk}\label{remnoderes}
All the removable nodes of $\varnothing^{+k}$ have the same residue because if there are any removable nodes in $\varnothing^{+k}$, they are the nodes corresponding to the beads in the new inserted runner and they have all the same residue $k\mod(e+1)$ as shown in Remark \ref{residue_newrun}.
\end{rmk}

This remark is helpful because it allows us to find an induction sequence from $\varnothing$ to $\varnothing^{+k}$, as shown in the following.
%
From now on, given an integer $j$ we denote by $\overline{j}$ the residue of $j$ modulo $e+1$ and by $\mathfrak{F}_{\overline{i}}$ for $i \in \Z$ the generator $f_{\overline{i}}$ of the quantised enveloping algebra $U_v(\widehat{\mathfrak{sl}}_{e+1})$.
For $\ell\geq 1$ and $i\in \mathbb{Z}$, set
\begin{equation}
    \mathcal{G}^{(\ell)}_{\overline{i}}:= \mathfrak{F}^{(\ell)}_{\overline{i}} \mathfrak{F}^{(\ell)}_{\overline{i-1}} \ldots \mathfrak{F}^{(\ell)}_{\overline{i-(e-1)}}.
\end{equation}
\begin{prop}\label{indseq}
Let $k$ be a non-negative integer. Write $k=k_1e+k_2$ with $k_1\geq0$ and $0 \leq k_2 \leq e-1$. Then
    $$\mathfrak{F}^{(k_1)}_{\overline{k}} \mathfrak{F}^{(k_1)}_{\overline{k-1}}\ldots\mathfrak{F}^{(k_1)}_{\overline{k-k_2+1}}\mathcal{G}^{(k_1-1)}_{\overline{k-k_2}}\mathcal{G}^{(k_1-2)}_{\overline{k-k_2+1}} \ldots \mathcal{G}^{(1)}_{\overline{k-k_2+k_1-2}}(\varnothing)=\varnothing^{+k}.$$
where $\mathfrak{F}^{(k_1)}_{\overline{k}} \mathfrak{F}^{(k_1)}_{\overline{k-1}}\ldots\mathfrak{F}^{(k_1)}_{\overline{k-k_2+1}}$ only occur if $k_2\neq0$.
\end{prop}

\begin{proof}
We prove this by induction on $k$.
\begin{itemize}
    \item If $k<e$, then $\varnothing^{+k}=\varnothing$ by Proposition \ref{empty+k}, so there is nothing to prove.
    \item Suppose $k\geq e$ and write $k=k_1e+k_2$ with $k_1\geq1$ and $0 \leq k_2 \leq e-1$. Then by Lemma 3.3 in \cite{Frunrem} we have $\varnothing^{+k} = \mathfrak{F}_{\bar{k}}^{(m)}(\varnothing^{+(k-1)})$ for some constant $m$. In particular, by the proof of Proposition \ref{empty+k} we know that
    \begin{itemize}
        \item if $k_2=0$, then $m=k_1-1$;
        \item if $k_2\neq 0$, then $m=k_1$.
    \end{itemize}
    We study the two cases separately.
    \begin{itemize}
        \item If $k_2=0$, then $k-1=(k_1-1)e+e-1$ and so by the induction hypothesis we have 
         \begin{align*}
         \begin{split}      
         \varnothing^{+(k-1)}&=\mathfrak{F}_{\overline{k-1}}^{(k_1-1)} \mathfrak{F}_{\overline{k-2}}^{(k_1-1)}\cdots \mathfrak{F}_{\overline{k-1-e+1+1}}^{(k_1-1)} \mathcal{G}_{\overline{k-1-e+1}}^{(k_1-1-1)} \\
            & \qquad \mathcal{G}_{\overline{k-1-e+1+1}}^{(k_1-1-2)} \ldots \mathcal{G}_{\overline{k-1-e+1+k_1-1-2}}^{(1)}(\varnothing)\end{split}\\            
            &=\mathfrak{F}_{\overline{k-1}}^{(k_1-1)} \mathfrak{F}_{\overline{k-2}}^{(k_1-1)}\cdots \mathfrak{F}_{\overline{k-e+1}}^{(k_1-1)} \mathcal{G}_{\overline{k+1}}^{(k_1-2)} \mathcal{G}_{\overline{k+2}}^{(k_1-3)} \ldots \mathcal{G}_{\overline{k+k_1-2}}^{(1)}(\varnothing).
         \end{align*}
         Hence, 
         \begin{align*}
         \varnothing^{+k} &= \mathfrak{F}_{\overline{k}}^{(k_1-1)}\mathfrak{F}_{\overline{k-1}}^{(k_1-1)} \mathfrak{F}_{\overline{k-2}}^{(k_1-1)}\cdots \mathfrak{F}_{\overline{k-e+1}}^{(k_1-1)} \mathcal{G}_{\overline{k+1}}^{(k_1-2)} \mathcal{G}_{\overline{k+2}}^{(k_1-3)} \ldots \mathcal{G}_{\overline{k+k_1-2}}^{(1)}(\varnothing)\\
         &= \mathcal{G}_{\overline{k}}^{(k_1-1)}\mathcal{G}_{\overline{k+1}}^{(k_1-2)} \mathcal{G}_{\overline{k+2}}^{(k_1-3)} \ldots \mathcal{G}_{\overline{k+k_1-2}}^{(1)}(\varnothing).
         \end{align*}
          \item If $k_2\neq0$, then $k-1=k_1e+k_2-1$ and so by induction hypothesis we have 
         \begin{align*}
         \begin{split}
            \varnothing^{+(k-1)}&=\mathfrak{F}_{\overline{k-1}}^{(k_1)} \mathfrak{F}_{\overline{k-1-1}}^{(k_1)}\cdots \mathfrak{F}_{\overline{k-1-k_2+1+1}}^{(k_1)} \mathcal{G}_{\overline{k-1-k_2+1}}^{(k_1-1)} \\
            & \qquad \mathcal{G}_{\overline{k-1-k_2+1+1}}^{(k_1-2)} \ldots \mathcal{G}_{\overline{k-1-k_2+1+k_1-2}}^{(1)}(\varnothing)
            \end{split}\\
            & =\mathfrak{F}_{\overline{k-1}}^{(k_1)} \mathfrak{F}_{\overline{k-2}}^{(k_1)}\cdots \mathfrak{F}_{\overline{k-k_2+1}}^{(k_1)} \mathcal{G}_{\overline{k-k_2}}^{(k_1-1)} \mathcal{G}_{\overline{k-k_2+1}}^{(k_1-2)} \ldots \mathcal{G}_{\overline{k-k_2+k_1-2}}^{(1)}(\varnothing).
         \end{align*}
         Hence, 
       $$\varnothing^{+k} = \mathfrak{F}_{\overline{k}}^{(k_1)}\mathfrak{F}_{\overline{k-1}}^{(k_1)} \mathfrak{F}_{\overline{k-2}}^{(k_1)}\cdots \mathfrak{F}_{\overline{k-k_2+1}}^{(k_1)} \mathcal{G}_{\overline{k-k_2}}^{(k_1-1)} \mathcal{G}_{\overline{k-k_2+1}}^{(k_1-2)} \ldots \mathcal{G}_{\overline{k-k_2+k_1-2}}^{(1)}(\varnothing).$$
    \end{itemize}
\end{itemize}
\end{proof}

We give some examples of the induction sequence we want to deal with.
\begin{exe}
With the same notation and choices of Example \ref{exempty+k}, the induction sequence given in Proposition \ref{indseq} from $\varnothing$ to $\varnothing^{+k}$ acts as follows. Notice that at each step we apply the induction operators starting from the rightmost one. 
\begin{itemize}
    \item[2.] For $k=9$, in terms of abacus configuration we have\\
\scalebox{0.85}
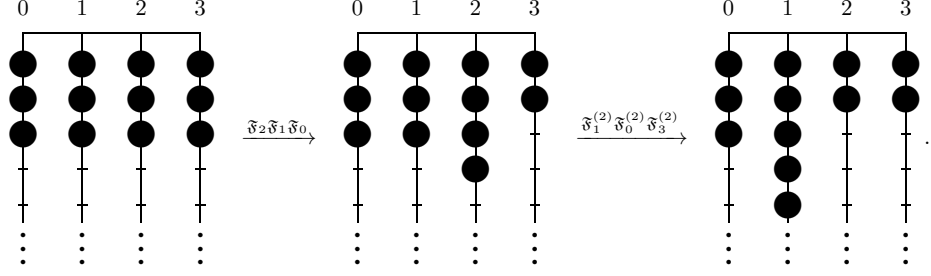
{
  \begin{tabular}{ccccc}
   $\begin{matrix}

        0 & 1 & 2 & 3\\

        \bigtl & \bigtm & \bigtm & \bigtr \\

        \bigbd & \bigbd & \bigbd & \bigbd \\

        \bigbd & \bigbd & \bigbd & \bigbd \\

        \bigbd & \bigbd & \bigbd & \bigbd \\
        
        \bignb & \bignb & \bignb & \bignb  \\
        
        \bignb & \bignb & \bignb & \bignb  \\

        \bigvd & \bigvd & \bigvd & \bigvd \\

        \end{matrix}$
        & $\xrightarrow{\mathfrak{F}_2\mathfrak{F}_1\mathfrak{F}_0}$
        &
        $\begin{matrix}

        0 & 1 & 2 & 3\\

        \bigtl & \bigtm & \bigtm & \bigtr \\

        \bigbd & \bigbd & \bigbd & \bigbd \\

        \bigbd & \bigbd & \bigbd & \bigbd \\

        \bigbd & \bigbd & \bigbd & \bignb \\
        
        \bignb & \bignb & \bigbd & \bignb  \\
        
        \bignb & \bignb & \bignb & \bignb  \\

        \bigvd & \bigvd & \bigvd & \bigvd \\
        \end{matrix}
        $
        & $\xrightarrow{\mathfrak{F}_1^{(2)}\mathfrak{F}_0^{(2)}\mathfrak{F}_3^{(2)}}$
        &
        $\begin{matrix}

        0 & 1 & 2 & 3\\

        \bigtl & \bigtm & \bigtm & \bigtr \\

        \bigbd & \bigbd & \bigbd & \bigbd \\

        \bigbd & \bigbd & \bigbd & \bigbd \\

        \bigbd & \bigbd & \bignb & \bignb \\
        
        \bignb & \bigbd & \bignb & \bignb  \\
        
        \bignb & \bigbd & \bignb & \bignb  \\

        \bigvd & \bigvd & \bigvd & \bigvd \\
        \end{matrix}
        $.
  \end{tabular}}

\item[3.] For $k=10$, in terms of abacus configuration we have\\
\scalebox{0.85}
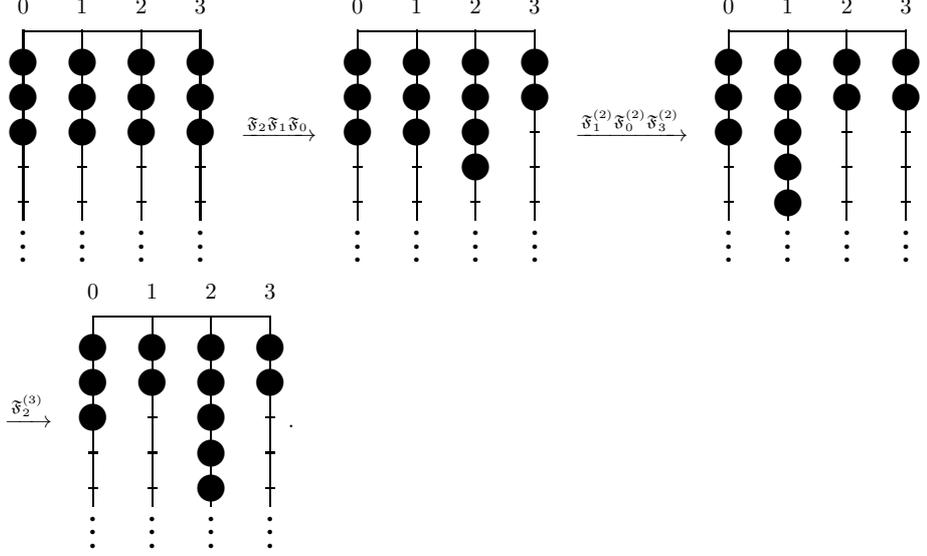
{
  \begin{tabular}{ccccc}
   $\begin{matrix}

        0 & 1 & 2 & 3\\

        \bigtl & \bigtm & \bigtm & \bigtr \\

        \bigbd & \bigbd & \bigbd & \bigbd \\

        \bigbd & \bigbd & \bigbd & \bigbd \\

        \bigbd & \bigbd & \bigbd & \bigbd \\
        
        \bignb & \bignb & \bignb & \bignb  \\
        
        \bignb & \bignb & \bignb & \bignb  \\

        \bigvd & \bigvd & \bigvd & \bigvd \\

        \end{matrix}$
        & $\xrightarrow{\mathfrak{F}_2\mathfrak{F}_1\mathfrak{F}_0}$
        &
        $\begin{matrix}

        0 & 1 & 2 & 3\\

        \bigtl & \bigtm & \bigtm & \bigtr \\

        \bigbd & \bigbd & \bigbd & \bigbd \\

        \bigbd & \bigbd & \bigbd & \bigbd \\

        \bigbd & \bigbd & \bigbd & \bignb \\
        
        \bignb & \bignb & \bigbd & \bignb  \\
        
        \bignb & \bignb & \bignb & \bignb  \\

        \bigvd & \bigvd & \bigvd & \bigvd \\
        \end{matrix}
        $
        & $\xrightarrow{\mathfrak{F}_1^{(2)}\mathfrak{F}_0^{(2)}\mathfrak{F}_3^{(2)}}$
        &
        $\begin{matrix}

        0 & 1 & 2 & 3\\

        \bigtl & \bigtm & \bigtm & \bigtr \\

        \bigbd & \bigbd & \bigbd & \bigbd \\

        \bigbd & \bigbd & \bigbd & \bigbd \\

        \bigbd & \bigbd & \bignb & \bignb \\
        
        \bignb & \bigbd & \bignb & \bignb  \\
        
        \bignb & \bigbd & \bignb & \bignb  \\

        \bigvd & \bigvd & \bigvd & \bigvd \\
        \end{matrix}
        $
  \end{tabular}}\\
  \scalebox{0.85}{
     \begin{tabular}{cc}
        $\xrightarrow{\mathfrak{F}_2^{(3)}}$  &         $\begin{matrix}

        0 & 1 & 2 & 3\\

        \bigtl & \bigtm & \bigtm & \bigtr \\

        \bigbd & \bigbd & \bigbd & \bigbd \\

        \bigbd & \bigbd & \bigbd & \bigbd \\

        \bigbd & \bignb & \bigbd & \bignb \\
        
        \bignb & \bignb & \bigbd & \bignb  \\
        
        \bignb & \bignb & \bigbd & \bignb  \\

        \bigvd & \bigvd & \bigvd & \bigvd \\
        \end{matrix}
        .$
     \end{tabular}}
\end{itemize}
\end{exe}

Another useful result will be the next lemma that gives a way to establish when a sequence of induction operators $\mathfrak{F}^{(m)}_{\bar{i}}$ for $m,i\in \Z$ with $m\geq 0$ acts non-zero.
\begin{lem}\label{seqindnot0}
Let $k\geq e+1$. Write $k=k_1e+k_2$ with $k_1\geq1$ and $0 \leq k_2 \leq e-1$. Consider
\begin{equation}\label{action0}
\mathfrak{F}_{\overline{k+e+1}}^{(q_{e+1})}\mathfrak{F}_{\overline{k+e}}^{(q_{e})} \ldots \mathfrak{F}_{\overline{k+2}}^{(q_{2})}\mathfrak{F}_{\overline{k+1}}^{(q_{1})} ( \varnothing^{+(k-e-1)} ),
\end{equation}

for some $q_j\geq 0$ for all $1 \leq j\leq e+1$. Then $\eqref{action0}\neq0$ if and only if 
\begin{enumerate}
\item $k_1-1 \geq q_1\geq q_2 \geq \ldots \geq q_{e+1-k_2}$, and $q_{e+2-k_2}\geq q_{e+3-k_2} \geq \ldots \geq q_{e+1}$, and,
    \item if $q_{e+1-k_2}=k_1-1$, then  $q_{e+1-k_2}\geq q_{e+2-k_2}-1$;\\
    if $q_{e+1-k_2}<k_1-1$, then  $q_{e+1-k_2}\geq q_{e+2-k_2}$.
\end{enumerate}
\end{lem}

\begin{proof}
We want to prove that $\eqref{action0}\neq 0$ if and only if conditions 1. and 2. hold. This is equivalent to proving that 
\begin{equation}\label{equivtoprove}
    \mathfrak{F}^{(q_j)}_{\overline{k+j}} \ldots \mathfrak{F}^{(q_1)}_{\overline{k+1}}( \varnothing^{+(k-e-1)} )\neq0 \text{ for all }1\leq j \leq e+1
\end{equation}
if and only if conditions 1. and 2. hold.
Recall that if $s$ is the runner corresponding to nodes of residue $\overline{k+j}$, then $\mathfrak{F}_{\overline{k+j}}^{(q_j)}$ moves $q_j$ beads from runner $s-1$ to runner $s$. Moreover, by Proposition \ref{empty+k}, the addable nodes of residue $\overline{k+1}$ of $\varnothing^{+(k-e-1)}$ are the nodes corresponding to beads in the new inserted runner and there are $k_1-1$ of them.

 In particular, given the abacus configuration of $\varnothing^{+(k-e-1)}$, for each $j \neq e+2-k_2$ the number of addable nodes of residue $\overline{k+j}$ of any term in $\mathfrak{F}^{(q_{j-1})}_{\overline{k+j-1}} \ldots \mathfrak{F}^{(q_1)}_{\overline{k+1}}( \varnothing^{+(k-e-1)} )$ is $q_{j-1}$. For $j=e+2-k_2$, the number of addable nodes of residue $\overline{k+j}=\overline{k+e+2-k_2}=\overline{k-k_2+1}$ of any term in $\mathfrak{F}^{(q_{e+1-k_2})}_{\overline{k-k_2}} \ldots \mathfrak{F}^{(q_1)}_{\overline{k+1}}( \varnothing^{+(k-e-1)} )$ is 
\begin{itemize}
    \item $q_{e+1-k_2}+1$ if $q_{e+1-k_2}=k_1-1$,
    \item $q_{e+1-k_2}$ if $q_{e+1-k_2}<k_1-1$.
\end{itemize}
We proceed by induction on $j$.
\begin{itemize}
\item If $j=1$, then $\mathfrak{F}^{(q_1)}_{\overline{k+1}}( \varnothing^{+(k-e-1)} ) \neq 0$ if and only if $q_1\leq k_1-1$ because the number of addable nodes of residue $\overline{k+1}$ of $\varnothing^{+(k-e-1)}$ is $k_1-1$.
\item If $j>1$, then by induction hypothesis we know 
$$\mathfrak{F}_{j-1}( \varnothing^{+(k-e-1)} ):=\mathfrak{F}^{(q_{j-1})}_{\overline{k+j-1}} \ldots \mathfrak{F}^{(q_1)}_{\overline{k+1}}( \varnothing^{+(k-e-1)} )\neq0$$
if and only if conditions 1. and 2. hold for $q_1, \ldots, q_{j-1}$.
We want to prove that
$$\mathfrak{F}^{(q_{j})}_{\overline{k+j}}\mathfrak{F}_{j-1}( \varnothing^{+(k-e-1)} )\neq0$$
if and only if the conditions 1. and 2. hold also for $q_j$.
\begin{itemize}
    \item If $j\leq e+1-k_2$, then the number of addable nodes of residue $\overline{k+j}$ of any term in $\mathfrak{F}_{j-1}( \varnothing^{+(k-e-1)} )$ is $q_{j-1}$. So, $\mathfrak{F}^{(q_{j})}_{\overline{k+j}}\mathfrak{F}_{j-1}( \varnothing^{+(k-e-1)} )\neq0$ if and only if $q_{j-1}\geq q_j$.
    \item If $j=e+2-k_2$, then the number of addable nodes of residue $\overline{k+j}=\overline{k-k_2+1}$ of any term in $\mathfrak{F}_{e+1-k_2}( \varnothing^{+(k-e-1)} )$ is
    \begin{center}
    $\begin{cases}
    q_{e+1-k_2}+1 & \text{if } q_{e+1-k_2}=k_1-1,\\
    q_{e+1-k_2} & \text{if } q_{e+1-k_2}<k_1-1.
    \end{cases}$  
    \end{center}
    So, if $q_{e+1-k_2}=k_1-1$ then $\mathfrak{F}^{(q_{e+2-k_2})}_{\overline{k-k_2+1}}\mathfrak{F}_{e+1-k_2}( \varnothing^{+(k-e-1)} )\neq0$ if and only if $q_{e+1-k_2}\geq q_{e+2-k_2}-1$. If $q_{e+1-k_2}<k_1-1$ then $\mathfrak{F}^{(q_{e+2-k_2})}_{\overline{k-k_2+1}}\mathfrak{F}_{e+1-k_2}( \varnothing^{+(k-e-1)} )\neq0$ if and only if $q_{e+1-k_2}\geq q_{e+2-k_2}$.
    \item If $j>e+2-k_2$, then the number of addable nodes of residue $\overline{k+j}$ of any term of $\mathfrak{F}_{j-1}( \varnothing^{+(k-e-1)} )$ is $q_{j-1}$. So, $\mathfrak{F}^{(q_{j})}_{\overline{k+j}}\mathfrak{F}_{j-1}( \varnothing^{+(k-e-1)} )\neq0$ if and only if $q_{j-1}\geq q_j$.
    \end{itemize}
\end{itemize}

\end{proof}

\begin{cor}\label{corseqindnot0}
In the same notation of Lemma \ref{seqindnot0}, if $\eqref{action0}\neq 0$ then it holds that
\begin{itemize}
    \item if $q_{e+1-k_2}<k_1-1$, then $q_{e+1}\leq q_1$;
    \item if $q_{e+1-k_2}=k_1-1$, then $q_{e+1}\leq q_1+1$.
\end{itemize}
\end{cor}

\begin{proof}
It follows directly from writing the inequalities of conditions 1. and 2. of Lemma \ref{seqindnot0} one next to the other for the two cases.
\end{proof}

\subsection{Addition of a full runner for $r\geq2$}\label{def+k}

{Let $\bm\la$ be a $r$-multipartition of $n$ and $\kappa = (\kappa_1, \ldots, \kappa_r)$ be a multicharge for $\mathcal{H}_{r,n}$. 
For each $j \in \{1, \ldots r\}$, we represent each component $\la^{(j)}$ with a truncated abacus consisting of $n_j$ beads such that $n_j\equiv\kappa_j \mod e$.}

Given  $\bm\la$ as above, we construct a new $r$-multipartition as follows. Let $0~\leq~d~\leq~e-1$. For each $j \in \{1, \ldots, r\}$,
\begin{itemize}
\item take $n_j$ determined as above and construct the abacus display for $\la^{(j)}$ with $n_j$ beads;
\item set $k^{(j)}$ a non-negative integer such that 
$n_j+k^{(j)} \equiv d \mod e$  for all $j\in \{1, \ldots, r\}$;
\item write $n_j + k^{(j)} = c_je+d$ for all $j\in \{1, \ldots, r\}$;
\item add a runner to each component of the abacus display immediately to the left of runner $d$;
\item for each $j$, put $c_j$ beads on the new inserted runner of each component $j$, in the top $c_j$ positions, i.e. the positions labelled by $d, d+e+1, \ldots, d + (c_j-1)(e+1)$ in the usual labelling for an abacus with $e+1$ runners.
\end{itemize}
The $r$-multipartition whose abacus is obtained is $\bm\la^{+\bm k}:=\bm\la^{+ (k^{(1)}, \ldots, k^{(r)})}.$

\begin{rmk}
Notice that choosing the number of beads $(n_1, \ldots, n_r)$ in order to construct the abacus display of a multipartition $\bm \la$ corresponds to the choice of the multicharge $\bm a = (n_1, \ldots, n_r)$.
\end{rmk}

\begin{exe}\label{exrunnadd}
Suppose $\bm\la = ((4,3,2),(2^2),(3))$ and $e=4$. Choose $\bm n=(n_1,n_2,n_3)=(11,9,12)$, that corresponds to choose $\bm n$ as multicharge. 
So we get the following abacus display:
$$
        \begin{matrix}

        0 & 1 & 2 & 3\\

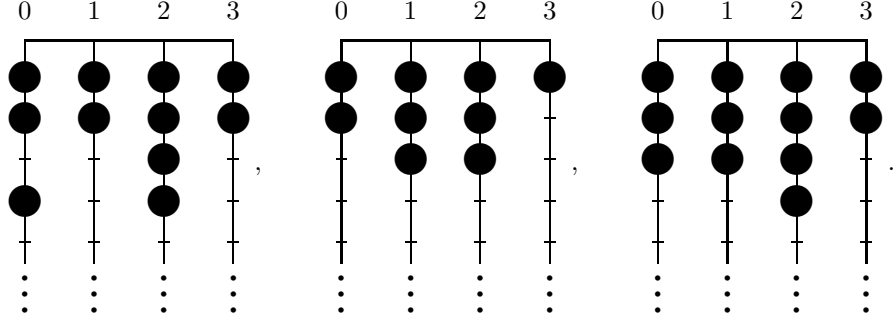
        \bigtl & \bigtm & \bigtm & \bigtr\\
        \bigbd & \bigbd & \bigbd & \bigbd\\
        \bigbd & \bigbd & \bigbd & \bigbd\\
        \bignb & \bignb & \bigbd & \bignb\\
        \bigbd & \bignb & \bigbd & \bignb\\
        \bignb & \bignb & \bignb & \bignb\\
        \bigvd & \bigvd & \bigvd & \bigvd

        \end{matrix}, \qquad
        \begin{matrix}
        
        0 & 1 & 2 & 3\\

        \bigtl & \bigtm & \bigtm & \bigtr\\
        \bigbd & \bigbd & \bigbd & \bigbd\\
        \bigbd & \bigbd & \bigbd & \bignb\\
        \bignb & \bigbd & \bigbd & \bignb\\
        \bignb & \bignb & \bignb & \bignb\\
        \bignb & \bignb & \bignb & \bignb\\
        \bigvd & \bigvd & \bigvd & \bigvd

        \end{matrix}, \qquad
        \begin{matrix}

        0 & 1 & 2 & 3\\

        \bigtl & \bigtm & \bigtm & \bigtr\\
        \bigbd & \bigbd & \bigbd & \bigbd\\
        \bigbd & \bigbd & \bigbd & \bigbd\\
        \bigbd & \bigbd & \bigbd & \bignb\\
        \bignb & \bignb & \bigbd & \bignb\\
        \bignb & \bignb & \bignb & \bignb\\
        \bigvd & \bigvd & \bigvd & \bigvd

        \end{matrix}.$$

Then,
\begin{itemize}

\item if $(k^{(1)},k^{(2)}, k^{(3)})=(5,3,0)$, we obtain $\bm\la^{+(5,3,0)}$ with abacus configuration\\
\scalebox{0.9}
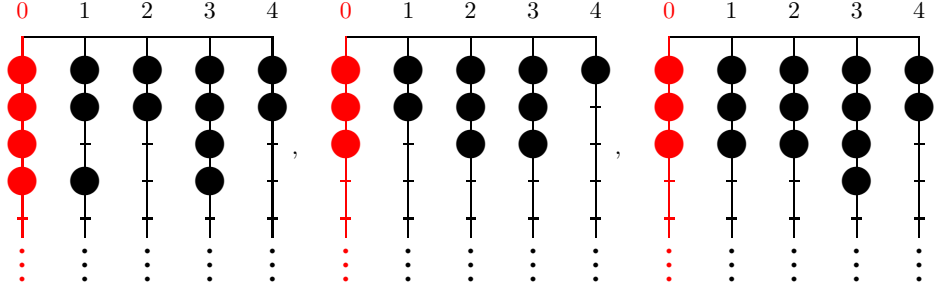
{
\begin{tabular}{c}
$\begin{matrix}

        \color{red}{0} & 1 & 2 & 3 & 4\\

        \bigtl & \bigtm & \bigtm & \bigtm & \bigtr\\
        \color{red}{\bigbd} & \bigbd & \bigbd & \bigbd & \bigbd\\
        \color{red}{\bigbd} & \bigbd & \bigbd & \bigbd & \bigbd\\
        \color{red}{\bigbd} & \bignb & \bignb & \bigbd & \bignb\\
        \color{red}{\bigbd} & \bigbd & \bignb & \bigbd & \bignb\\
        \color{red}{\bignb} & \bignb & \bignb & \bignb & \bignb\\
        \color{red}{\bigvd} & \bigvd &\bigvd & \bigvd & \bigvd

        \end{matrix}, \quad
        \begin{matrix}
        
        \color{red}{0} & 1 & 2 & 3 & 4\\

        \bigtl & \bigtm & \bigtm & \bigtm & \bigtr\\
        \color{red}{\bigbd} & \bigbd & \bigbd & \bigbd & \bigbd\\
        \color{red}{\bigbd} & \bigbd & \bigbd & \bigbd & \bignb\\
        \color{red}{\bigbd} & \bignb & \bigbd & \bigbd & \bignb\\
        \color{red}{\bignb} & \bignb & \bignb & \bignb & \bignb\\
        \color{red}{\bignb} & \bignb & \bignb & \bignb & \bignb\\
        \color{red}{\bigvd} & \bigvd &\bigvd & \bigvd & \bigvd

        \end{matrix}, \quad
        \begin{matrix}

        \color{red}{0} & 1 & 2 & 3 & 4\\

        \bigtl & \bigtm & \bigtm & \bigtm & \bigtr\\
        \color{red}{\bigbd} & \bigbd & \bigbd & \bigbd & \bigbd\\
        \color{red}{\bigbd} & \bigbd & \bigbd & \bigbd & \bigbd\\
        \color{red}{\bigbd} & \bigbd & \bigbd & \bigbd & \bignb\\
        \color{red}{\bignb} & \bignb & \bignb & \bigbd & \bignb\\
        \color{red}{\bignb} & \bignb & \bignb & \bignb & \bignb\\
        \color{red}{\bigvd} & \bigvd & \bigvd & \bigvd & \bigvd

        \end{matrix}$
\end{tabular}}\\
         and multicharge equals to $(15,12,15)$;

\item if $(k^{(1)},k^{(2)}, k^{(3)})=(3,9,6)$, we obtain $\bm\la^{+(3,9,6)}$ with abacus configuration\\
\scalebox{0.9}
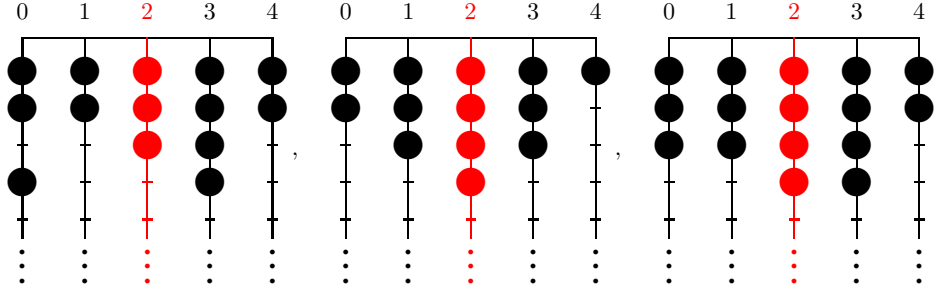
{
\begin{tabular}{c}$\begin{matrix}

        0 & 1 & \color{red}{2} & 3 & 4\\

        \bigtl & \bigtm & \bigtm & \bigtm & \bigtr\\
        \bigbd & \bigbd & \color{red}{\bigbd} & \bigbd & \bigbd\\
        \bigbd & \bigbd & \color{red}{\bigbd} & \bigbd & \bigbd\\
        \bignb & \bignb & \color{red}{\bigbd} & \bigbd & \bignb\\
        \bigbd & \bignb & \color{red}{\bignb} & \bigbd & \bignb\\
        \bignb & \bignb & \color{red}{\bignb} & \bignb & \bignb\\
        \bigvd & \bigvd & \color{red}{\bigvd} & \bigvd & \bigvd

        \end{matrix}, \quad
        \begin{matrix}
        
        0 & 1 & \color{red}{2} & 3 & 4\\

        \bigtl & \bigtm & \bigtm & \bigtm & \bigtr\\
        \bigbd & \bigbd & \color{red}{\bigbd} & \bigbd & \bigbd\\
        \bigbd & \bigbd & \color{red}{\bigbd} & \bigbd & \bignb\\
        \bignb & \bigbd & \color{red}{\bigbd} & \bigbd & \bignb\\
        \bignb & \bignb & \color{red}{\bigbd} & \bignb & \bignb\\
        \bignb & \bignb & \color{red}{\bignb} & \bignb & \bignb\\
        \bigvd & \bigvd & \color{red}{\bigvd} & \bigvd & \bigvd

        \end{matrix}, \quad
        \begin{matrix}

        0 & 1 & \color{red}{2} & 3 & 4\\

        \bigtl & \bigtm & \bigtm & \bigtm & \bigtr\\
        \bigbd & \bigbd & \color{red}{\bigbd} & \bigbd & \bigbd\\
        \bigbd & \bigbd & \color{red}{\bigbd} & \bigbd & \bigbd\\
        \bigbd & \bigbd & \color{red}{\bigbd} & \bigbd & \bignb\\
        \bignb & \bignb & \color{red}{\bigbd} & \bigbd & \bignb\\
        \bignb & \bignb & \color{red}{\bignb} & \bignb & \bignb\\
        \bigvd & \bigvd & \color{red}{\bigvd} & \bigvd & \bigvd

        \end{matrix}$
        \end{tabular}}\\ and multicharge equals to $(14,13,16)$.
\end{itemize}

Notice that since we are considering another multicharge for $\bm\la^{+\bm k}$, namely the multicharge given by $(n_1+c_1, \ldots, n_r+c_r)$ and we are labelling the runners in the usual way for an abacus display with $e+1$ runners, the label of the each runner corresponds to the residue of the nodes corresponding to the beads in each runner.
\end{exe}

{
\begin{lem}
Let $e\geq2$, and let $\bm n=(n_1, \ldots, n_r)$ and $\tilde{\bm n}=(\tilde{n}_1, \ldots, \tilde{n}_r)$ be two choices of number of beads for the $e$-abacus display of $\bm\la$. Denote by $\bm\la_{\bm n}^{+\bm k}$ (respectively, $\bm\la_{\tilde{\bm n}}^{+\bm k}$) the multipartition obtained applying $^{+\bm k}$ to $\mathrm{Ab}_e(\bm\la)_{\bm n}$ (respectively, $\mathrm{Ab}_e(\bm\la)_{\tilde{\bm n}}$). If $(\tilde{n}_1, \ldots, \tilde{n}_r)=(n_1 + h_1, \ldots, n_r+h_r)$ with $h_j \in \mathbb{Z}$ and $h_j\equiv h_1$ mod $e$ for all $j \in \{1, \ldots, r\}$, then
\begin{enumerate}
\item $\bm\la^{+\bm k}_{\bm n}=\bm\la^{+\bm k}_{\tilde{\bm n}}$;
\item the underlying multicharges are the same up to a shift, that is, if $\bm a = (a_1, \ldots, a_r)$ is the multicharge corresponding to $\bm\la^{+\bm k}_{\bm n}$ and $\tilde{\bm a}= (\tilde{a}_1, \ldots, \tilde{a}_r)$ is the multicharge corresponding to $\bm\la_{\tilde{\bm n}}^{+\bm k}$, then there exists $s\in \{0, \ldots, e-1\}$ such that $$a_j \equiv \tilde{a}_j+s \mod (e+1).$$ 
Therefore, there is an isomorphism between the Ariki-Koike algebra with multicharge $\bm a$ and the Ariki-Koike algebra with multicharge $\tilde{\bm a}$.
\end{enumerate}
\end{lem}

\begin{proof}
Consider an $r$-multipartition $\bm\la$ and construct two abacus displays for $\bm\la$, one with $\bm n=(n_1, \ldots, n_r)$ beads and the other with $\tilde{\bm n}=(\tilde{n}_1, \ldots, \tilde{n}_r)$ beads where $\tilde{n}_j=n_j + h_j$ and $h_j \in \mathbb{Z}$ for each $j \in \{1, \ldots, r\}$.

Suppose that $h_1\equiv \ldots \equiv h_r$ mod $e$. Thus, write $h_j=c_{h_j}e + d_h$ with $0 \leq d_h \leq e-1$ for $j \in \{1, \ldots, r\}$. We want to show that $\bm\la_{\bm n}^{+\bm k}$ and $\bm\la_{\tilde{\bm n}}^{+\bm k}$ represent the same multipartition.
From the construction of $\bm\la_{\bm n}^{+\bm k}$ and $\bm\la_{\tilde{\bm n}}^{+\bm k}$, we have for $j\in \{1, \ldots, r\}$
\begin{align*}
n_j+k_j &=c_je + d;\\
\tilde{n}_j+k_j &= \tilde{c}_je + \tilde{d}.
\end{align*}
Hence, using $\tilde{n}_j=n_j + h_j$, we get
\begin{equation}\label{numbeads}
n_j + h_j + k_j = c_je + d + h_j = \tilde{c}_je + \tilde{d}
\end{equation}
for each $j \in \{1, \ldots, r\}$. 
Now, we distinguish two cases:
\begin{itemize}
\item if $h_j\equiv0$ mod $e$, i.e. $h_j=c_{h_j} e$ for some $c_{h_j} \in \mathbb{Z}$, the abacus display $\mathrm{Ab}_e(\bm \la)_{\tilde{\bm n}}$ is obtained from the abacus display $\mathrm{Ab}_e(\bm\la)_{\bm n}$ by adding or removing $c_{h_j}$ entire rows of beads at the top of the abacus. Thus, by \eqref{numbeads} we get
$$\tilde{n}_j +k_j = c_je + d + c_{h_j} e = (c_j + c_{h_j})e + d = \tilde{c}_je + \tilde{d}.$$ This implies $\tilde{c}_j=c_j+c_{h_j}$ and $\tilde{d}=d$, so we add a runner to the abacus display $\mathrm{Ab}_e(\bm \la)_{\tilde{\bm n}}$ immediately to the left of runner $d$ and put $c_j+c_{h_j}$ beads on this new runner, in the top $c_j+c_{h_j}$ positions. Hence, $\bm\la_{\bm n}^{+\bm k}$ and $\bm\la_{\tilde{\bm n}}^{+\bm k}$ represent the same multipartition because we add the runner in the same position of each abacus and the level of the last bead of the inserted runner is the same in both abacuses.
\item if $h_j=c_{h_j}e + d_h$ for some $c_{h_j} \in \mathbb{Z}$ and $d_h \in \{1, \ldots, e-1\}$, the abacus display $\mathrm{Ab}_e(\bm \la)_{\tilde{\bm n}}$ is obtained from the abacus $\mathrm{Ab}_e(\bm\la)_{\bm n}$ by adding or removing $c_{h_j}$ entire rows of beads and a row with $d_h$ beads at the top of the abacus. Thus, by \eqref{numbeads} we get
$$\tilde{n}_j +k_j = c_je + d + c_{h_j}e + d_h = (c_j + c_{h_j})e + (d+d_h) = \tilde{c}_je + \tilde{d}.$$ This implies $\tilde{c}_j=c_j+c_{h_j}$ and $\tilde{d}=d + d_h$, so we add a runner to the abacus display $\mathrm{Ab}_e(\bm \la)_{\tilde{\bm n}}$ immediately to the left of runner $d + d_h$ and put $c_j+c_{h_j}$ beads on this new runner, in the top $c_j+c_{h_j}$ positions. Hence, $\bm\la_{\bm n}^{+\bm k}$ and $\bm\la_{\tilde{\bm n}}^{+\bm k}$ represent the same multipartition because, in order to get $\bm\la_{\tilde{\bm n}}^{+\bm k}$, we add the runner in a position translated of $d_h$ compared to the one in $\bm\la_{\bm n}^{+\bm k}$ and the level of the last bead of the inserted runner is the same in both abacuses.

\end{itemize}
To conclude, the multicharge corresponding to $\bm\la_{\bm n}^{+\bm k}$ is $\bm a = (n_1+c_1, \ldots, n_r +c_r)$ and the multicharge corresponding to $\bm\la_{\tilde{\bm n}}^{+\bm k}$ is $\tilde{\bm a} = (n_1+h_1+c_1+c_{h_1}, \ldots, n_r+h_r+c_r+c_{h_r}).$ 
For all $j \in \{1, \ldots, r\}$, we can write $h_j=c_{h_j}e + s$ for some $c_{h_j} \in \mathbb{Z}$ and $s \in \{0, \ldots, e-1\}$ since $h_j\equiv h_1\mod e$. Hence, for all $j \in \{1, \ldots, r\}$
\begin{align*}
(n_j+h_j+c_j+c_{h_j})-(n_j+c_j) &= h_j + c_{h_j} \\
&= c_{h_j}e+s+c_{h_j} \\
&= c_{h_j}(e+1)+s\\
&\equiv s \mod (e+1).
\end{align*}
Thus, the isomorphism between the Ariki-Koike algebra with multicharge $\bm a$ and the one with multicharge $\tilde{\bm a}$ is given by
\begin{align*}
T_0 &\mapsto q^{-s}T_0\\
T_i &\mapsto T_i \text{ for all }i=1, \ldots, n-1.
\end{align*}
\end{proof}}

\begin{exe}
Let $\bm\la = ((4,3,2),(2^2),(3))$ and $e=4$, as in Example \ref{exrunnadd}. For $\bm k=(3,9,6)$,
\begin{itemize}
\item if we choose $\bm n=(n_1,n_2,n_3)=(11,9,12)$, we get the abacus display of $\bm\la^{+(3,9,6)}$ as in Example \ref{exrunnadd} that represents the multipartition $((5,3,2^2), (5,2^3), (3^2))$.
\item if we choose $\tilde{\bm n}=(\tilde{n}_1,\tilde{n}_2,\tilde{n}_3)=(10,8,11)$, we have the following truncated abacus for $\bm\la$:
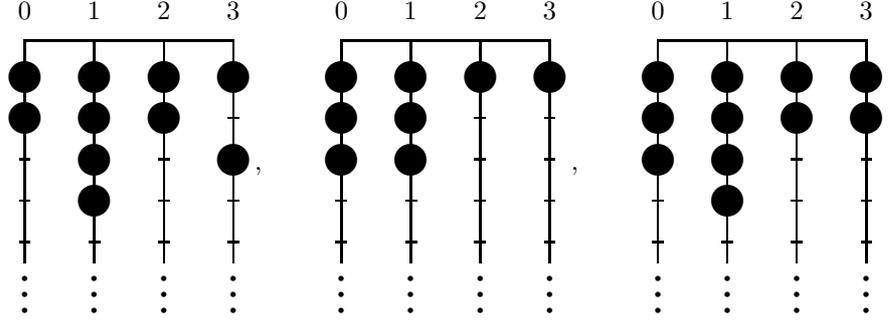
$$
        \begin{matrix}

        0 & 1 & 2 & 3\\

        \bigtl & \bigtm & \bigtm & \bigtr\\
        \bigbd & \bigbd & \bigbd & \bigbd\\
        \bigbd & \bigbd & \bigbd & \bignb\\
        \bignb & \bigbd & \bignb & \bigbd\\
        \bignb & \bigbd & \bignb & \bignb\\
        \bignb & \bignb & \bignb & \bignb\\
        \bigvd & \bigvd & \bigvd & \bigvd

        \end{matrix}, \qquad
        \begin{matrix}
        
        0 & 1 & 2 & 3\\

        \bigtl & \bigtm & \bigtm & \bigtr\\
        \bigbd & \bigbd & \bigbd & \bigbd\\
        \bigbd & \bigbd & \bignb & \bignb\\
        \bigbd & \bigbd & \bignb & \bignb\\
        \bignb & \bignb & \bignb & \bignb\\
        \bignb & \bignb & \bignb & \bignb\\
        \bigvd & \bigvd & \bigvd & \bigvd

        \end{matrix}, \qquad
        \begin{matrix}

        0 & 1 & 2 & 3\\

        \bigtl & \bigtm & \bigtm & \bigtr\\
        \bigbd & \bigbd & \bigbd & \bigbd\\
        \bigbd & \bigbd & \bigbd & \bigbd\\
        \bigbd & \bigbd & \bignb & \bignb\\
        \bignb & \bigbd & \bignb & \bignb\\
        \bignb & \bignb & \bignb & \bignb\\
        \bigvd & \bigvd & \bigvd & \bigvd

        \end{matrix}
.$$
Hence, $\bm\la^{+(3,9,6)}$ has the following abacus configuration

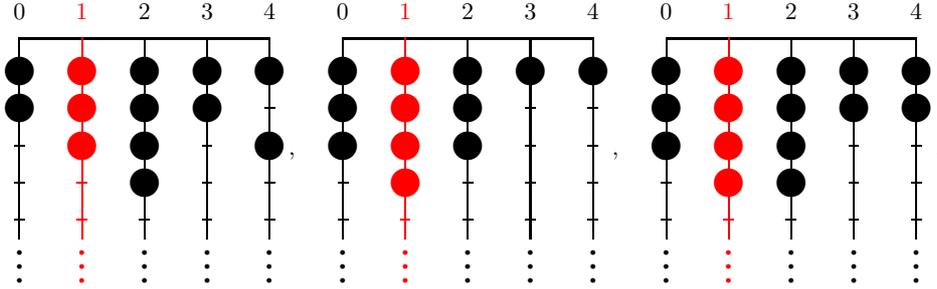
\scalebox{0.9}{
\begin{tabular}{c}
$\begin{matrix}

        0 & \color{red}{1} & 2 & 3 & 4\\

        \bigtl & \bigtm & \bigtm & \bigtm & \bigtr\\
        \bigbd & \color{red}{\bigbd} & \bigbd & \bigbd & \bigbd\\
        \bigbd & \color{red}{\bigbd} & \bigbd & \bigbd & \bignb\\
        \bignb & \color{red}{\bigbd} & \bigbd & \bignb & \bigbd\\
        \bignb & \color{red}{\bignb} & \bigbd & \bignb & \bignb\\
        \bignb & \color{red}{\bignb} & \bignb & \bignb & \bignb\\
        \bigvd & \color{red}{\bigvd} & \bigvd & \bigvd & \bigvd

        \end{matrix}, \quad
        \begin{matrix}
        
        0 & \color{red}{1} & 2 & 3 & 4\\

        \bigtl & \bigtm & \bigtm & \bigtm & \bigtr\\
        \bigbd & \color{red}{\bigbd} & \bigbd & \bigbd & \bigbd\\
        \bigbd & \color{red}{\bigbd} & \bigbd & \bignb & \bignb\\
        \bigbd & \color{red}{\bigbd} & \bigbd & \bignb & \bignb\\
        \bignb & \color{red}{\bigbd} & \bignb & \bignb & \bignb\\
        \bignb & \color{red}{\bignb} & \bignb & \bignb & \bignb\\
        \bigvd & \color{red}{\bigvd} & \bigvd & \bigvd & \bigvd

        \end{matrix}, \quad
        \begin{matrix}

        0 & \color{red}{1} & 2 & 3 & 4\\

        \bigtl & \bigtm & \bigtm & \bigtm & \bigtr\\
        \bigbd & \color{red}{\bigbd} & \bigbd & \bigbd & \bigbd\\
        \bigbd & \color{red}{\bigbd} & \bigbd & \bigbd & \bigbd\\
        \bigbd & \color{red}{\bigbd} & \bigbd & \bignb & \bignb\\
        \bignb & \color{red}{\bigbd} & \bigbd & \bignb & \bignb\\
        \bignb & \color{red}{\bignb} & \bignb & \bignb & \bignb\\
        \bigvd & \color{red}{\bigvd} & \bigvd & \bigvd & \bigvd

        \end{matrix},
$
\end{tabular}}\\
that represents again the multipartition $((5,3,2^2), (5,2^3), (3^2))$.
\end{itemize}
\end{exe}

We extend the operator $^{+\bm k}$ linearly to the whole Fock space. 
\begin{lem}\label{multireg_pres}
Let $\bm k$ be an $r$-tuple of non-negative integers. If $\bm\la$ is an $e$-multiregular $r$-multipartition, then $\bm\la^{+\bm k}$ is an $(e+1)$-multiregular $r$-multipartition.
\end{lem}

\begin{proof}
Suppose that $\bm\la$ is an $e$-multiregular $r$-multipartition. Then $\la^{(j)}$ is $e$-regular for each $j \in \{1, \ldots, r\}$, that is $\la^{(j)}$ has at most $e-1$ equal parts, equivalently it has at most $e-1$ consecutive beads in its abacus display. Hence, by construction each ${\la^{(j)}}^{+k^{(j)}}$ has at most $e$ consecutive beads in its abacus display and so ${\la^{(j)}}^{+k^{(j)}}$ is $(e+1)$-regular for each $j \in \{1, \ldots r\}$.
\end{proof}

\subsection{Induction operators and addition of a full runner}
In this section, we give some results that will be really helpful in the proof of our main theorem (Theorem \ref{canbasis_runrem}) in which we prove a runner removal-type theorem for Ariki-Koike algebras.

We will work in the following setting. Let $\bm\la$ be an $r$-multipartition. Consider the truncated abacus configuration of $\bm\la$ with $\bm n = (n_1, \ldots, n_r)$ beads. Set $\bm k=(k^{(1)}, \ldots, k^{(r)})$ to be an $r$-tuple of non-negative integers such that for all $j\in \{1, \ldots, r\}$ 
\begin{equation}\label{positiond}
n_j+k^{(j)} = c_je + d   
\end{equation} with $0 \leq d \leq e-1$. Then $\bm\la^{+\bm k}$ is the $r$-multipartition obtained from $\bm\la$ by adding a new runner with $c_j$ beads (as described in Section \ref{def+k}) to the left of runner $d$ in each component of $\bm\la$. This determines that the new inserted runner is labelled by $d$ mod $(e+1)$.
For our first lemma, following the proof of Lemma 3.5 in \cite{Frunrem}, we need to define the following function:
$$\begin{array}{cccc}
g \colon& (\zez)\setminus \{d\} &\rightarrow& (\mathbb{Z}/(e + 1)\mathbb{Z}) \setminus \{d,d + 1\}\\
 &d-i \mod e &\mapsto & d-i\mod (e+1)
\end{array}$$ 
for any $d \in I$ and $i = 1, \ldots, e - 1$.
\begin{lem}\label{newrunnernotdd+1}
Let $\bm\la$ and $\bm\xi$ be $r$-multipartitions. Let $\bm k=(k^{(1)}, \ldots, k^{(r)})$ be an $r$-tuple of non-negative integers and $d \in \{0, \ldots, e-1\}$ be the label of the new inserted runner of $\bm\la^{+\bm k}$ such that \eqref{positiond} holds. Suppose $i \in I\setminus \{d\}$. Then $\bm\la \xrightarrow{m:i} \bm\xi$ if and only if $\bm\la^{+\bm k} \xrightarrow{m:g(i)} \bm\xi^{+\bm k}$, and if this happens then we have $N_i(\bm\la, \bm\xi) = N_{g(i)}(\bm\la^{+\bm k}, \bm\xi^{+\bm k})$.
\end{lem}
\begin{proof}
We have $\bm\la \xrightarrow{m:i} \bm\xi$ if and only if the abacus display for $\bm\xi$ may be obtained by moving $m$ beads from runner $i-1$ to runner $i$, and in this case $N_i(\bm\la, \bm\xi)$ is determined by the configurations of these two runners in the two abacus displays. The fact that $i\neq d$ means that in constructing the abacus displays for $\bm\la^{+\bm k}$ and $\bm\xi^{+\bm k}$ the new runner is not added in between these two runners, and so the condition $\bm\la^{+\bm k} \xrightarrow{m:g(i)} \bm\xi^{+\bm k}$ and the coefficient $N_i(\bm\la, \bm\xi) = N_{g(i)}(\bm\la^{+\bm k}, \bm\xi^{+\bm k})$ are determined from these two runners in exactly the same way.
\end{proof}

\begin{cor}\label{corineqd}
Suppose $i \in I \setminus \{d\}$, $m \geq 1$ and $\bm\la$ is any multipartition. Then $\left(f_{i}^{(m)}(\bm\la)\right)^{+\bm k} =\mathfrak{F}_{g(i)}^{(m)}(\bm\la^{+\bm k})$.
\end{cor}
\begin{proof}
This is immediate from Lemma \ref{newrunnernotdd+1} and the description of the action of $f_{i}^{(m)}$ in Section \ref{algsec}.
\end{proof}

For the next lemma, we introduce the following notation. If $\bm\la$ and $\bm\xi$ are $r$-multipartitions, then we write $\bm\la \xrightarrow[m:i]{m:i+1} \bm\xi$ to indicate that $\bm\xi$ is obtained from $\bm\la$ by adding first $m$ addable $(i+1)$-nodes and then $m$ addable $i$-nodes. This notation is just a shorter version of the following one:
$$\bm\la \xrightarrow{m:i+1}\bm\nu\xrightarrow{m:i} \bm\xi$$
where $\bm\nu$ is an $r$-multipartition.
 
\begin{lem}\label{newrunnerdd+1}
Let $\bm\la$ and $\bm\xi$ be $r$-multipartitions. Let $\bm k=(k^{(1)}, \ldots, k^{(r)})$ be an $r$-tuple of non-negative integers and $d \in \{0, \ldots, e-1\}$ be the label of the new inserted runner of $\bm\la^{+\bm k}$ such that \eqref{positiond} holds. If the last bead in runner $d-1$ of each component $\la^{(j)}$ is at most in position $(c_j-1)e+d-1$ for all $j=1, \ldots, r$, 
then $\bm\la \xrightarrow{m:d} \bm\xi$ if and only if $\bm\la^{+\bm k} \xrightarrow[m:d]{m:d+1} \bm\xi^{+\bm k}$.
\end{lem}

Before proving this lemma, it is worth noticing the following.

\begin{rmk}\label{moveimplies}
In the assumptions of Lemma \ref{newrunnerdd+1}, every move of a bead $\mathfrak{b}$ at level $\ell$ from runner $d$ to runner $d+1$ in a component $J$ of $\bm\la^{+\bm k}$ determines uniquely a move of the bead $\bar{\mathfrak{b}}$ at level $\ell$ from runner $d-1$ to runner $d$ in the component $J$ of $\bm\nu$ where $\bm\nu$ is an $r$-multipartition such that $\bm\la^{+\bm k} \xrightarrow{m:d+1}\bm\nu$.
\end{rmk}

\begin{proof}[Proof of Lemma \ref{newrunnerdd+1}]
Suppose that $\bm\la \xrightarrow{m:d} \bm\xi$. This means that the abacus display for $\bm\xi$ may be obtained by moving $m$ beads from runner $d-1$ to runner $d$ of $\bm\la$. Notice that this moving of beads from runner $d-1$ to runner $d$ can occur simultaneously in different components of $\bm\la$, say $\la^{(j_1)}, \ldots, \la^{(j_s)}$ are the components of $\bm\la$ involved to get $\bm\xi$. Let $m^{(j_t)}$ be the number of beads moved in $\la^{(j_t)}$, for $t=1, \ldots, s$.  So, in each of these components there are at least $m^{(j_t)}$ levels $\ell_1, \ldots, \ell_{m^{(j_t)}}$ of the abacus that present a configuration of the type $ 
    \begin{matrix}
        \bd & \nb \\
    \end{matrix}$
in runners $d-1$ and $d$.

When we apply the operator $^{+\bm k}$ to $\bm\la$, by assumption we are going to add a new runner in between runner $d-1$ and runner $d$ in each component of $\bm\la$. This new runner by hypothesis has the last bead that is at least at the same level of the last bead in runner $d-1$ of every component of $\bm\la$. This implies that in each component $(\la^{(j_1)})^{+k^{(j_1)}}, \ldots, (\la^{(j_s)})^{+k^{(j_t)}}$, at each level $\ell_1, \ldots, \ell_{m^{(j_t)}}$ for $t=1, \ldots, s$ we have an abacus configuration of the type $$ \begin{matrix}
        d-1 & d & d+1 \\
        \bd & \bd & \nb \\
        \end{matrix}.$$ 
Now, consider the $r$-multipartition $\bm\mu$ such that $\bm\la^{+\bm k} \xrightarrow[m:d]{m:d+1}\bm\mu$ where the $a$ beads that we move from runner $d$ to runner $d+1$ are in the components $j_1, \ldots,j_s$ and at levels $\ell_1, \ldots, \ell_{m^{(j_t)}}$ for $t=1, \ldots, s$. Then $\bm\mu$ is exactly the multipartition $\bm\xi^{+\bm k}$. Indeed, given our assumption on the position of the last bead in the new runner, after moving $a$ beads from runner $d$ to runner $d+1$ as above, the only beads that we can move from runner $d-1$ to runner $d$ are the beads in the components $j_1, \ldots,j_s$ and at levels $\ell_1, \ldots, \ell_{m^{(j_t)}}$ for $t=1, \ldots, s$. This last move fills the empty spaces that we create with the first move in some runners $d$, giving then an abacus display with all runners $d$ full of beads and with the last bead in position $(c_j-1)(e+1)+d$ for $j=1, \ldots, r$.\\

Conversely, suppose that $\bm\la^{+\bm k} \xrightarrow[m:d]{m:d+1} \bm\xi^{+\bm k}$. We want to show that $\bm\la \xrightarrow{m:d} \bm\xi$. This, again, follows by our assumption on the position of the last bead in each runner $d-1$ of $\bm\la$. Indeed, this hypothesis implies that if we look at the abacus configuration of an addable $(d+1)$-node in any component of $\bm\la^{+\bm k}$, it will be one of the following:
\begin{enumerate}
\item $$\begin{matrix}
        d-1 & d & d+1 \\
        \bd & \bd & \nb \\
        \end{matrix},$$
\item $$\begin{matrix}
        d-1 & d & d+1 \\
        \nb & \bd & \nb \\
        \end{matrix}.$$
\end{enumerate}

We do not need to consider the second type of abacus configuration, because moving a bead from runner $d$ to runner $d+1$ cannot be followed by moving a bead from runner $d-1$ to runner $d$. Thus, there is no chance of getting $\bm\xi^{+\bm k}$ from this abacus configuration since every move of a bead from runner $d$ to runner $d+1$ determines uniquely a move of a bead from runner $d-1$ to runner $d$. 
Hence, suppose that the abacus configuration of all addable $(d+1)$-nodes of $\bm\la^{+\bm k}$ we move to get $\bm\xi^{+\bm k}$ is of type 1. In this case, if we move a bead $\mathfrak{b}$ from runner $d$ to runner $d+1$ in $\bm\la^{+\bm k}$, then we can only move the bead from runner $d-1$ to runner $d$ in the same component and at the same level of the bead $\mathfrak{b}$. Hence, if $j_1, \ldots, j_s$ are the components of $\bm\la^{+\bm k}$ involved to get $\bm\xi^{+\bm k}$ and $m^{(j_t)}$ is the number of beads moved in the component $j_t$, for $t=1, \ldots, s$, then $\bm\xi$ can be obtained from $\bm\la$ moving $m^{(j_t)}$ beads from runner $d-1$ to runner $d$ in the components $j_t$, for $t=1, \ldots, s$. 
\end{proof}

\begin{lem}\label{Ndd+1}
With the same assumption of Lemma \ref{newrunnerdd+1}. We have $$N_d(\bm\la, \bm\xi) = N_{d+1}(\bm\la^{+\bm k}, \bm\nu) + N_{d}(\bm\nu,\bm\xi^{+\bm k})$$ where $\bm\nu$ is the unique $r$-multipartition such that $\bm\la^{+\bm k} \xrightarrow{m:d+1}\bm\nu\xrightarrow{m:d} \bm\xi^{+\bm k}$.
\end{lem}

\begin{proof}
For $j=1, \ldots, r$, set $g^{(j)}e + d-1$ to be the position of the last bead of $\bm\la$ in runner $d-1$ of the component $\la^{(j)}$. 
Notice that, instead of considering all the $r$-multipartitions $\bm\nu$ such that $\bm\la^{+\bm k} \xrightarrow{m:d+1}\bm\nu$, we can restrict to consider only those $r$-multipartitions $\bm\nu$ where no one of the $m$ addable $(d+1)$-nodes of $\bm\la^{+\bm k}$, added to obtain $\bm\nu$, corresponds to a bead in position $x(e+1)+d$ with $g^{(j)}< x \leq c_j-1$ for each $j=1, \ldots, r$.
Indeed, if we consider $\bm\nu$ such that $\bm\la^{+\bm k} \xrightarrow{m:d+1}\bm\nu$ where at least one of the $m$ addable $(d+1)$-nodes added to $\bm\la^{+\bm k}$ is in position $x(e+1)+d$ with $g^{(J)}< x \leq c_J-1$ for a component $J$, then there is no multipartition $\bm \pi$ for which $\bm\nu\xrightarrow{m:d} \bm\pi$ because by Remark \ref{moveimplies} we should move the bead in position $x(e+1)+d-1$ in the component $J$ of $\bm\nu$, but there is no bead in that position. This can be seen in terms of abacus display in the following way: the abacus configuration at the level $x$ of the component $\nu^{(J)}$ in runners $d-1$, $d$, $d+1$ is
$$ \begin{matrix}
        d-1 & d & d+1 \\
        \nb & \nb & \bd \\
        \end{matrix},$$
from which it is clear that we have no chance of moving an addable $d$-node at level $x$ from runner $d-1$ to runner $d$. Thus, when we restrict to such $r$-multipartition $\bm\nu$, then we have
\begin{enumerate}
\item $\#\{\mathfrak{n} \in \bm\nu \setminus \bm\la^{+\bm k}\} = \#\{\mathfrak{n}\in \bm\xi\setminus \bm\la\}$;
\item $\#\{\mathfrak{n} \in \bm\xi^{+\bm k} \setminus \bm\nu\} = \#\{\mathfrak{n}\in \bm\xi\setminus \bm\la\}$ that is equivalent to 
$$\#\{\text{addable $d$-nodes of $\bm\nu$}\}=\#\{\text{addable $d$-nodes of $\bm\la$}\}.$$
\end{enumerate}
Indeed, provided that we are considering the $r$-multipartitions $\bm\nu$ satisfying the condition above we can conclude that these equalities holds for the following reasons.
\begin{enumerate}
\item Since the new inserted runner of $\bm\la^{+\bm k}$ is full of beads and with the last bead in a higher position than any last bead in every runner $d-1$ of $\bm\la$, the addable $(d+1)$-nodes of $\bm\la^{+\bm k}$ consists of the addable $d$-nodes of $\bm\la$ and the addable $(d+1)$-nodes of $\bm\la^{+\bm k}$ at level $x$ with $g^{(j)}< x \leq c_j-1$ for each $j\in\{1, \ldots r\}$, that is
\begin{align*}
\#&\{\text{addable $(d+1)$-nodes of $\bm\la^{+\bm k}$}\}\\
&=\#\{\text{addable $(d+1)$-nodes of $\bm\la^{+\bm k}$ at level $x\leq g^{(j)}$}\} \\
&\quad + \#\{\text{addable $(d+1)$-nodes of $\bm\la^{+\bm k}$ at level $g^{(j)}< x \leq c_j-1$}\}\\
&=\#\{\text{addable $d$-nodes of $\bm\la$}\} \\
&\quad+ \#\{\text{addable $(d+1)$-nodes of $\bm\la^{+\bm k}$ at level $g^{(j)}< x \leq c_j-1$}\}.
\end{align*}
Anyway, restricting to the $r$-multipartitions $\bm\nu$ as above means that we are excluding the $r$-multipartitions obtained from $\bm\la^{+\bm k}$ by adding addable $(d+1)$-nodes at level $x$ with $g^{(j)}< x \leq c_j-1$. So, in this case
$$\#\{\text{addable $(d+1)$-nodes of $\bm\la^{+\bm k}$}\}= \#\{\text{addable $d$-nodes of $\bm\la$}\}.$$
Hence, $\#\{\mathfrak{n} \in \bm\nu \setminus \bm\la^{+\bm k}\} = \#\{\mathfrak{n}\in \bm\xi\setminus \bm\la\}$.
\item By Remark \ref{moveimplies}, there is a correspondence between the addable $(d+1)$-nodes of $\bm\la^{+\bm k}$ and the addable $d$-nodes of $\bm\nu$ and so we can state that $\#\{\mathfrak{n} \in \bm\xi^{+\bm k} \setminus \bm\nu\}=\#\{\mathfrak{n} \in \bm\nu \setminus \bm\la^{+\bm k}\} = \# \{\mathfrak{n}\in \bm\xi\setminus \bm\la\}$.
\end{enumerate}
Recall that by definition
\begin{align*}
N_d(\bm\la, \bm\xi) = \sum_{\mathfrak{n} \in \bm\xi \setminus \bm\la} (&\#\{\text{addable $d$-nodes of $\bm\xi$ above }\mathfrak{n}\}\\
  -&\#\{\text{removable $d$-nodes of $\bm\la$ above }\mathfrak{n}\}),
\end{align*}
\begin{align*}
N_{d+1}(\bm\la^{+\bm k}, \bm\nu) = \sum_{\mathfrak{n} \in \bm\nu \setminus \bm\la^{+\bm k}} (&\#\{\text{addable $(d+1)$-nodes of $\bm\nu$ above }\mathfrak{n}\}\\
  -&\#\{\text{removable $(d+1)$-nodes of $\bm\la^{+\bm k}$ above }\mathfrak{n}\}),
\end{align*}
\begin{align*}
N_d(\bm\nu, \bm\xi^{+\bm k}) = \sum_{\mathfrak{n} \in \bm\xi^{+\bm k} \setminus \bm\nu} (&\#\{\text{addable $d$-nodes of $\bm\xi^{+\bm k}$ above }\mathfrak{n}\}\\
  -&\#\{\text{removable $d$-nodes of $\bm\nu$ above }\mathfrak{n}\}).
\end{align*}
Now, consider $\bm\la$ and $\mathfrak{n}\in \bm\xi\setminus \bm\la$. Let $J$ be the component of $\bm\xi$ of the node $\mathfrak{n}$. Set 
\begin{itemize}
\item $r_1$ to be the number of rows of the abacus of $\la^{(J)}$ above $\mathfrak{n}$ with a configuration of the following type in runners $d-1$ and $d$:
$$ \begin{matrix}
        d-1 & d  \\
        \bd & \bd \\
        \end{matrix},$$
\item $r_2$ to be the number of rows of the abacus of $\la^{(J)}$ above $\mathfrak{n}$ with a configuration of the following type in runners $d-1$ and $d$:
$$ \begin{matrix}
        d-1 & d  \\
        \bd & \nb \\
        \end{matrix},$$
\item $r_3$ to be the number of rows of the abacus of $\la^{(J)}$ above $\mathfrak{n}$ with a configuration of the following type in runners $d-1$ and $d$:
$$ \begin{matrix}
        d-1 & d  \\
        \nb & \bd \\
        \end{matrix},$$
\item $b$ to be the number of removable $d$-nodes in $\bm\xi\setminus\bm\la$ above $\mathfrak{n}$ in the component $J$.
\end{itemize}
In order to make easier to visualise the abacus display, we can assume that the rows of the abacus of $\la^{(J)}$ above $\mathfrak{n}$ with the same configuration between the runner $d-1$ and $d$ occurs as shown below:
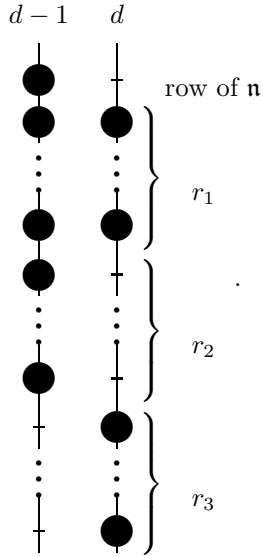
\NiceMatrixOptions{nullify-dots}
$$ \begin{NiceMatrix}[last-col=3]
    d-1 & d \\
    &&\\
    \bigbd & \bignb & \text{row of }\mathfrak{n}\\
    \bigbd & \bigbd & \\
    \bigvd & \bigvd & \quad r_1\\
    \bigbd & \bigbd & \\[1mm]
    \bigbd & \bignb & \\
    \bigvd & \bigvd & \quad r_2\\
    \bigbd & \bignb & \\[1mm]
    \bignb & \bigbd & \\
    \bigvd & \bigvd & \quad r_3\\
    \bignb & \bigbd & \\
    \CodeAfter
    \SubMatrix.{4-1}{6-2}\}
    \SubMatrix.{7-1}{9-2}\}
    \SubMatrix.{10-1}{12-2}\}
   \end{NiceMatrix}.$$
Then we can assume that the abacus display of component $J$ of $\bm\xi$ has the following configuration in the row of $\mathfrak{n}$ (with $\mathfrak{n}$ in \textcolor{red}{red}) and those above $\mathfrak{n}$:
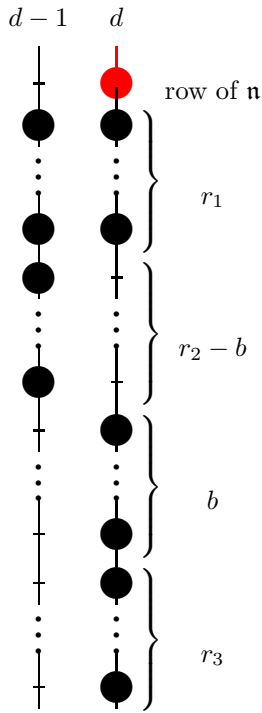
$$ \begin{NiceMatrix}
    d-1 & d \\
    &&\\
    \bignb & \color{red}{\bigbd} & \text{row of }\mathfrak{n}\\
    \bigbd & \bigbd & \\
    \bigvd & \bigvd & r_1\\
    \bigbd & \bigbd & \\[1mm]
    \bigbd & \bignb & \\
    \bigvd & \bigvd & r_2-b\\
    \bigbd & \bignb & \\[1mm]
    \bignb & \bigbd &\\
    \bigvd & \bigvd & b\\
    \bignb & \bigbd & \\[1mm]
    \bignb & \bigbd &\\
    \bigvd & \bigvd & r_3\\
    \bignb & \bigbd & \\
    \CodeAfter
    \SubMatrix.{4-1}{6-2}\}
    \SubMatrix.{7-1}{9-2}\}
    \SubMatrix.{10-1}{12-2}\}
    \SubMatrix.{13-1}{15-2}\}
     \end{NiceMatrix}.$$
Notice that the assumption on the order of the rows is not necessary for the scope of the proof, but it makes easier to represent the abacus display.
Thus, for $\mathfrak{n} \in \bm\xi\setminus\bm\la$ we have that $$N_d(\la^{(J)},\xi^{(J)})=(r_2-b)-r_3.$$
Thus, when we consider the component $J$ of the corresponding $r$-multipartitions $\bm\la^{+\bm k} \xrightarrow{m:d+1}\bm\nu\xrightarrow{m:d} \bm\xi^{+\bm k}$ in terms of abacus display above the nodes $\mathfrak{n}'$ and $\mathfrak{n}''$ (in \textcolor{red}{red}) that are the nodes corresponding to the node $\mathfrak{n}$, we have the following configurations between runners $d-1$, $d$, $d+1$:\\
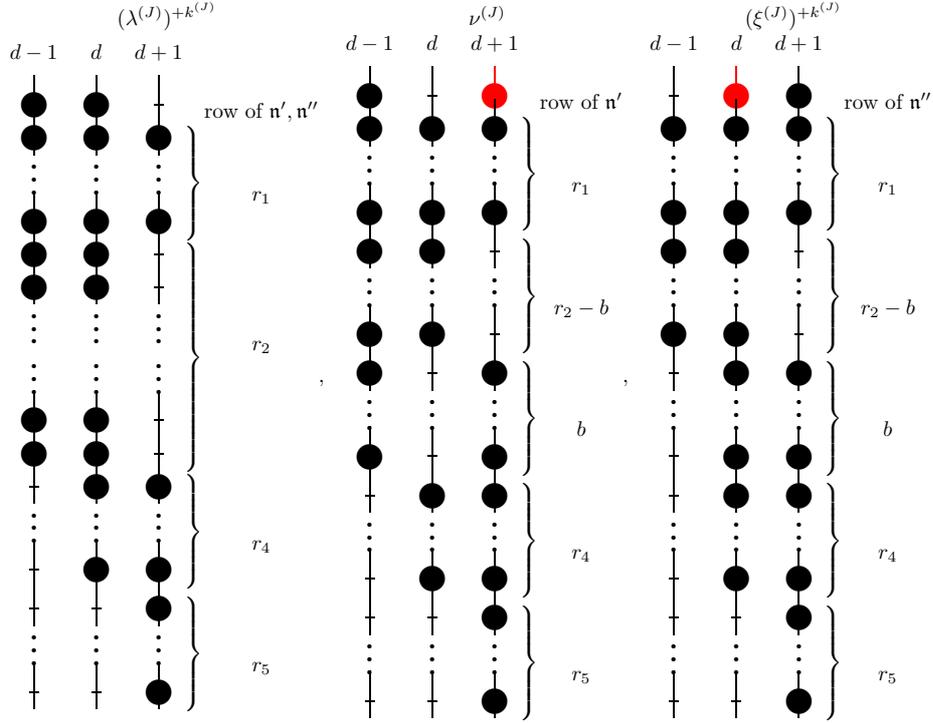
\begin{center}
\scalebox{0.80}{
\begin{tabular}{ccc}
$(\la^{(J)})^{+k^{(J)}}$ & $\nu^{(J)}$ & $(\xi^{(J)})^{+k^{(J)}}$\\
$ \begin{NiceMatrix}
    d-1 & d & d+1\\
    &&&\\
    \bigbd & \bigbd & \bignb & \text{row of } \mathfrak{n}', \mathfrak{n}''\\
    \bigbd & \bigbd & \bigbd & \\
    \bigvd & \bigvd & \bigvd & r_1\\
    \bigbd & \bigbd & \bigbd & \\
    \bigbd & \bigbd & \bignb & \\
    \bigbd & \bigbd & \bignb & \\
    \bigvd & \bigvd & \bigvd & r_2\\
    \bigvd & \bigvd & \bigvd & \\
    \bigbd & \bigbd & \bignb & \\
    \bigbd & \bigbd & \bignb & \\
    \bignb & \bigbd & \bigbd &\\
    \bigvd & \bigvd & \bigvd & r_4\\
    \bignb & \bigbd & \bigbd & \\[1mm]
    \bignb & \bignb & \bigbd & \\
    \bigvd & \bigvd & \bigvd & r_5\\
    \bignb & \bignb & \bigbd \\
    \CodeAfter
    \SubMatrix.{4-1}{6-3}\}
    \SubMatrix.{7-1}{12-3}\}
    \SubMatrix.{13-1}{15-3}\}
    \SubMatrix.{16-1}{18-3}\}
   \end{NiceMatrix},$
   &
$ \begin{NiceMatrix}

    d-1 & d & d+1\\
    &&&\\
    \bigbd & \bignb & \color{red}{\bigbd} & \text{row of }\mathfrak{n}'\\
    \bigbd & \bigbd & \bigbd & \\
    \bigvd & \bigvd & \bigvd & r_1\\
    \bigbd & \bigbd & \bigbd & \\[1mm]
    \bigbd & \bigbd & \bignb & \\
    \bigvd & \bigvd & \bigvd & r_2-b\\
    \bigbd & \bigbd & \bignb & \\[1mm]
    \bigbd & \bignb & \bigbd & \\
    \bigvd & \bigvd & \bigvd & b\\
    \bigbd & \bignb & \bigbd & \\[1mm]
    \bignb & \bigbd & \bigbd &\\
    \bigvd & \bigvd & \bigvd & r_4\\
    \bignb & \bigbd & \bigbd & \\[1mm]
    \bignb & \bignb & \bigbd & \\
    \bigvd & \bigvd & \bigvd & r_5\\
    \bignb & \bignb & \bigbd \\
    \CodeAfter
    \SubMatrix.{4-1}{6-3}\}
    \SubMatrix.{7-1}{9-3}\}
    \SubMatrix.{10-1}{12-3}\}
    \SubMatrix.{13-1}{15-3}\}
    \SubMatrix.{16-1}{18-3}\}
   \end{NiceMatrix},$
   &
   $ \begin{NiceMatrix}

    d-1 & d & d+1\\
    &&&\\
    \bignb & \color{red}{\bigbd} & \bigbd & \text{row of }\mathfrak{n}''\\
    \bigbd & \bigbd & \bigbd & \\
    \bigvd & \bigvd & \bigvd & r_1\\
    \bigbd & \bigbd & \bigbd & \\[1mm]
    \bigbd & \bigbd & \bignb & \\
    \bigvd & \bigvd & \bigvd & r_2-b\\
    \bigbd & \bigbd & \bignb & \\[1mm]
    \bignb & \bigbd & \bigbd & \\
    \bigvd & \bigvd & \bigvd & b\\
    \bignb & \bigbd & \bigbd & \\[1mm]
    \bignb & \bigbd & \bigbd &\\
    \bigvd & \bigvd & \bigvd & r_4\\
    \bignb & \bigbd & \bigbd & \\[1mm]
    \bignb & \bignb & \bigbd & \\
    \bigvd & \bigvd & \bigvd & r_5\\
    \bignb & \bignb & \bigbd \\
    \CodeAfter
    \SubMatrix.{4-1}{6-3}\}
    \SubMatrix.{7-1}{9-3}\}
    \SubMatrix.{10-1}{12-3}\}
    \SubMatrix.{13-1}{15-3}\}
    \SubMatrix.{16-1}{18-3}\}
   \end{NiceMatrix},$   
\end{tabular}}
\end{center}
where $r_4+r_5=r_3$.

Thus, for $\mathfrak{n}' \in \bm\nu\setminus\bm\la^{+\bm k}$ and $\mathfrak{n}'' \in \bm\xi^{+\bm k}\setminus\bm\nu$ we have that
\begin{gather*}
N_{d+1}((\la^{(J)})^{+k^{(J)}},\nu^{(J)})=(r_2-b)-r_5,\\
N_{d}(\nu^{(J)},(\xi^{(J)})^{+k^{(J)}})=0-r_4.
\end{gather*}
Hence, for the component $J$ we have $$N_d(\la^{(J)},\xi^{(J)})= N_{d+1}((\la^{(J)})^{+k^{(J)}},\nu^{(J)}) + N_{d}(\nu^{(J)},(\xi^{(J)})^{+k^{(J)}}).$$
By Proposition \ref{Ncomps}, to conclude we just need to show that for each $\mathfrak{n}\in \bm\xi\setminus\bm\la$ and for each component $j<J$
\begin{equation}\label{cpmj}
    N_d(\la^{(j)},\xi^{(j)})= N_{d+1}((\la^{(j)})^{+k^{(j)}},\nu^{(j)}) + N_{d}(\nu^{(j)},(\xi^{(j)})^{+k^{(j)}}).
\end{equation}
We can extend easily the previous argument to every component $j<J$ of $\bm\xi$. Indeed, if $j<J$ is a component of $\bm\xi$ then, by the total order of all addable and removable nodes of a multipartition, all the nodes in such a component $j$ of $\bm\xi$ are above $\mathfrak{n}$ and so in order to prove \eqref{cpmj} we can use exactly the same argument explained for component $J$, considering all the nodes in the component $j$, instead of just the ones above $\mathfrak{n}$.
Thus, we can conclude.
\end{proof}

\begin{cor}\label{corieqd}
Let $\bm\la$ be an $r$-multipartition. Let $\bm k=(k^{(1)}, \ldots, k^{(r)})$ be an $r$-tuple of non-negative integers and $d \in \{0, \ldots, e-1\}$ be the label of the new inserted runner of $\bm\la^{+\bm k}$ such that \eqref{positiond} holds. If the last bead in runner $d-1$ of each component $\la^{(j)}$ is at most in position $(c_j-1)e+d-1$ for all $j=1, \ldots, r$, then $\left(f_{d}^{(m)}(\bm\la)\right)^{+\bm k} =\mathfrak{F}_{d}^{(m)}\mathfrak{F}_{d+1}^{(m)}(\bm\la^{+\bm k})$.
\end{cor}
\begin{proof}
This is immediate from Lemmas \ref{newrunnerdd+1} and \ref{Ndd+1} and the description of the action of $f_{i}^{(m)}$ in Section \ref{algsec}.
\end{proof}

\subsection{Case $r\geq2$}
In this section we prove the so called `full runner removal' theorem for Ariki-Koike algebras for $r\geq2$. In particular, we will compute canonical basis vectors $G^{\bm s}(\bm\mu)$ and, respectively, $G^{\bm s^+}(\bm\mu^{+\bm k})$ in the Fock spaces $\mathcal{F}^{\bm s}$ and, respectively  $\mathcal{F}^{\bm s^+}$, rather than working with the Ariki-Koike algebras directly. So, we will apply induction operators $\mathfrak{f}$ (respectively, $\mathfrak{F}$)  on multipartitions with $r\geq 2$ components and thus we need to take in to account the role that each of these components is playing.

 Let $r\geq2$. We start imposing some conditions on the non-negative integers $k^{(1)}, \ldots, k^{(r)}$. Let $k^{(1)}, \ldots, k^{(r)}$ be non-negative integers such that
\begin{align}
    & k^{(r)} \geq \mu_1^{(r)}, \text{ and;}\nonumber \\
    & k^{(j)}-k^{(h)} \geq \mu_1^{(j)}+e-1+\sum\limits_{t=h+1}^{j-1}|\mu^{(t)}| \text{ for all }1\leq h<j \leq r. \label{conditiononk^j}
\end{align}

\begin{prop}\label{emptyto+k_r}
Let $\bm s =(s_1, \dots, s_r)\in I^r$ be a multicharge and $(\varnothing, \bm\mu)=(\varnothing, \mu^{(2)}, \ldots, \mu^{(r)})$ be a $r$-multipartition. Let $(k^{(1)}, \bm k) = (k^{(1)}, k^{(2)} \ldots, k^{(r)})$ be an $r$-tuple of non-negative integers such that conditions \eqref{conditiononk^j} hold.
Write $k^{(1)}=k_1^{(1)}e + k_2^{(1)}$ with $k_1^{(1)}\geq 0$ and $0 \leq k_2^{(1)}\leq e-1$. Denote by $\bm a^{+(k^{(1)}, \bm k)}=(a_1, \ldots, a_r)$ the multicharge associated to the multipartition $(\varnothing^{+k^{(1)}}, {\mu^{(2)}}^{+k^{(2)}}, \ldots, {\mu^{(r)}}^{+k^{(r)}})$. Then, setting $\alpha = k^{(1)}+a_1$ and $\beta= k^{(1)}-k_2^{(1)}+a_1$,
\begin{equation*}
\mathfrak{F}^{(k_1^{(1)})}_{\overline{\alpha}} \mathfrak{F}^{(k_1^{(1)})}_{\overline{\alpha -1}}\ldots\mathfrak{F}^{(k_1^{(1)})}_{\overline{\beta+1}}\mathcal{G}^{(k_1^{(1)}-1)}_{\overline{\beta}}\mathcal{G}^{(k_1^{(1)}-2)}_{\overline{\beta+1}} \ldots \mathcal{G}^{(1)}_{\overline{\beta +k_1^{(1)}-2}}((\varnothing, \bm\mu^{+\bm k}))=(\varnothing^{+k^{(1)}}, \bm\mu^{+\bm k}).    
\end{equation*}
where $\mathfrak{F}^{(k_1^{(1)})}_{\overline{\alpha}} \mathfrak{F}^{(k_1^{(1)})}_{\overline{\alpha-1}}\ldots\mathfrak{F}^{(k_1^{(1)})}_{\overline{\beta+1}}$ occur if $\alpha \neq \beta$.
\end{prop}

\begin{proof}
{Construct the truncated $e$-abacus configuration of the $r$-multipartition $(\varnothing, \bm\mu)$ consisting of $\bm n=(n_1,\ldots, n_r)$ beads with $n_j\equiv s_j \mod e$ for all $j$.} For $j=1, \ldots, r$, write
\begin{equation*}
   n_j+k^{(j)}=c_je+d 
\end{equation*}
with $0 \leq d\leq e-1$. 
Construct now the truncated $(e+1)$-abacus display of the $r$-multipartition $(\varnothing, \bm\mu^{+\bm k})$ with $(n_1+c_1, \ldots, n_r+c_r)$ beads. Notice that the multicharge associated to this multipartition is $\bm a^{+(k^{(1)},\bm k)}=(n_1+c_1, \ldots, n_r+c_r)$. Thus $a_j = n_j + c_j$ for $1 \leq j \leq r$.
We proceed now by induction on $k^{(1)}$ {with the following induction hypothesis. If $k = k_1 e + k_2$ with $k_1\geq 0$ and $0 \leq k_2\leq e-1$ is a non-negative integer such that $k<k^{(1)}$, then it holds that
\begin{align*}
\mathfrak{F}^{(k_1)}_{\overline{a}} \ldots\mathfrak{F}^{(k_1)}_{\overline{b+1}}\mathcal{G}^{(k_1-1)}_{\overline{b}} \ldots \mathcal{G}^{(1)}_{\overline{b +k_1-3}}((\varnothing, \bm\mu^{+\bm k}))=(\varnothing^{+k}, \bm\mu^{+\bm k}),
\end{align*}
where $a = k+a_1$, $b= k-k_2+a_1$, and $\mathfrak{F}^{(k_1)}_{\overline{a}} \ldots\mathfrak{F}^{(k_1)}_{\overline{b+1}}$ only occur if $a\neq b$.}

If $k^{(1)}<e$, then $\varnothing^{+k^{(1)}}= \varnothing$ by Proposition \ref{empty+k}. So, $$(\varnothing^{+k^{(1)}}, \bm\mu^{+\bm k}) = (\varnothing, \bm\mu^{+\bm k}).$$

Suppose $k^{(1)} \geq e$ and write $k^{(1)}=k_1^{(1)}e + k_2^{(1)}$ with $k_1^{(1)}\geq 1$ and $0 \leq k_2^{(1)}\leq e-1$. Then by Proposition \ref{indseq} we have $$\mathfrak{F}^{(k^{(1)}_1)}_{\overline{k^{(1)}}} \mathfrak{F}^{(k_1^{(1)})}_{\overline{k^{(1)}-1}}\ldots\mathfrak{F}^{(k^{(1)}_1)}_{\overline{k^{(1)}-k_2^{(1)}+1}}\mathcal{G}^{(k_1^{(1)}-1)}_{\overline{k^{(1)}-k_2^{(1)}}}\mathcal{G}^{(k_1^{(1)}-2)}_{\overline{k^{(1)}-k_2^{(1)}+1}} \ldots \mathcal{G}^{(1)}_{\overline{k^{(1)}-k_2^{(1)}+k_1^{(1)}-2}}(\varnothing)=\varnothing^{+k^{(1)}},$$
where $\mathfrak{F}^{(k^{(1)}_1)}_{\overline{k^{(1)}}} \mathfrak{F}^{(k_1^{(1)})}_{\overline{k^{(1)}-1}}\ldots\mathfrak{F}^{(k^{(1)}_1)}_{\overline{k^{(1)}-k_2^{(1)}+1}}$ only occur if $k_2^{(1)}\neq 0$.
Hence, we want to apply this induction sequence to $(\varnothing, \bm\mu^{+\bm k})$ and show that the only resulting multipartition is $(\varnothing^{+k^{(1)}}, \bm\mu^{+\bm k})$. The first fact that we need to consider when we deal with multipartitions is the multicharge. Therefore, the residues involved in the induction sequence need to be translated by the multicharge.
So, the induction sequence that we want to apply to $(\varnothing, \bm\mu^{+\bm k})$ is the following:
$$
\mathfrak{F}^{(k_1^{(1)})}_{\overline{\alpha}} \mathfrak{F}^{(k_1^{(1)})}_{\overline{\alpha -1}}\ldots\mathfrak{F}^{(k_1^{(1)})}_{\overline{\beta+1}}\mathcal{G}^{(k_1^{(1)}-1)}_{\overline{\beta}}\mathcal{G}^{(k_1^{(1)}-2)}_{\overline{\beta+1}} \ldots \mathcal{G}^{(1)}_{\overline{\beta +k_1^{(1)}-2}}$$
with $\alpha = k^{(1)}+a_1$ and $\beta= k^{(1)}-k_2^{(1)}+a_1$.

If $k^{(1)}=e$, then by Proposition \ref{empty+k} $\varnothing^{+e}= \varnothing$. So, $$(\varnothing^{+e}, {\mu^{(2)}}^{+k^{(2)}}, \ldots, {\mu^{(r)}}^{+k^{(r)}}) = (\varnothing, {\mu^{(2)}}^{+k^{(2)}}, \ldots, {\mu^{(r)}}^{+k^{(r)}}).$$

{If $k^{(1)} > e$, then for $\alpha\neq \beta$ consider 
\begin{equation}\label{indk-e}
  \mathfrak{F}^{(k_1^{(1)})}_{\overline{\alpha}} \mathfrak{F}^{(k_1^{(1)})}_{\overline{\alpha -1}}\ldots\mathfrak{F}^{(k_1^{(1)})}_{\overline{\beta+1}}\mathfrak{F}^{(k_1^{(1)}-1)}_{\overline{\beta}} \mathfrak{F}^{(k_1^{(1)}-1)}_{\overline{\beta -1}}\ldots\mathfrak{F}^{(k_1^{(1)}-1)}_{\overline{\alpha+1}}(\varnothing^{+(k^{(1)}-e-1)}, \bm\mu^{+\bm k}) 
\end{equation}
while for $\alpha=\beta$ consider
\begin{equation}\label{indk-ea=b}
  \mathfrak{F}^{(k_1^{(1)}-1)}_{\overline{\alpha}} \mathfrak{F}^{(k_1^{(1)}-1)}_{\overline{\alpha -1}}\ldots\mathfrak{F}^{(k_1^{(1)}-1)}_{\overline{\alpha-e-1}}\mathfrak{F}^{(k_1^{(1)}-2)}_{\overline{\alpha+1}}(\varnothing^{+(k^{(1)}-e-1)}, \bm\mu^{+\bm k}).
\end{equation}
Notice that in \eqref{indk-e} and in \eqref{indk-ea=b} the number of induction operators $\mathfrak{F}_i^{(a)}$ is exactly $e+1$.} Moreover, we have 
$$\overline{\alpha}=\overline{k^{(1)}+a_1}= \overline{k^{(1)}+n_1+c_1} = \overline{c_1e+d + c_1} = \overline{c_1(e+1)+d}.$$
Thus $\alpha \equiv d$ mod $(e+1)$. 
{We want to show now that 
$$\eqref{indk-e}=(\varnothing^{+k^{(1)}}, \bm\mu^{+\bm k});$$ 
$$\eqref{indk-ea=b}=(\varnothing^{+k^{(1)}}, \bm\mu^{+\bm k}).$$
However, note that in \eqref{indk-e} and \eqref{indk-ea=b}
\begin{itemize}
\item the residues involved in the induction sequences are the same and in the same order;
\item the first operator is applied one time less than the last operator. 
\end{itemize}
Hence, the argument we give in the following works exactly in the same way in the two cases. So, we will show only that $\eqref{indk-e}=(\varnothing^{+k^{(1)}}, \bm\mu^{+\bm k})$ because a similar argument applies to prove that $\eqref{indk-ea=b}=(\varnothing^{+k^{(1)}}, \bm\mu^{+\bm k})$.

Showing that $\eqref{indk-e}=(\varnothing^{+k^{(1)}}, \bm\mu^{+\bm k})$ is equivalent to showing that this induction sequence is applied only to the first component of $(\varnothing^{+(k^{(1)}-e-1)}, \bm\mu^{+\bm k})$.} Indeed, if this sequence acts only on the first component we have $\eqref{indk-e} =(\varnothing^{+k^{(1)}}, \bm\mu^{+\bm k})$ by Proposition \ref{indseq}. 
Recall that the action of an induction operator $\mathfrak{F}_i^{(a)}$ for $a\geq1$ on the abacus display of a multipartition consists of moving $a$ beads from runner $i-1$ to an empty position in runner $i$ in some components. Suppose that this induction sequence acts not only on the first component of the multipartition $(\varnothing^{+(k^{(1)}-e-1)}, \bm\mu^{+\bm k})$ and it does not give 0. In particular, suppose that 
\begin{itemize}
    \item $\mathfrak{F}_{\overline{\alpha}}^{(p_{e+1})} \ldots \mathfrak{F}_{\overline{\alpha+1}}^{(p_{1})}$ is the induction subsequence acting on $\bm\mu^{+\bm k}$,
    \item $\mathfrak{F}_{\overline{\alpha}}^{(q_{e+1})} \ldots \mathfrak{F}_{\overline{\alpha+1}}^{(q_{1})}$ is the induction subsequence acting on $\varnothing^{+(k^{(1)}-e-1)}$,
\end{itemize}
with $q_{e+1}=k_1^{(1)}-p_{e+1}$, $q_1=k_1^{(1)}-1-p_1$ and at least one $p_i\neq0$. Since we apply $\mathfrak{F}_{\overline{\alpha+1}}^{(p_1)}$ to $\bm\mu^{+\bm k}$, this means that we are moving $p_1$ beads from runner $d$ to runner $d+1$ in some components of $\bm\mu^{+\bm k}$. Call such components $j_1, \ldots, j_b$. By construction of $\bm\mu^{+\bm k}$, this corresponds in terms of abacus display to having $p_1$ empty spaces in the runners $d$ of components $j_1, \ldots, j_b$ of $\bm\mu^{+\bm k}$. The action of the terms $\mathfrak{F}_{\overline{\alpha-1}}^{(p_e)} \ldots \mathfrak{F}_{\overline{\alpha+2}}^{(p_2)}$ involves {moving beads from runner $d+1$ to runner $d+2$, and then from runner $d+2$ to runner $d+3$, and so forth until we move beads from runner $d-2$ to runner $d-1$.
For $j \in \{2,\dots,r\}$, let $l_j$ be the level of the last bead of the new inserted runner in the component $j$ of $\bm\mu^{+\bm k}$. There are now two cases to consider:
\begin{enumerate}
    \item the beads moved by this induction in runner $d-1$ are at most at level $l_j$ for all components $j$;
    \item one of the beads moved by this induction in runner $d-1$ is at level $l_j+1$ for some $j$.
\end{enumerate}}
\textbf{Case 1.} If the moved beads are at most at the same level of the last bead of each new inserted runner of $\bm\mu^{+\bm k}$, then the number of beads that the last induction operator $\mathfrak{F}_{\overline{\alpha}}^{(p_{e+1})}$ can move from runner $d-1$ to runner $d$ can be at most $p_1$. {If not (i.e., $p_{e+1}>p_1$), then the action of the induction subsequence of $\bm\mu^{+\bm k}$ is $0$ because there are no enough gaps in runner $d$ of $\mathfrak{F}_{\overline{\alpha-1}}^{(p_e)} \ldots \mathfrak{F}_{\overline{\alpha+2}}^{(p_2)}\mathfrak{F}_{\overline{\alpha+1}}^{(p_1)}(\bm\mu^{+\bm k})$ where to move beads from runner $d-1$ to runner $d$. 
In fact, the gaps in runner $d$ of $\mathfrak{F}_{\overline{\alpha-1}}^{(p_e)} \ldots \mathfrak{F}_{\overline{\alpha+2}}^{(p_2)}\mathfrak{F}_{\overline{\alpha+1}}^{(p_1)}(\bm\mu^{+\bm k})$ are $p_1$. Indeed, the first operator $\mathfrak{F}^{(p_1)}_{\overline{\alpha+1}}$ moves $p_1$ beads from runner $d$ to runner $d+1$ of $\bm\mu^{+\bm k}$. Also, the beads moved by $\mathfrak{F}_{\overline{\alpha-1}}^{(p_e)} \ldots \mathfrak{F}_{\overline{\alpha+2}}^{(p_2)}$ in runner $d-1$ are, for all $j$, at most at level $l_j$, that is the level of the last bead in the new inserted runner of component $j$ of $\bm\mu^{+\bm k}$.} Thus, we have $p_{e+1}\leq p_1$. Hence,
\begin{align*}
    q_1&=k_1^{(1)}-1-p_1 \leq k_1^{(1)}-1,\\
    q_{e+1}&=k_1^{(1)}-p_{e+1} \geq k_1^{(1)}-p_1=q_1+1.
\end{align*}
We claim that the action of this induction sequence on $\varnothing^{+(k^{(1)}-e-1)}$ is $0$. {If $q_{e+1}>q_1+1$, then we can conclude that $\mathfrak{F}_{\overline{\alpha}}^{(q_{e+1})} \ldots \mathfrak{F}_{\overline{\alpha+1}}^{(q_{1})}(\varnothing^{+(k^{(1)}-e-1)})=0$ by Lemma \ref{seqindnot0}. If $q_{e+1}=q_1+1$, by contradiction we assume that $\mathfrak{F}_{\overline{\alpha}}^{(q_{e+1})} \ldots \mathfrak{F}_{\overline{\alpha+1}}^{(q_{1})}(\varnothing^{+(k^{(1)}-e-1)})$ is non-zero. By Lemma \ref{seqindnot0} this means that $q_{e+1}=\ldots=q_{e+2-k_2}=k_1$ and $q_{e+1-k_2}=\ldots=q_1=k_1-1$. However, this is a contradiction because this would imply that $p_i=0$ for all $i$, while we are assuming that at least one of the $p_i$'s is non-zero.
Therefore, we proved that the induction sequence acts non-zero when all the induction operators are applied to the first component of $(\varnothing^{+(k^{(1)}-e-1)}, \bm\mu^{+\bm k})$. Hence we conclude in this case.} \\

\textbf{Case 2.} If one of the moved beads is at a higher level than the last bead of one of the new inserted runners of $\bm\mu^{+\bm k}$, then the number of beads that the last induction operator $\mathfrak{F}_{\overline{\alpha}}^{(p_{e+1})}$ can move from runner $d-1$ to runner $d$ can be at most $p_1+1$. {If not (i.e. $p_{e+1}>p_1+1$), then the action of the induction subsequence of $\bm\mu^{+\bm k}$ is $0$ because there are no enough gaps in runners $d$ of $\mathfrak{F}_{\overline{\alpha-1}}^{(p_e)} \ldots \mathfrak{F}_{\overline{\alpha+2}}^{(p_2)}\mathfrak{F}_{\overline{\alpha+1}}^{(p_1)}(\bm \mu^{+\bm k})$ where to move $p_{e+1}>p_1+1$ beads from runner $d-1$ to runner $d$. 
In fact, the total number of gaps in runner $d$ of $\mathfrak{F}_{\overline{\alpha-1}}^{(p_e)} \ldots \mathfrak{F}_{\overline{\alpha+2}}^{(p_2)}\mathfrak{F}_{\overline{\alpha+1}}^{(p_1)}(\bm \mu^{+\bm k})$ created by the first operator $\mathfrak{F}^{(p_1)}_{\overline{\alpha+1}}$ is $p_1$. Also, by the assumption of Case 2., there is a component $J$ of $\mathfrak{F}_{\overline{\alpha-1}}^{(p_e)} \ldots \mathfrak{F}_{\overline{\alpha+2}}^{(p_2)}\mathfrak{F}_{\overline{\alpha+1}}^{(p_1)}(\bm\mu^{+\bm k})$ with a bead in runner $d-1$ at level $l_J+1$ and with a gap on its right. This is because $l_J$ is the level of the last bead in the new inserted runner of component $J$ of $\bm\mu^{+\bm k}$. Thus, in total, there are at most $p_1+1$ possible positions where to move beads from runner $d-1$ to runner $d$ in $\mathfrak{F}_{\overline{\alpha-1}}^{(p_e)} \ldots \mathfrak{F}_{\overline{\alpha+2}}^{(p_2)}\mathfrak{F}_{\overline{\alpha+1}}^{(p_1)}(\bm\mu^{+\bm k})$.} Thus, we have $p_{e+1}\leq p_1+1$. Hence,
\begin{align*}
    q_1&=k_1^{(1)}-1-p_1 \leq k_1^{(1)}-1,\\
    q_{e+1}&=k_1^{(1)}-p_{e+1} \geq k_1^{(1)}-p_1-1=q_1.
\end{align*}
We claim that the action of this induction sequence on $\varnothing^{+(k^{(1)}-e-1)}$ is $0$. In order to prove this, we deal with the following three cases separately:
\begin{enumerate}
    \item[a.] $q_{e+1}=q_1$;
    \item[b.] $q_{e+1}=q_1+1$;
    \item[c.] $q_{e+1}>q_1+1$.
\end{enumerate}
For cases b. and c. we can conclude as in Case 1. with analogous arguments. As for case a., we need to be more cautious. In this case we have $q_{e+1}=q_1$ and so by Lemma \ref{seqindnot0} the induction subsequence acting on $\varnothing^{+k^{(1)}-e-1}$ is non-zero if and only if $$q_1=q_2= \ldots =q_x=q_{x+1}= \ldots =q_e=q_{e+1}.$$
We want to show that the induction subsequence acting on $\varnothing^{+k^{(1)}-e-1}$ acts as 0.
If $q_1=q_2= \ldots =q_{e+1-k_2^{(1)}}=q_{e+2-k_2^{(1)}}= \ldots =q_e=q_{e+1}$, then 
\begin{align*}
    &p_1= \ldots =p_{e+1-k_2^{(1)}}=k_1^{(1)}-q_1-1,\\
    &p_{e+2-k_2^{(1)}}= \ldots =p_{e+1}=k_1^{(1)}-q_{e+1}=k_1^{(1)}-q_1.
\end{align*}
In this case 
we have no problem with the induction subsequence acting on $\varnothing^{+k^{(1)}-e-1}$, but the induction subsequence on $\bm\mu^{+\bm k}$ is the one that acts as 0. We need again to distinguish two cases.
\begin{itemize}
    \item If the first $p_1$ beads moved from runner $d$ to runner $d+1$ involve at least one bead not in the last $p_1$ positions of the new inserted runner, then the action of the induction subsequence on $\bm\mu^{+\bm k}$ is 0. Indeed, the last induction operator requires to move $p_{e+1}$ beads from runner $d-1$ to runner $d$, but in runner $d$ we have at most $p_{e+1}-1=p_1$ empty positions.
    {\item If the first $p_1$ beads moved from runner $d$ to runner $d+1$ are exactly in the last $p_1$ positions of the new inserted runners of component $J$ of $\bm\mu^{+\bm k}$, then the action of the induction subsequence on $\bm\mu^{+\bm k}$ is 0. Notice that our assumption 
    $$k^{(j)}-k^{(h)} \geq \mu_1^{(j)}+\sum\limits_{t=h+1}^{j-1}|\mu^{(t)}|+e-1 \text{ for all }1\leq h<j \leq r$$
    implies that the last bead of $\mu^{(J)}$ is at most at level $c_J-k^{(1)}_1$ of the abacus of $\mu^{(J)}$. Indeed, the difference in height between the last bead of the new inserted runner and the last bead of $\mu^{(J)}$ is given by
    $c_J-1-x$ where $\mu^{(J)}_1+n_J-1=xe+i$ for $x\geq0$ and $0 \leq i \leq e-1$ and so we get that
    \begin{align*}
        c_J-1-x &= \dfrac{1}{e}(k^{(J)}+n_J-d)-1-\dfrac{1}{e}(\mu^{(J)}_1+n_J-1-i)\\
        & = \dfrac{1}{e}(k^{(J)}-\mu^{(J)}_1)+ \dfrac{1}{e}(-d-e+1+i)\\
        & \geq \dfrac{1}{e}(k^{(1)}+\sum\limits_{t=2}^{J-1}|\mu^{(t)}|+e-1) + \dfrac{1}{e}(-d-e+1+i) \\
        & \geq  \dfrac{1}{e}(k^{(1)}+e-1)+ \dfrac{1}{e}(-d-e+1+i)\\
        & =  \dfrac{1}{e}(k_1^{(1)}e+k_2^{(1)})+ \dfrac{1}{e}(-d-e+1+i+e-1)\\
        & =  k_1^{(1)}+ \dfrac{1}{e}(k_2^{(1)}-d+i)\\
        & >  k_1^{(1)}+ \dfrac{1}{e}(-e)\\
        & =  k_1^{(1)}-1.
    \end{align*}
    If the induction subsequence on $\bm\mu^{+\bm k}$ is 0 at some step before the last induction operator then we are done. Suppose that we can apply all the induction operators until the last one $\mathfrak{F}_{\bar{\alpha}}^{(p_{e+1})}$, and that we get something non-zero. Then the action of this last operator is 0 because we do not have enough beads in runner $d-1$ to move in runner $d$ in the abacus display of $\mathfrak{F}_{\overline{\alpha-1}}^{(p_e)} \ldots \mathfrak{F}_{\overline{\alpha+2}}^{(p_2)}\mathfrak{F}_{\overline{\alpha+1}}^{(p_1)}(\bm\mu^{+\bm k})$. In fact, the position $(c_J-1)(e+1)+(d-1)$ is empty since $x < c_J-k^{(1)}_1\leq c_J-1$ and so there are only $p_{e+1}-1$ addable nodes of residue $d$.}
    \end{itemize}
Therefore, we proved that the induction sequence acts non-zero when all the induction operators are applied to the first component of $(\varnothing^{+(k^{(1)}-e-1)}, \bm\mu^{+\bm k})$ also in this case.
So, we get that $\eqref{indk-e}= (\varnothing^{+k^{(1)}}, \bm\mu^{+\bm k})$.
By induction hypothesis we know also that
\begin{align*}
(\varnothing^{+(k^{(1)}-e-1)}, \bm\mu^{+\bm k})=
\mathfrak{F}^{(k_1^{(1)}-1)}_{\overline{\alpha}} \ldots\mathfrak{F}^{(k_1^{(1)}-1)}_{\overline{\beta+1}}\mathcal{G}^{(k_1^{(1)}-2)}_{\overline{\beta}} \ldots \mathcal{G}^{(1)}_{\overline{\beta +k_1^{(1)}-3}}((\varnothing, \bm\mu^{+\bm k})).
\end{align*}
Hence, we can conclude.    
\end{proof}

The above proposition gives the following result about the canonical basis coefficients of $(\varnothing,\bm\mu^{+\bm k})$ and $(\varnothing^{+k^{(1)}},\bm\mu^{+\bm k})$.

\begin{prop}\label{canbasis_empty}
Let $(\varnothing, \bm\mu)=(\varnothing, \mu^{(2)}, \ldots, \mu^{(r)})$ be a $r$-multipartition and $(k^{(1)}, \bm k) = (k^{(1)}, k^{(2)} \ldots, k^{(r)})$ be an $r$-tuple of non-negative integers such that conditions \eqref{conditiononk^j} hold.
Let $\bm s^+ \in (\Z/(e+1)\Z)^{r}$. Suppose that
$$G^{\bm s^+}_{e+1}((\varnothing,\bm\mu^{+\bm k}))=\sum_{\bm\mu\trianglerighteq\bm\la}d_{\bm\la\bm\mu}(v) (\varnothing,\bm\la^{+\bm k})$$
where $d_{\bm\la\bm\mu}(v) \in v\mathbb{N}[v]$ for $\bm\la \neq \bm\mu$. Then
$$G^{\bm s^+}_{e+1}((\varnothing^{+k^{(1)}},\bm\mu^{+\bm k}))=\sum_{\bm\mu\trianglerighteq\bm\la}d_{\bm\la\bm\mu}(v) (\varnothing^{+k^{(1)}},\bm\la^{+\bm k}).$$
\end{prop}

\begin{proof}
Let $\bm s^+=(s_1, \ldots, s_r)$ be a multicharge. Suppose that $G^{\bm s^+}_{e+1}((\varnothing,\bm\mu^{+\bm k}))=\sum_{\bm\mu\trianglerighteq\bm\la}d_{\bm\la\bm\mu}(v) (\varnothing,\bm\la^{+\bm k})$ with $d_{\bm\la\bm\mu}(v) \in v\mathbb{N}[v]$ for $\bm\la \neq \bm\mu$. We want to apply the induction sequence from $\varnothing$ to $\varnothing^{+k^{(1)}}$ given by Proposition \ref{indseq} to $G^{\bm s^+}_{e+1}((\varnothing,\bm\mu^{+\bm k}))$. Hence, we want to use Proposition \ref{emptyto+k_r}. So, writing $k^{(1)}=k_1^{(1)}e + k_2^{(1)}$ with $k_1^{(1)}\geq 0$ and $0 \leq k_2^{(1)}\leq e-1$, we want to apply the following induction sequence
\begin{equation}\label{actionG}
    \mathfrak{F}^{(k_1^{(1)})}_{\overline{\alpha}} \mathfrak{F}^{(k_1^{(1)})}_{\overline{\alpha -1}}\ldots\mathfrak{F}^{(k_1^{(1)})}_{\overline{\beta+1}}\mathcal{G}^{(k_1^{(1)}-1)}_{\overline{\beta}}\mathcal{G}^{(k_1^{(1)}-2)}_{\overline{\beta+1}} \ldots \mathcal{G}^{(1)}_{\overline{\beta +k_1^{(1)}-2}}(G^{\bm s^+}_{e+1}(\varnothing, \bm\mu^{+\bm k})),
\end{equation}
where
\begin{itemize}
    \item $\alpha = k^{(1)}+s_1$ and $\beta= k^{(1)}-k_2^{(1)}+s_1$;
    \item $\mathfrak{F}^{(k_1^{(1)})}_{\overline{\alpha}} \mathfrak{F}^{(k_1^{(1)})}_{\overline{\alpha-1}}\ldots\mathfrak{F}^{(k_1^{(1)})}_{\overline{\beta+1}}$ occur if $\alpha \neq \beta$.
\end{itemize}

Our assumptions are exactly the hypotheses of Proposition \ref{emptyto+k_r}. Thus, we have
\begin{align*}
    \eqref{actionG} &=\sum_{\bm\mu\trianglerighteq\bm\la}d_{\bm\la\bm\mu}(v) \mathfrak{F}^{(k_1^{(1)})}_{\overline{\alpha}} \mathfrak{F}^{(k_1^{(1)})}_{\overline{\alpha -1}}\ldots\mathfrak{F}^{(k_1^{(1)})}_{\overline{\beta+1}}\mathcal{G}^{(k_1^{(1)}-1)}_{\overline{\beta}}\mathcal{G}^{(k_1^{(1)}-2)}_{\overline{\beta+1}} \ldots \mathcal{G}^{(1)}_{\overline{\beta +k_1^{(1)}-2}}((\varnothing,\bm\la^{+\bm k}))\\
    &=\sum_{\bm\mu\trianglerighteq\bm\la}d_{\bm\la\bm\mu}(v) (\varnothing^{+k^{(1)}},\bm\la^{+\bm k}),
\end{align*}
which is of the form $(\varnothing^{+k^{(1)}},\bm\mu^{+\bm k})+\sum\limits_{\bm\mu\neq\bm\la}d_{\bm\la\bm\mu}(v)(\varnothing^{+k^{(1)}},\bm\la^{+\bm k})$.
Hence, this shows that $G^{\bm s^+}_{e+1}((\varnothing^{+k^{(1)}},\bm\mu^{+\bm k}))=\sum\limits_{\bm\mu\trianglerighteq\bm\la}d_{\bm\la\bm\mu}(v) (\varnothing^{+k^{(1)}},\bm\la^{+\bm k})$ by the uniqueness of the canonical basis.
\end{proof}
Write $\mathfrak{f}_{i_l}^{(h_l)} \cdots \mathfrak{f}_{i_1}^{(h_1)}$ with $h_1, \ldots, h_l$ non-negative integers and $i_1, \ldots, i_l \in I$. Fix $d\in \{0, \ldots, e-1\}$. Then define $\mathfrak{F}$ to be the induction sequence obtained in the following way: for all $j=1, \ldots, l$
\begin{itemize}
    \item if $i_j\neq d$, replace $\mathfrak{f}_{i_j}^{(h_j)}$ with $\mathfrak{F}^{(h_j)}_{g(i_j)}$;
    \item if $i_j=d$, replace $\mathfrak{f}_d^{(h_j)}$ with $\mathfrak{F}^{(h_j)}_{d}\mathfrak{F}^{(h_j)}_{d+1}$.
\end{itemize}

\begin{prop}\label{same_coeff}

Let $(\varnothing,\bm\mu)=(\varnothing,\mu^{(2)}, \ldots, \mu^{(r)})$ be an $e$-multiregular $r$-multipartition and $(k^{(1)},\bm k )= (k^{(1)}, k^{(2)}\ldots, k^{(r)})$ be an $r$-tuple of non-negative integers such that conditions \eqref{conditiononk^j} hold.
Let $d \in \{0, \ldots, e-1\}$ be the label of the new inserted runner of $(\varnothing,\bm\mu)^{+(k^{(1)},\bm k )}$ such that \eqref{positiond} holds. Suppose that 
$$\mathfrak{f} \cdot G_e^{\bm s}((\varnothing, \bm\mu))=\sum_{\bm\nu\in \mathcal{P}^{r}}g_{\bm\nu}\bm\nu,$$
where $g_{\bm\nu}\in \mathbb{Z}[q, q^{-1}]$ and $\mathfrak{f}=f_{i_l}^{(h_l)} \cdots f_{i_1}^{(h_1)}$ for some $h_1, \ldots, h_l \in \mathbb{Z}_{\geq0}$ and $i_1, \ldots, i_l \in I$ is such that $\mathfrak{f}\cdot\varnothing = \mu^{(1)} +\sum\limits_{\mu^{(1)} \triangleright  \tau} t_{\tau}\tau$ for $t_{\tau}\in \Z[q, q^{-1}]$. Then $$\mathfrak{F} \cdot G_{e+1}^{\bm s^{+}}((\varnothing, \bm\mu)^{+(k^{(1)},\bm k )})=\sum_{\bm\nu\in \mathcal{P}^{r}}g_{\bm\nu}\bm\nu^{+(k^{(1)},\bm k )},$$
where $\mathfrak{F}$ is defined as above.
\end{prop}

\begin{proof}
By definition of $^{+\bm k}$ we have that
$$\left(\mathfrak{f} \cdot G_e^{\bm s}((\varnothing, \bm\mu))\right)^{+(k^{(1)},\bm k )}=\sum_{\bm\nu\in \mathcal{P}^{r}}g_{\bm\nu}\bm\nu^{+(k^{(1)},\bm k )}.$$
Moreover, we get that
\begin{align*}
    \left(\mathfrak{f} \cdot G_e^{\bm s}((\varnothing, \bm\mu))\right)^{+(k^{(1)},\bm k )}&=\left(\mathfrak{f} \cdot \sum_{\mu\trianglerighteq\la}d_{\bm\la\bm\mu}(v) (\varnothing,\bm\la) \right)^{+(k^{(1)},\bm k )}\\
    &=\mathfrak{F} \cdot \sum_{\bm\mu\trianglerighteq\bm\la}d_{\bm\la\bm\mu}(v)(\varnothing,\bm\la)^{+(k^{(1)},\bm k )} && \text{by Cor. \ref{corineqd}, \ref{corieqd}}\\
    &=\mathfrak{F} \cdot G^{\bm s^+}_{e+1}((\varnothing,\bm\mu)^{+(k^{(1)},\bm k )}) && \text{by Prop. \ref{canbasis_empty}}
\end{align*}
Hence,
$$\left(\mathfrak{f} \cdot G_e^{\bm s}((\varnothing,\bm\mu))\right)^{+(k^{(1)},\bm k )} = \mathfrak{F} \cdot G_{e+1}^{\bm s^{+}}((\varnothing,\bm\mu)^{+(k^{(1)},\bm k )}).$$
Thus, we can conclude $$\mathfrak{F} \cdot G_{e+1}^{\bm s^{+}}((\varnothing,\bm\mu)^{+(k^{(1)},\bm k )})=\sum_{\bm\nu\in \mathcal{P}^{r}}g_{\bm\nu}\bm\nu^{+(k^{(1)},\bm k )}.$$
\end{proof}

\begin{thm}\label{canbasis_runrem}
Let $(\mu^{(1)},\bm\mu)=(\mu^{(1)}, \mu^{(2)}\ldots, \mu^{(r)})\in \mathcal{P}^r$ is an $e$-multiregular multipartition of $n$ and $\bm s = (s_1, \ldots, s_r)\in I^{r}$. Let $(k^{(1)},\bm k) =(k^{(1)},k^{(2)} \ldots, k^{(r)})$ with $k^{(j)}$ a non-negative integer for $j=1,\ldots,r$ such that conditions \eqref{conditiononk^j} hold.
Then $$G^{\bm s^{+}}_{e+1}((\mu^{(1)},\bm\mu)^{+(k^{(1)},\bm k)})=G^{\bm s}_{e}((\mu^{(1)},\bm\mu))^{+(k^{(1)},\bm k)}.$$
\end{thm}

\begin{proof}
Suppose that $(\mu^{(1)},\bm\mu)=(\mu^{(1)}, \ldots, \mu^{(r)})$ is an $e$-multiregular $r$-multipartition. We proceed by induction on the number $r$ of components.

If $r=1$, then $(\mu^{(1)},\bm\mu)=\mu$ and $(k^{(1)},\bm k)=k$, i.e., we are in the partition case. Thus, by Theorem \ref{Fullrunrem_level1} we have 
\begin{equation*}
G_{e+1}(\mu^{+k})=G_{e}(\mu)^{+k},
\end{equation*}
since $k\geq\mu_1$.

Suppose $r>1$. By induction on $r$, we know that for an $e$-multiregular $(r-1)$-multipartition $\bm\mu$
\begin{equation}\label{firstrel}
    G^{\bm s^{+}}_{e+1}(\bm\mu^{+\bm k})=G^{\bm s}_{e}(\bm\mu)^{+\bm k}.
\end{equation}
More explicitly, we have that $$G_{e}^{\bm s}(\bm \mu)=\bm\mu + \sum_{\bm\mu\triangleright\bm\la}d_{\bm\la\bm\mu}(v) \bm\la,$$
then \eqref{firstrel} implies that $$G_{e+1}^{\bm s^+}(\bm\mu^{+\bm k})=\bm\mu^{+\bm k} + \sum_{\bm\mu\triangleright\bm\la}d_{\bm\la\bm\mu}(v) \bm\la^{+\bm k}.$$
Now, we want to show that this is also true for the $r$-multipartition $(\mu^{(1)},\bm\mu)$.
By Corollary \ref{fcomp} it holds that
$$G^{\bm s^+}_{e+1}((\varnothing,\bm\mu^{+\bm k}))=(\varnothing,\bm \mu^{+\bm k}) + \sum_{\bm\mu\triangleright\bm\la}d_{\bm\la\bm\mu}(v) (\varnothing,\bm\la^{+\bm k}).$$
By Proposition \ref{canbasis_empty} we have that
$$G^{\bm s^+}_{e+1}((\varnothing^{+k^{(1)}},\bm\mu^{+\bm k}))=(\varnothing^{+k^{(1)}},\bm\mu^{+\bm k}) + \sum_{\bm\mu\triangleright\bm\la}d_{\bm\la\bm\mu}(v) (\varnothing^{+k^{(1)}},\bm\la^{+\bm k}).$$
Using the LLT algorithm on partitions, we can write $G^{(s_1)}(\mu^{(1)})$ as $\mathfrak{f}\cdot\varnothing$ in the Fock space $\mathcal{F}^{(s_1)}$, for some $\mathfrak{f}\in \mathcal{U}$.
Applying the induction sequence $\mathfrak{f}$ to $G_e^{\bm s}((\varnothing,\bm\mu))$ we can write
\begin{equation}\label{Ge(0,mu)}
   \mathfrak{f} \cdot G_e^{\bm s}((\varnothing, \bm\mu))=\sum_{\bm\nu\in \mathcal{P}^{r}}g_{\bm\nu}\bm\nu, 
\end{equation}
where $g_{\bm\nu}\in \mathbb{Z}[q, q^{-1}]$ because $\mathfrak{f} \cdot G_e^{\bm s}((\varnothing, \bm\mu))\in M^{\otimes \bm s}$ and $M^{\otimes \bm s}$ is a $\mathcal{U}$-submodule of $\mathcal{F}^{\bm s}$.
Performing step (c) of the LLT algorithm for multipartitions in \cite{Fay10} we get 
\begin{equation}\label{fGe(0,mu)}
  \mathfrak{f} \cdot G_e^{\bm s}((\varnothing, \bm\mu)) - \sum_{(\mu^{(1)},\bm\mu) \triangleright \bm\sigma} a_{\bm\sigma \bm\mu}(v) G_e^{\bm s}(\bm\sigma) = G_e^{\bm s}((\mu^{(1)}, \bm\mu)) 
\end{equation}
where $a_{\bm\sigma \bm\mu}(v) \in \mathbb{Z}[q+q^{-1}]$.

Now consider the induction sequence $\mathfrak{F}$ that is obtained translating the induction sequence $\mathfrak{f}$ from $e$ to $e+1$. This means that for each $i\in I$
\begin{itemize}
    \item if $i\neq d$, we replace $f_{i}^{(h)}$ with $\mathfrak{F}^{(h)}_{g(i)}$,
    \item if $i=d$, we replace $f_{d}^{(h)}$ with $\mathfrak{F}^{(h)}_{d}\mathfrak{F}^{(h)}_{d+1}$. 
\end{itemize}
 
Then apply $\mathfrak{F}$ to $G_{e+1}^{\bm s^{+}}((\varnothing^{+k^{(1)}},\mu^{+k^{(2)}}))$. By Proposition \ref{same_coeff} we have 
\begin{equation}\label{Ge+(0+,mu+)}
   \mathfrak{F} \cdot G_{e+1}^{\bm s^{+}}((\varnothing^{+k^{(1)}},\bm\mu^{+\bm k}))=\sum_{\bm\nu\in \mathcal{P}^{r}}g_{\bm\nu}\bm\nu^{+(k^{(1)}, \bm k)}. 
\end{equation}

We now consider $({\color{red}3.2.23})$. 
Since the coefficients occurring in the sum are exactly the same of \eqref{Ge(0,mu)}, we perform the following subtraction of terms

\begin{equation}\label{alice}
    \mathfrak{F} \cdot G_{e+1}^{\bm s^{+}}((\varnothing^{+k^{(1)}},\bm\mu^{+\bm k})) - \sum_{(\mu^{(1)},\bm\mu) \triangleright \bm\sigma} a_{\bm\sigma (\mu^{(1)},\bm\mu)}(v) G_{e+1}^{\bm s^{+}}(\bm\sigma^{+(k^{(1)},\bm k)}).
\end{equation}

We proceed by induction on the dominance order. Suppose that $G_{e+1}^{\bm s^{+}}(\bm\sigma^{+(k^{(1)},\bm k)}) = (G_e^{\bm s}(\bm\sigma))^{+(k^{(1)},\bm k)}$ for all $\bm\sigma \triangleleft (\mu^{(1)},\bm\mu)$. Then
\begin{align}\label{indFG-G+}
({\color{red}3.2.24})= \mathfrak{F} \cdot G_{e+1}^{\bm s^{+}}((\varnothing^{+k^{(1)}},\bm\mu^{+\bm k})) - \sum_{(\mu^{(1)},\bm\mu) \triangleright \bm\sigma} a_{\bm\sigma (\mu^{(1)},\bm\mu)}(v) (G_{e}^{\bm s}(\bm\sigma))^{+(k^{(1)},\bm k)}.
\end{align}
Since we are performing exactly the same operations of \eqref{fGe(0,mu)} and the starting coefficients are the same in \eqref{Ge(0,mu)} and $({\color{red}3.2.23})$, 
by definition of $^{+ \bm k}$ we have 
\begin{equation*}
({\color{red}3.2.25}) = G_e^{\bm s}((\mu^{(1)},\bm\mu))^{+(k^{(1)},\bm k)}.   
\end{equation*}
Moreover, by uniqueness of the canonical basis of $M^{\otimes\bm s^+}$ we can state
$$ 
({\color{red}3.2.25})= G_{e+1}^{\bm s^+}(({\mu^{(1)}}^{+k^{(1)}}, \bm\mu^{+\bm k})).$$
Hence,
$$G^{\bm s^{+}}_{e+1}((\mu^{(1)},\bm\mu)^{+(k^{(1)},\bm k)})=G^{\bm s}_{e}((\mu^{(1)},\bm\mu))^{+(k^{(1)},\bm k)}.$$
\end{proof}

\section{Example of empty runner removal for Ariki-Koike algebras}
This last section aims to motivate why we believe that an empty runner removal should hold for Ariki-Koike algebras.
An empty runner removal theorem is established for the Iwahori-Hecke algebras of the symmetric groups and the $q$-Schur algebras by James and Mathas in \cite{JM02}. Also, we have seen so far how much the representation theory of Ariki-Koike algebras resembles the one of the symmetric groups and how often results for the symmetric groups (or, the Iwahori-Hecke algebras of the symmetric groups) can be extended to Ariki-Koike algebras. Furthermore, all the examples we have examined with the help of GAP confirmed our belief. An instance of them is the following.

\begin{exe}
Let $r=2$, $e=3$, and write the set $I=\mathbb{Z}/3\mathbb{Z}$ as $\{0,1,2\}$.  Take $\bm s=(2,1)$. Consider $\bm\mu = ((2,1),(1))$ and its 3-abacus display

$$
        \begin{matrix}

        0 & 1 & 2\\

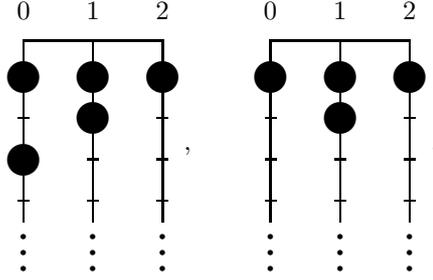
        \bigtl & \bigtm & \bigtr\\
        \bigbd & \bigbd & \bigbd\\
        \bignb & \bigbd & \bignb\\
        \bigbd & \bignb & \bignb\\
        \bignb & \bignb & \bignb\\
        \bigvd & \bigvd & \bigvd

        \end{matrix}, \qquad
        \begin{matrix}
        
        0 & 1 & 2\\

        \bigtl & \bigtm & \bigtr\\
        \bigbd & \bigbd & \bigbd\\
        \bignb & \bigbd & \bignb\\
        \bignb & \bignb & \bignb\\
        \bignb & \bignb & \bignb\\
        \bigvd & \bigvd & \bigvd

        \end{matrix}
.$$

The canonical basis element $G^{\bm s}_3(((2,1),(1)))$ is given by
\begin{align*}
G^{\bm s}_3(((2,1),(1)))&=((2,1),(1))+q((1^3),(1))+q^2((1^2),(1^2))\\
&= \bm\mu + q \bm\la_1 + q^2 \bm\la_2.
\end{align*}

Now, consider the multipartition with 4-abacus display
$$
        \begin{matrix}

        0 & 1 & 2 & 3\\

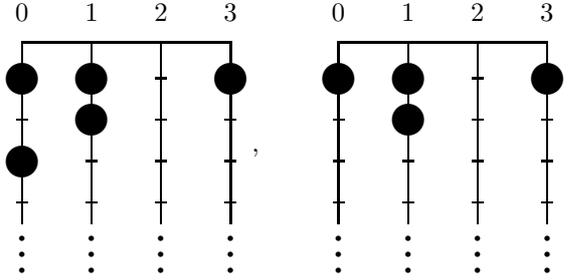
        \bigtl & \bigtm & \bigtm & \bigtr\\
        \bigbd & \bigbd & \bignb & \bigbd\\
        \bignb & \bigbd & \bignb & \bignb\\
        \bigbd & \bignb & \bignb & \bignb\\
        \bignb & \bignb & \bignb & \bignb\\
        \bigvd & \bigvd & \bigvd & \bigvd

        \end{matrix}, \qquad
        \begin{matrix}
        
        0 & 1 & 2 & 3\\

        \bigtl & \bigtm & \bigtm & \bigtr\\
        \bigbd & \bigbd & \bignb & \bigbd\\
        \bignb & \bigbd & \bignb & \bignb\\
        \bignb & \bignb & \bignb & \bignb\\
        \bignb & \bignb & \bignb & \bignb\\
        \bigvd & \bigvd & \bigvd & \bigvd

        \end{matrix}
,$$
that corresponds to the multipartition $((4,2,1),(2,1))$. Notice that the abacus display of $((4,2,1),(2,1))$ is obtained from the 3-abacus display of $\bm\mu$ by adding a runner with no beads to the left of each runner $2$ and relabelling the runners in the usual way. Write $\bm\mu^{+\varnothing}$ for $((4,2,1),(2,1))$ and $\bm s^{+\varnothing}$ for the corresponding multicharge, that is $(1,0)$.

Now, we compute the canonical basis element for $\bm\mu^{+\varnothing}$. We get
\begin{equation*}
G^{\bm s^{+\varnothing}}_4(((4,2,1),(2,1)))=((4,2,1),(2,1))+q((3,2^2),(2,1))+q^2((3,2,1),(2^2)).
\end{equation*}
Hence, we notice that
$$G^{\bm s^{+\varnothing}}_4(\bm\mu^{+\varnothing}) = \bm\mu^{+\varnothing} + q \bm\la_1^{+\varnothing} + q^2 \bm\la_2^{+\varnothing},$$
where $\bm\la_t^{+\varnothing}$ for $t\in\{1,2\}$ is obtained from the 3-abacus display of $\bm\la_t$ by adding a runner with no beads to the left of runner 2 in each component. Thus, in this case we have that the $v$-decomposition numbers match up, that is for $t\in\{1,2\}$
$$d_{\bm\la_t\bm\mu}^e(v)=d_{\bm\la_t^{+\varnothing}\bm\mu^{+\varnothing}}^{e+1}(v).$$
\end{exe}

Therefore, we conjecture that the following claim should hold.

\begin{conj}\label{conj1}
Let $\bm\lambda, \bm\mu \in \mathcal{P}^r$ be in a block $B$ of $\mathcal{H}_{r,n}$ with $\bm\mu$ $e$-multiregular. Suppose that $\bm\lambda^{+\varnothing}$ and $\bm\mu^{+\varnothing}$ are the multipartitions obtained from the $e$-abacus display of $\bm\lambda$ and $\bm\mu$ by adding an `empty' runner in each of their components. Then
$$d_{\bm\lambda\bm\mu}^e(v)=d_{\bm\lambda^{+\varnothing}\bm\mu^{+\varnothing}}^{e+1}(v).$$
\end{conj}

However, the argument we used to prove the `full' runner removal theorem is most likely not adaptable in a straightforward way to the proof of the above conjecture, because it heavily relies on the fact that the runner we add to each component is long enough. 

Moreover, for Iwahori-Hecke algebras (where one usually looks to take inspiration to prove results for Ariki-Koike algebras), the proof of the `empty' runner removal theorem provided by James and Mathas in \cite{JM02} makes use of $q$-Schur algebras. They first show that certain decomposition numbers of the $q$-Schur algebras $\mathcal{S}_{\C,q}(n)$ and $\mathcal{S}_{\C,q'}(m)$ are equal for specified $m>n$. Then, since the decomposition matrix for the Iwahori-Hecke algebras $H_{\C,q}(\mathfrak{S}_n)$ is a submatrix of the decomposition matrix of $\mathcal{S}_{\C,q}(n)$, they deduce the analogous result for the $v$-decomposition numbers for Iwahori-Hecke algebras. 

Although $q$-Schur algebras are defined also for Ariki-Koike algebras, a proper checking on how much of the proof for Iwahori-Hecke algebras can be extended to the Ariki-Koike algebras case needs to be done.  Some further work is therefore required to adapt this proof for Ariki-Koike algebras.

\section*{Acknowledgements}
This paper forms part of work towards a PhD degree under the supervision of Dr. Sin\'ead Lyle at University of East Anglia, and the author wishes to thank her for her direction and support throughout.


\addcontentsline{toc}{chapter}{Bibliography}
\bibliographystyle{alpha}
\bibliography{bibfile}

\begin{thebibliography}{WKM17}

\bibitem[AK94]{AK94}
S.~Ariki and K.~Koike.
\newblock A {Hecke} algebra of {$(\mathbb{Z}/r\mathbb{Z}) \wr \mathfrak{S}_n$} and construction of its irreducible representations.
\newblock {\em Adv. Math.}, 106(2):216--243, 1994.

\bibitem[Ari96]{A96}
S.~Ariki.
\newblock On the decomposition numbers of the {Hecke} algebra of {$G(r,1,n)$}.
\newblock {\em J. Math. Kyoto Univ.}, 36:789--808, 1996.

\bibitem[Ari01]{A01}
S.~Ariki.
\newblock On the classification of simple modules for cyclotomic {Hecke algebras of type $G(m, 1, n)$} and {Kleshchev} multipartitions.
\newblock {\em Osaka J. Math.}, 38:827--837, 2001.

\bibitem[BK09]{BK09}
J.\ Brundan and A.\ Kleshchev.
\newblock Graded decomposition numbers for cyclotomic {Hecke algebras}.
\newblock {\em Adv.\ Math.}, 222:1883--1942, 2009.

\bibitem[BKW11]{BKW11}
J.~Brundan, A.~Kleshchev, and W.~Wang.
\newblock Graded specht modules.
\newblock {\em Journal f\"ur die reine und angewandte {Mathematik}}, 2011(655):61 -- 87, 2011.

\bibitem[Del23]{mythesis}
Alice Dell'Arciprete.
\newblock {\em Decomposition numbers of Ariki-Koike algebras}.
\newblock Phd thesis, University of East Anglia, May 2023.

\bibitem[DJM98]{DJM98}
R.~Dipper, G.~D. James, and A.~Mathas.
\newblock Cyclotomic $q$-{Schur} algebras.
\newblock {\em Math. Z.}, 229(3):385--416, 1998.

\bibitem[DM02]{dm02}
R.~Dipper and A.~Mathas.
\newblock Morita equivalences of {Ariki-Koike} algebras.
\newblock {\em Math. Z.}, 240(3):579--610, 2002.

\bibitem[Fay07a]{Frunrem}
M.\ Fayers.
\newblock Another runner removal theorem for $v$-decomposition numbers of {Iwahori-Hecke algebras and $q$-Schur algebras}.
\newblock {\em J.\ Algebra}, 310:396--404, 2007.

\bibitem[Fay07b]{Fay07}
M.~Fayers.
\newblock Core blocks of {Ariki-Koike} algebras.
\newblock {\em J. Algebraic Combin.}, 26(1):47--81, 2007.

\bibitem[Fay07c]{Fay07conj}
M.~Fayers.
\newblock James's conjecture holds for weight four blocks of {Iwahori–Hecke} algebras.
\newblock {\em Journal of Algebra}, 317(2):593 -- 633, 2007.

\bibitem[Fay08]{Fay08conj}
M.~Fayers.
\newblock Decomposition numbers for weight three blocks of symmetric groups and {Iwahori-Hecke} algebras.
\newblock {\em Trans. Amer. Math. Soc.}, 360:1341 -- 1377, 2008.

\bibitem[Fay10]{Fay10}
M.\ Fayers.
\newblock An {LLT}-type algorithm for computing higher-level canonical bases.
\newblock {\em J.\ Algebra}, 214:2186--2198, 2010.

\bibitem[Gec92]{Gec92}
M.~Geck.
\newblock Brauer trees of {Hecke} algebras.
\newblock {\em Comm. in Algebra}, 20:2937--2973, 1992.

\bibitem[Gec98]{Gec98}
M.~Geck.
\newblock Representations of {Hecke} algebras at roots of unity.
\newblock {\em Seminaire Bourbaki}, 836, 1998.

\bibitem[GL96]{GL96}
J.~J. Graham and G.~I. Lehrer.
\newblock Cellular algebras.
\newblock {\em Inventiones mathematicae}, 123(1):1--34, 1996.

\bibitem[Jac05]{Jac05}
N.~Jacon.
\newblock An algorithm for the computation of the decomposition matrices for {Ariki-Koike} algebras.
\newblock {\em J. Algebra}, 292:100--109, 2005.

\bibitem[JK81]{JK81}
G.~James and A.~Kerber.
\newblock {\em The representation theory of the symmetric group}, volume~16 of {\em Encyclopedia of Mathematics and its Applications}.
\newblock Cambridge University Press, 1981.

\bibitem[JM02]{JM02}
G.~James and A.~Mathas.
\newblock Equating decomposition numbers for different primes.
\newblock {\em J. Algebra}, 258:599--614, 2002.

\bibitem[Kas02]{Kas02}
M.\ Kashiwara.
\newblock {\em Bases cristallines des groupes quantiques}, volume~9 of {\em Cours sp\'ecialis\'ees}.
\newblock 2002.

\bibitem[KL09]{KL09}
M.~Khovanov and A.~D. Lauda.
\newblock A diagramatic approach to categorification of quantum groups {I}.
\newblock {\em Represent. Theory}, 13:309--347, 2009.

\bibitem[LLT96]{LLT96}
A.\ Lascoux, B.\ Leclerc, and J.-Y.\ Thibon.
\newblock Hecke algebras at roots of unity and crystal bases of quantum affine algebras.
\newblock {\em Comm.\ Math.\ Phys.}, 181:205--263, 1996.

\bibitem[Mat99]{Mat99}
A.~Mathas.
\newblock {\em Iwahori-Hecke algebras and {Schur} algebras of the symmetric group}, volume~15 of {\em University Lecture Series}.
\newblock American Mathematical Society, 1999.

\bibitem[Mat03]{Mat03}
A.~Mathas.
\newblock Tilting modules for cyclotomic {Schur} algebras.
\newblock {\em J. Reine Angew. Math.}, 262:137--169, 2003.

\bibitem[Mat04]{Mat04}
A.\ Mathas.
\newblock The representation theory of the {Ariki-Koike} and cyclotomic {$q$-Schur} algebras.
\newblock {\em Representation theory of algebraic groups and quantum groups, Adv. Stud. Pure Math}, 40:261--320, 2004.

\bibitem[Ric96]{Ric96}
M.~J. Richards.
\newblock Some decomposition numbers for {Hecke} algebras of general linear groups.
\newblock {\em Mathematical Proceedings of the Cambridge Philosophical Society}, 119:383 -- 402, 1996.

\bibitem[Rou08]{Rou08}
R.~Rouquier.
\newblock 2-{Kac-Moody} algebras.
\newblock \href{https://arxiv.org/pdf/0812.5023.pdf}{arXiv:0812.5023}, 2008.

\bibitem[WKM17]{will17}
G.~Williamson, A.~Kontorovich, and P.~McNamara.
\newblock Schubert calculus and torsion explosion.
\newblock {\em J. Amer. Math. Soc.}, 30, 2017.

\bibitem[Yvo07]{Yvo07}
X.~Yvonne.
\newblock An algorithm for computing the canonical bases of higher-level $q$-deformed {Fock} spaces.
\newblock {\em J. Algebra}, 309:760--785, 2007.

\end{thebibliography}

\end{document}